\documentclass[12pt]{article}

\usepackage{amsmath,enumerate,amsbsy,amsfonts,amssymb,mathabx,amscd,graphicx,algorithm}
\usepackage[round,authoryear]{natbib}
\usepackage{authblk}
\usepackage{mathrsfs}
\usepackage{epsfig}
\newtheorem{definition}{Definition}
\newtheorem{theorem}{Theorem}

\newtheorem{lemma}{Lemma}

\newtheorem{proposition}{Proposition}
\usepackage{multirow}

\newcommand{\beps}{\boldsymbol \epsilon}

\newcommand{\bSigma}{\boldsymbol \Sigma}

\newcommand{\bsigma}{\boldsymbol \sigma}
\newcommand{\bGamma}{\boldsymbol \Gamma}

\newcommand{\bDelta}{\boldsymbol \Delta}
\newcommand{\bXi}{\boldsymbol \Xi}
\newcommand{\bpsi}{\boldsymbol \psi}
\newcommand{\bPsi}{\boldsymbol \Psi}
\newcommand{\bphi}{\boldsymbol \phi}
\newcommand{\bPhi}{\boldsymbol \Phi}

\newcommand{\bxi}{\boldsymbol \xi}
\newcommand{\btheta}{\boldsymbol \theta}
\newcommand{\bgamma}{\boldsymbol \gamma}
\newcommand{\bdelta}{\boldsymbol \delta}

\newcommand{\bfeta}{\boldsymbol \eta}
\newcommand{\bvar}{\boldsymbol \varepsilon}

\newcommand{\bvarepsilon}{\boldsymbol \varepsilon}
\newcommand{\bzeta}{\boldsymbol \zeta}

\newcommand{\ba}{{\mathbf a}}
\newcommand{\be}{{\mathbf e}}
\newcommand{\bbf}{{\mathbf f}}
\newcommand{\bx}{{\mathbf x}}
\newcommand{\by}{{\mathbf y}}

\newcommand{\bv}{{\mathbf v}}

\newcommand{\br}{{\mathbf r}}
\newcommand{\bw}{{\mathbf w}}
\newcommand{\bg}{{\mathbf g}}
\newcommand{\bs}{{\mathbf s}}
\newcommand{\bbb}{{\mathbf b}}

\newcommand{\bD}{{\bf D}}
\newcommand{\bA}{{\bf A}}
\newcommand{\bB}{{\bf B}}
\newcommand{\bC}{{\bf C}}
\newcommand{\bE}{{\bf E}}
\newcommand{\bI}{{\bf I}}

\newcommand{\bP}{{\bf P}}
\newcommand{\bS}{{\bf S}}
\newcommand{\bX}{{\bf X}}
\newcommand{\bY}{{\bf Y}}
\newcommand{\bZ}{{\bf Z}}
\newcommand{\bR}{{\bf R}}
\newcommand{\bU}{{\bf U}}
\newcommand{\bV}{{\bf V}}
\newcommand{\bQ}{{\bf Q}}
\newcommand{\bJ}{{\bf J}}

\newcommand{\bcD}{\boldsymbol{\cal D}}
\newcommand{\bcI}{\boldsymbol{\cal I}}

\newcommand{\cD}{{\cal D}}

\newcommand{\cL}{{\cal L}}
\newcommand{\cM}{{\cal M}}
\newcommand{\cK}{{\cal K}}

\newcommand{\cU}{{\cal U}}
\newcommand{\cS}{{\cal S}}

\newcommand{\mF}{{\mathscr F}}

\newcommand{\eZ}{\mathbb{Z}}
\newcommand{\eR}{\mathbb{R}}
\newcommand{\cH}{\mathbb{H}}
\newcommand{\eS}{\mathbb{S}}

\newcommand{\tA}{\text{A}}

\newcommand{\tE}{\text{E}}
\newcommand{\tF}{\text{F}}

\newcommand{\tX}{\text{X}}
\newcommand{\tR}{\text{R}}

\newcommand{\cov}{\text{Cov}}
\newcommand{\var}{\text{Var}}

\newcommand{\LS}{\text{LS}}

\newcommand{\aic}{\text{AIC}}
\newcommand{\bic}{\text{BIC}}

\newcommand{\llangle}{\langle \langle}
\newcommand{\rrangle}{\rangle \rangle}

\newcommand{\0}{{\bf 0}}
\newcommand{\T}{\text{$\scriptscriptstyle T$}}
\newtheorem{condition}{Condition}

\makeatletter
\DeclareRobustCommand\widecheck[1]{{\mathpalette\@widecheck{#1}}}
\def\@widecheck#1#2{%
	\setbox\z@\hbox{\m@th$#1#2$}%
	\setbox\tw@\hbox{\m@th$#1%
		\widehat{%
			\vrule\@width\z@\@height\ht\z@
			\vrule\@height\z@\@width\wd\z@}$}%
	\dp\tw@-\ht\z@
	\@tempdima\ht\z@ \advance\@tempdima2\ht\tw@ \divide\@tempdima\thr@@
	\setbox\tw@\hbox{%
		\raise\@tempdima\hbox{\scalebox{1}[-1]{\lower\@tempdima\box
				\tw@}}}%
	{\ooalign{\box\tw@ \cr \box\z@}}}
\makeatother
\RequirePackage[colorlinks,citecolor=blue,urlcolor=blue]{hyperref}

\begin{document}

\title{\bf \Large A General Theory for Large-Scale Curve Time \\ Series via Functional Stability Measure}
\author[1]{Shaojun Guo}
\author[2]{Xinghao Qiao}
\affil[1]{Institute of Statistics and Big Data, Renmin University of China, P.R. China}
\affil[2]{Department of Statistics, London School of Economics and Political Science, U.K.}
\setcounter{Maxaffil}{0}
\renewcommand\Affilfont{\itshape\small}
\date{\vspace{-5ex}}
\maketitle
\begin{abstract}
Modelling a large bundle of curves arises in a broad spectrum of real applications. However, existing literature relies primarily on the critical assumption of independent curve observations. In this paper, we provide a general theory for large-scale Gaussian curve time series, where the temporal and cross-sectional dependence across multiple curve observations exist and the number of functional variables, $p,$ may be large relative to the number of observations, $n.$ We propose a novel {\it functional stability measure} for multivariate stationary processes based on their spectral properties and use it to establish some useful concentration bounds on the sample covariance matrix function. These concentration bounds serve as a fundamental tool for further theoretical analysis, in particular, for deriving nonasymptotic upper bounds on the errors of the regularized estimates in high dimensional settings. As {\it functional principle component analysis} (FPCA) is one of the key techniques to handle functional data, we also investigate the concentration properties of the relevant estimated terms under a FPCA framework. To illustrate with an important application, we consider {\it vector functional autoregressive models} and develop a regularization approach to estimate autoregressive coefficient functions under the sparsity constraint. Using our derived nonasymptotic results, we investigate the theoretical properties of the regularized estimate in a ``large $p,$ small $n$" regime. The finite sample performance of the proposed method is examined through simulation studies.
\end{abstract}


\section{Introduction}
In functional data analysis, it is commonly assumed that each measured curve, treated as the unit of observation, is independently sampled from some realization of an underlying stochastic process. Curve time series, on the other hand, refers to a collection of curves observed consecutively over time, where the temporal dependence across curve observations exhibits. Existing literature mainly focuses on modelling univariate or bivariate curve time series, see, e.g., \cite{bathia2010,cho2013,panaretos2013} and \cite{hormann2015}. Recent advances in technology have made large-scale curve time series datasets become increasingly common in many applications. Examples include cumulative intraday return trajectories \cite[]{horvath2014} and functional volatility processes \cite[]{muller2011} for a large number of stocks, daily concentration curves of particulate matter and gaseous pollutants at different cities of China \cite[]{li2017}, and intraday energy consumption curves for thousands of London households (available at \url{https://data.london.gov.uk/dataset/smartmeter-energy-use-data-in-london-households}). 
These applications require understanding the relationships among a large bundle of curves based on a relatively small to moderate number of serially dependent observations.

Throughout the paper, suppose we have $n$ observed $p$-dimensional vector of random curves,
$\bX_t(\cdot)=\big(X_{t1}(\cdot), \dots, X_{tp}(\cdot)\big)^{\T}, t=1,\dots,n,$ defined on a compact interval $\cU.$ Let $\bX_t(\cdot)$ have mean zero and covariance matrix function,
$\bSigma_0(u,v)=E\big\{\bX_t(u)\bX_t(v)^{\T}\big\}, (u,v)\in \cU^2.$ Estimating $\bSigma_0$ is not only of interest by itself but also involved in subsequent analysis, such as dimension reduction and modelling of multivariate functional data. In the classical setting where $p$ is fixed, $n \rightarrow \infty$ and $\big\{\bX_t(\cdot)\big\}$ are independent and identically distributed (i.i.d.), it is well known that the
sample covariance matrix function
\begin{equation}
\label{sam.covmat.op}
\widehat \bSigma_0(u,v)=n^{-1}\sum_{t=1}^n\bX_t(u)\bX_t(v)^{\T}, \ (u,v)\in \cU^2,
\end{equation}
is a consistent estimator. 
The existing literature on $\widehat \bSigma_0$ replies on the key assumption of i.i.d. observations and its diagonalwise nonasymptotic properties have been studied in \cite{qiao2018a} and \cite{qiao2018b}. However, theoretical work on $\widehat\bSigma_0$ based on fixed or even high-dimensional curve time series, where the data exhibit not only cross-sectional but also temporal dependence, is largely untouched. It is of practical and theoretical interest to ask the following questions: {\it (i) What is the general dependence condition to be satisfied such that $\widehat\bSigma_0$ is consistent under some functional matrix norms? (ii) How does the underlying dependence structure affect the theoretical properties of $\widehat\bSigma_0?$} Attempting to answer both questions facilitates our understanding of addressing large-scale curve time series problems and forms the core of our paper.

Statistical learning in high dimensional scenarios is often impossible unless some lower-dimensional structure is imposed on the model parameter space. One large class assumes various functional sparsity patterns, where different regularized estimation procedures can be developed, see e.g., under the i.i.d. setting, functional additive regression \cite[]{fan2015} and functional graphical models \cite[]{qiao2018a}. To investigate the theoretical properties of such regularized estimates in a time series context, one need to develop suitable nonasymptotic results for $\widehat\bSigma_0.$ The existing theoretical work on $\widehat\bSigma_0$ for curve time series has focused on studying either its entrywise asymptotic properties or diagonalwise nonasymptotic properties under a special autoregressive framework \cite[]{Bbosq1}. However, these results are not sufficient to evaluate the performance of the regularized estimates under high-dimensional scaling with a general dependence structure. 

One main purpose of this paper is to establish some useful concentration bounds on $\widehat\bSigma_0$ for a large class of stationary processes,  which serve as a fundamental tool for further nonasymptotic analysis. In particular,
we focus on multivariate stationary Gaussian curve time series and introduce a novel {\it functional stability measure} for multivariate stationary processes based on their spectral properties, which provides insights into the effect of dependence on $\widehat\bSigma_0.$ With the proposed functional stability measure, we establish concentration bounds on the quadratic and bilinear forms of the operator induced from $\widehat\bSigma_0,$ and then use these results to derive some concentration bounds on $\widehat\bSigma_0$ under different functional matrix norms. Due to the infinite dimensional nature of functional data, in practice, one usually adopt a dimension reduction approach, e.g. {\it functional principal component analysis} (FPCA) based on the Karhunen-Lo\'eve expansion of a random curve, to approximate each curve by a finite representation before subsequent analysis. We also develop some concentration results for relevant estimated terms under a FPCA framework, which provide a powerful tool for the nonasymptotic analysis of the FPCA-based regularized estimates in ``large $p,$ small $n$" settings.

To illustrate the usefulness of our derived nonasymptotic results with an important application, we consider {\it vector functional autoregressive} (VFAR) models, which characterize the complex temporal and cross-sectional interrelationship among a vector of $p$ curve time series and can lead to simultaneous forecasting of multiple curves. One advantage of a VFAR model is that it accommodates dynamic linear interdependencies among multiple curve time series into a static framework within a Hilbert space. In a high dimensional VFAR model, not only $p$ is large but each curve itself is an infinite dimensional object, we assume the sparsity structure embedded on the transition matrix functions and implement FPCA on each curve, which results in fitting a higher dimensional {\it vector autoregressive} (VAR) model. Then we propose the regularized estimators for the block transition matrices, on which a block sparsity constraint is enforced via a standardized group lasso penalty \cite[]{simon2012}. Using our derived concentration results, we establish nonasymptotic upper bounds on the estimation errors and show that the regularized estimation procedure can produce consistent estimates under high-dimensional scaling.

{\bf Related literature}. High-dimensional scalar time series have been extensively studied in recent years. \cite{chen2013} considered estimation of covariance matrices and their inverses for stationary and locally stationary processes. \cite{basu2015a} established useful concentration bounds for dependent data based on their proposed stability measure for stationary Gaussian processes. \cite{chen2016} studied a Dantzig-selector type regularized estimator for linear functionals of high-dimensional linear processes.  \cite{loh2012} and \cite{wu2016} investigated the theoretical properties of the lasso estimates. The statistical inference via Gaussian approximations were considered in \cite{zhang2017}. Examples of recent developments in high dimensional VAR models include the following. \cite{kock2015} established oracle inequalities for high dimensional VAR models. \cite{basu2015a} studied the theoretical properties of $\ell_1$ regularized estimates. \cite{han2015} alternatively proposed a Dantzig-type penalization and formulated the estimation procedure as a linear program. \cite{guo2016} considered a class of VAR models with banded coefficient matrices. For examples of research on functional autoregressive models, see \cite{Bbosq1,koko2013,aue2015} and the reference therein.

{\bf Outline of the paper}. The remainder of the paper is structured as follows. In Section~\ref{sec.main}, we first introduce a functional stability measure and use it to establish concentration bounds on the quadratic/bilinear forms of the operator induced from $\widehat\bSigma_0,$ and $\widehat\bSigma_0$ itself under different functional matrix norms. We then consider the Karhunen-Lo\'eve expansion of each random curve and develop some useful concentration results under a FPCA framework. In Section~\ref{sec.vfar}, we propose a sparse VFAR model, develop a FPCA-based penalization approach, present theoretical analysis of the regularized estimates and finally illustrate the superiority of the proposed method to its competitors through simulation studies. We conclude our paper by discussing several possible future works in Section~\ref{sec.discussion}. All technical proofs are relegated to the Appendix and the  Supplementary Material. 

{\bf Notation}. We summarize here some notation to be used throughout the paper. Let $\eZ$ and $\eR$ denote the sets of integers and real numbers, respectively. We denote the indicator function by $I(\cdot).$ For a finite set $J,$ we denote its cardinality by $|J|.$ For two positive sequences $\{a_n\}$ and $\{b_n\}$, we write $a_n \succsim b_n$ if
there exists an absolute constant $c,$ such that $a_n \geq c b_n$ for all $n.$ We use $a_n \asymp b_n$ to denote $a_n \succsim b_n$ and $b_n \succsim a_n.$ We write $a_n=o(1)$ if $a_n \rightarrow 0$ as $n \to \infty.$
For a vector $\bx \in {\eR}^p,$ we denote the $\ell_q$ norm by $||\bx||_q=(\sum_{j=1}^p|x_j|^q)^{1/q}$ for $q>0.$ For two $p_1$ by $p_2$ matrices, $\bA$ and $\bB,$ we let $\langle\langle \bA,\bB\rangle\rangle=\text{trace}(\bA^{\T}\bB)$ and denote the Frobenius, operator, elementwise maximum norms of $\bB$ by $||\bB||_F=\big(\sum_{j,k} B_{jk}^2\big)^{1/2},$ $||\bB||={\sup}_{||\bx||_2\leq 1}||\bB\bx||_2$, $||\bB||_{\max}=\max_{j,k}|B_{jk}|,$ respectively. Let $L_2(\cU)$ denote a Hilbert space of square integrable functions defined on the compact set $\cU$ equipped with the inner product $\langle f,g \rangle=\int_{\cU}f(u)g(u)du$ for $f,g \in L_2(\cU)$ and the induced norm $\|\cdot\|=\langle \cdot,\cdot \rangle^{1/2}.$ For a Hilbert space $\cH \subseteq \cL_2(\cU),$ we denote the $p$-fold Cartesian product  by $\cH^p=\cH \times  \dots \times \cH$ and the tensor product by $\eS=\cH \otimes \cH.$ For two $p$-dimensional function-valued vectors, $\bbf=(f_1, \ldots, f_p)^{\T}$ and $\bg=(g_1, \dots,g_p)^{\T}$ in $\cH^p$ with $\cH=L_2(\cU),$ we denote the inner product by
$\langle\bbf,\bg\rangle=\sum_{j=1}^p\langle f_j, g_j\rangle.$ We use $\|\bbf\|_0=\sum_{j=1}^p I(\|f_j\| \neq 0)$ and $\|\bbf\|=\langle \bbf,\bbf \rangle^{1/2}$ to denote the functional versions of vectors $\ell_0$ and $\ell_2$ norms, respectively. For an integral operator $\cK: \cH \rightarrow \cH$ induced from the kernel $K \in \eS$ through 
\begin{equation}
\label{df.op}
\cK(f)(u) = \int_{\cU} K(u, v)f(v)dv = \langle K(u,\cdot),f(\cdot)\rangle,
\end{equation}
for any given $f\in \cH,$ we denote the operator norm by $\|\cK\|_{\cL}=\sup_{\|f\| \leq 1}\|\cK(f)\|.$ Without any ambiguity, we denote the Hilbert-Schmidt norm of both $K$ and $\cK$ by 
$\|K\|_{\cS} = \|\cK\|_{\cS}=\big(\int\int K(u,v)^2dudv\big)^{1/2}.$ For a $p$ by $p$ bivariate-function-valued matrix, $\bA=\big(\tA_{jk}(\cdot,\cdot)\big)_{1 \leq j,k \leq p}$ in $\eS^{p \times p},$ we define the functional versions of maximum eigenvalue, Frobenius, elementwise maximum and matrix $\ell_{\infty}$ norms by $\lambda_{\max}(\bA),$ $
\|\bA\|_{\tF} = \big(\sum_{j,k} \|\tA_{jk}\|_{\cS}^2\big)^{1/2},$ $\|\bA\|_{\max} =  \max_{j,k} \|\tA_{jk}\|_{\cS}$ and $\|\bA\|_{\infty} = \max_{j} \sum_{k} \|\tA_{jk}\|_{\cS},
$ 
respectively. For a block matrix $\bB = (\bB_{jk}) \in \eR^{p_1q \times p_2 q}$ with $\bB_{jk} \in {\eR}^{q \times q}$ for each $(j,k),$ we define its $q$-block versions of Frobenius, elementwise maximum and matrix $\ell_1$ norms by $\|\bB\|_F = \big(\sum_{j,k}\|\bB_{jk}\|_{\tF}^2\big)^{1/2},$  $\|\bB\|_{\max}^{(q)} = \max_{j,k} \|\bB_{jk}\|_{\tF}$ and 
$
\|\bB\|_{1}^{(q)} = {\max}_{k}\sum_{j} \|\bB_{jk}\|_{\tF}, 
$ respectively. 

\section{Main results}
\label{sec.main}
Suppose that $\{\bX_t(\cdot)\}_{t \in \eZ},$ defined on $\cU,$ is a sequence of $p$-dimensional centered and covariance-stationary Gaussian processes with mean zero and $p$ by $p$ autocovariance matrix function,  $\bSigma_{h}(u,v)=\cov\big\{\bX_t(u),\bX_{t+h}(v)\big\}=\big(\Sigma_{h,jk}(u,v)\big)_{1\leq j,k\leq p},t,h \in \eZ,(u,v) \in \cU^2.$ 
In particular when $h=0,$ one typically refers to $\Sigma_{0,jk}$ as marginal-covariance functions for $j=k,$ and cross-covariance functions for $j \neq k.$ 
In an analogy to (\ref{df.op}), we define an autocovariance matrix operator at lag $h,$ $\bXi_{h}:\cH^p \to \cH^p$ with the kernel $\bSigma_{h}$ such that, for any given $\bPhi \in \cH^p,$
\begin{equation}
\label{autocov.op}
    \bXi_h(\bPhi)(u): = \int_{\cU} \bSigma_{h}(u,v) \bPhi(v)dv=\left[\begin{array}{c}\langle\bsigma_{h1}(u,\cdot),\bPhi(\cdot)\rangle\\\vdots\\\langle\bsigma_{hp}(u,\cdot),\bPhi(\cdot)\rangle\end{array}\right], 
\end{equation} 
where $\bsigma_{hj}(u,\cdot)=\big(\Sigma_{h,j1}(u,\cdot), \dots, \Sigma_{h,jp}(u,\cdot)\big)^{\T}$ for $j=1, \dots, p.$ In the special case of $h=0,$ the covariance matrix function $\bSigma_0$ is symmetric and non-negative definite in the sense that $\bSigma_0(u,v)=\bSigma_0(v,u)^{\T}$ for any $(u,v) \in \cU^2$ and $\langle\bPhi, \bXi_0(\bPhi)\rangle \geq 0$ for any $\bPhi\in \cH^p.$  Note that the induced operator $\bXi_{h}$ and its kernel $\bSigma_{h}$ are in one-to-one correspondence through (\ref{autocov.op}), for notational consistency we will express terms, defined from $\bSigma_{h},$ using $\bXi_{h},$ where the context is clear. 
\subsection{Functional stability measure}
\label{sec.fsm}
Before introducing the functional stability measure, we first consider the {\it spectral density matrix operator} of  $\{\bX_t(\cdot)\}_{t \in \eZ},$ defined from the  Fourier transform of autocovariance matrix operators $\{\bXi_{h}\}_{h \in \eZ},$ which encodes the second-order dynamical properties of $\{\bX_t(\cdot)\}_{t \in \eZ}.$

\begin{definition}
\label{def.sd}
We define the spectral density matrix function of $\{\bX_t(\cdot)\}_{t \in \eZ}$ at frequency $\theta \in [-\pi,\pi]$ by
$$
	f_{\bX,\theta}(\cdot,\cdot) = \frac{1}{2 \pi} \sum_{h \in \eZ} \bSigma_{h}(\cdot,\cdot) \exp(- i h \theta),
$$
and the induced spectral density matrix operator ${\mF}_{\bX,\theta}: \cH^p \rightarrow \cH^p$ with the kernel $f_{\bX,\theta},$
such that, for any given $\bPhi \in \cH^p,$
\begin{equation}
\label{def.sd.operator}
	{\mF}_{\bX,\theta}(\bPhi) = \frac{1}{2 \pi} \sum_{h \in \eZ} \bXi_{h}(\bPhi) \exp(- i h \theta).
\end{equation}
\end{definition}

The spectral density matrix function generalizes the notion of the spectral density matrix \cite[]{basu2015a} to the functional domain, and the spectral density matrix operator can be viewed as a generalization of the spectral density operator \cite[]{panaretos2013} to the multivariate setting. The existence of $\mF_{\bX,\theta}, \theta \in [-\pi,\pi]$ defined in (\ref{def.sd.operator}) is guaranteed if $\sum_{h =0}^\infty \|\bXi_{h}\|^2_{\cL} < \infty,$ where we denote by $\|\bXi_{h}\|_{\cL}={\text{sup}}_{\|\bPhi\| \leq 1, \bPhi \in \cH^p} \|\bXi_{h}(\bPhi)\|$ the operator norm of $\bXi_{h}.$ Furthermore, if $\sum_{h =0}^\infty \|\bXi_{h}\|_{\cL} < \infty,$ then ${\mF}_{\bX,\theta}$ is uniformly bounded and continuous in $\theta$ with respect to $\|\cdot\|_{\cL},$ and the following inversion formula holds. 
\begin{equation}
\label{inv.formula}
	\bXi_{h}(\cdot) = \int_{-\pi}^\pi {\mF}_{\bX,\theta}(\cdot) \exp( i h \theta)d\theta, ~~~\text{ for all }h \in \eZ.
\end{equation} 
The inversion relationships in (\ref{def.sd.operator}) and (\ref{inv.formula}) indicate that spectral density matrix operators  and autocovariance matrix operators comprise a Fourier transform pair. Hence, we can use $\{\mF_{\bX, \theta}, \theta \in [-\pi,\pi]\}$ to study the second-order dynamics of $\{\bX_t(\cdot)\}_{t \in \eZ},$ as stated in the following condition.

\begin{condition}
\label{cond.bd.fsm}
(i) The spectral density matrix operator $\mF_{\bX,\theta}, \theta \in [-\pi,\pi]$ exists;
(ii) The functional stability measure of $\{\bX_t(\cdot)\}_{t \in \eZ},$ defined as follows, is bounded, i.e.
\begin{equation}
\label{eq.bd.fsm}
	\cM(\mF_{\bX}) =  2 \pi \cdot \underset{\theta\in [-\pi, \pi], \bPhi \in \cH_0^p}{\text{esssup}} \frac{ \big\langle \bPhi, \mF_{\bX,\theta}(\bPhi) \big\rangle}{\big\langle\bPhi, \bXi_0(\bPhi)\big\rangle} < \infty,
\end{equation}
where $\cH_0^p = \big\{\bPhi \in \cH^p: \langle\bPhi, \bXi_0(\bPhi)\rangle \in (0,\infty) \big\}.$
\end{condition}

We have several comments for the above condition. First, the functional stability measure, $\cM(\mF_{\bX}),$ is expressed as a term proportional to the essential supremum of the functional Rayleigh quotient of $\mF_{\bX,\btheta}$ relative to $\bXi_0$ over $\theta \in [-\pi,\pi].$ In particular, under the non-functional setting with $\bPhi \in \eR^p$ and $f_{\bX,\theta}, \bSigma_0 \in \eR^{p \times p},$ (\ref{eq.bd.fsm}) reduces to
$$
2 \pi \cdot \underset{\theta\in [-\pi, \pi], \bPhi \neq {\bf 0}}{\text{ess}\sup} \frac{\bPhi^{\T} f_{\bX,\theta}\bPhi}{\bPhi^{\T}\bSigma_0 \bPhi} < \infty,
$$
which is equivalent to the upper bound condition for the stability measure, $\widetilde\cM(f_{\bX}),$ introduced by \cite{basu2015a}, i.e.
$$\widetilde\cM(f_{\bX})=\underset{\theta\in [-\pi, \pi], \bPhi \neq {\bf 0}}{\text{ess}\sup} \frac{\bPhi^{\T} f_{\bX,\theta}\bPhi}{\bPhi^{\T}\bPhi} < \infty.$$
Second, we denote the functional analog of $\widetilde\cM(f_{\bX})<\infty$ by
\begin{equation}
\label{eq.bd.basu}
\widetilde\cM(\mF_{\bX})=\underset{\theta\in [-\pi, \pi], \bPhi \neq {\bf 0}}{\text{ess}\sup}\frac{ \big\langle \bPhi, \mF_{\bX,\theta}(\bPhi) \big\rangle}{\big\langle\bPhi, \bcI_p(\bPhi)\big\rangle} < \infty,
\end{equation}
where $\bcI_p$ is the identity matrix operator induced from the kernel, $\bI_p=\big(C_{jk}(\cdot,\cdot)\big)_{1 \leq j,k \leq p}\in {\eS}^{p \times p}$ with $C_{jk}(u,v)=I(j=k)I(u=v)$ for $(u,v) \in \cU^2.$ By contrast, our proposed functional stability measure, $\cM(\mF_{\bX}),$ makes more sense than $\widetilde\cM(\mF_{\bX}),$ since it can account for the effect of eigenvalues of $\bXi_0$ at different magnitude levels.  Specially, if $X_{t1}(\cdot), \dots, X_{tp}(\cdot)$ are finite dimensional objects, the upper bound conditions in (\ref{eq.bd.fsm}) and (\ref{eq.bd.basu}) would be equivalent. Third, it is clear that, unlike $\widetilde\cM(f_{\bX})$ or $\widetilde\cM(\mF_{\bX}),$ $\cM(\mF_{\bX})$ is a scale-free stability measure. In the special case of no temporal dependence, we have $\cM(\mF_\bX)=1.$ Fourth, since the autocovariance matrix function characterizes a multivariate Gaussian processes, it can be used to quantify the temporal and cross-sectional dependence for this class of models. In particular, the spectral density matrix function provides insights into the stability of the process. In our analysis of large-scale curve time series, we will use $\cM(\mF_{\bX})$ as a stability measure of the process of $\{\bX_t(\cdot)\}_{t \in \eZ}.$ Larger values of $\cM(\mF_{\bX})$ would correspond to a less stable process.

Condition~\ref{cond.bd.fsm} is satisfied by a large class of functional linear processes including, e.g. VFAR models and functional factor models, where the underlying spectral density matrix function/operator can be expressed. See Section~\ref{sec.vfar} for details on VFAR models. To illustrate using an example, we consider a VFAR model of order $1$, denoted as VFAR(1), as follows
\begin{equation}
\label{vfar1.ex}
\bX_t(u) = \int_{\cU} \bA(u,v) \bX_{t-1}(v)dv + \bvar_t(u), \ (u,v)\in\cU^2. 
\end{equation}
In the special case of a symmetric $\bA$, i.e. $\bA(u,v)=\bA(v,u)^{\T},$ 
equation (\ref{vfar1.ex}) has a stationary solution if and only if $||\bA||_{\cL} <1.$ See Theorem~3.5 of \cite{Bbosq1} for the discussion of this condition when $p=1.$ However, this restrictive condition is violated by many stable VFAR(1) models with nonsymmetric $\bA.$ Moreover, it does not generalize beyond VFAR(1) models. See \cite{basu2015a} for the remark on the restrictiveness of an analogous condition in the non-functional setting. 

\begin{figure*}[t]
 	\centering
   	\includegraphics[width=7.2cm,height=7.2cm]{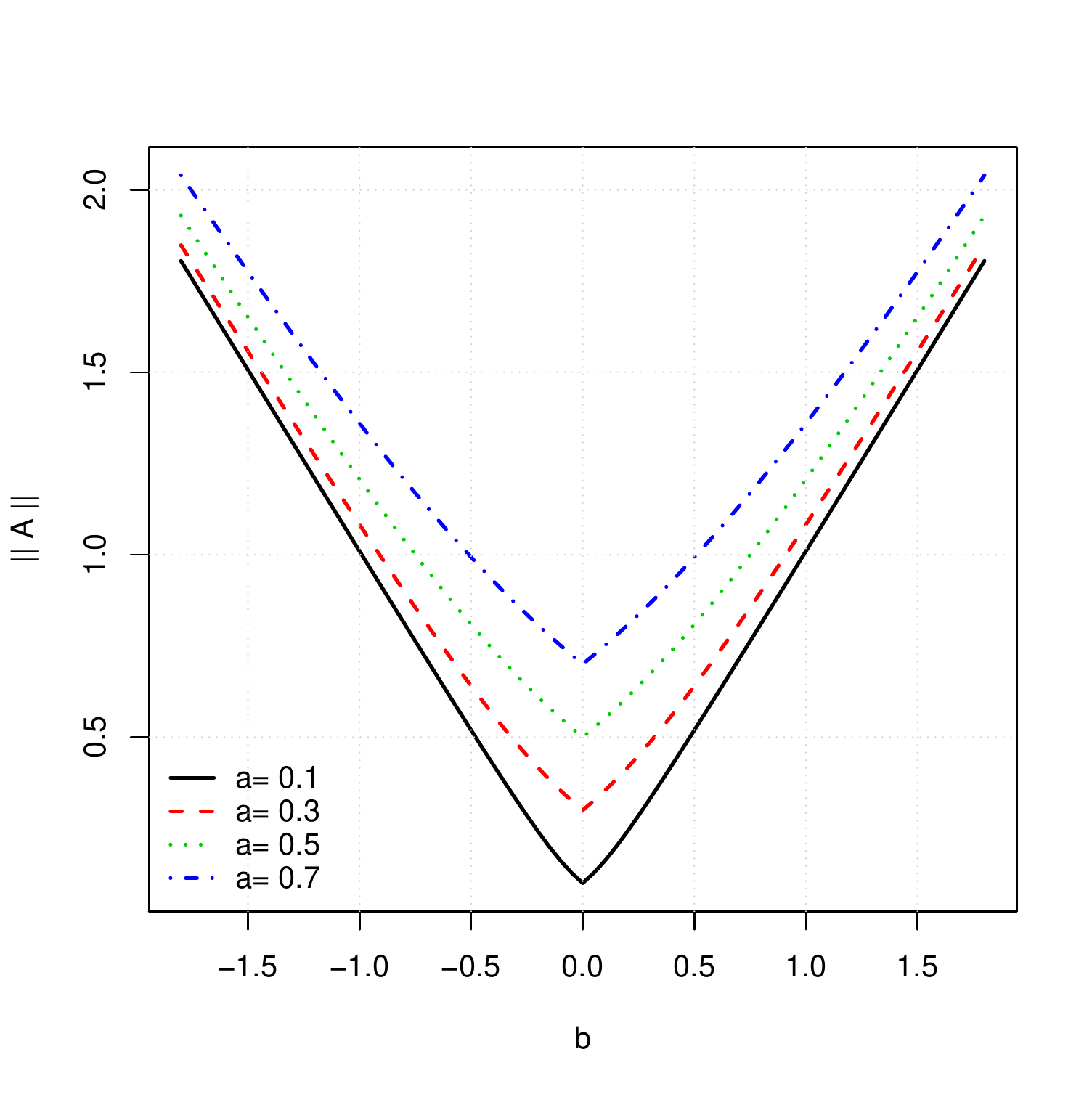}\hspace{.5cm}\includegraphics[width=7.2cm,height=7.2cm]{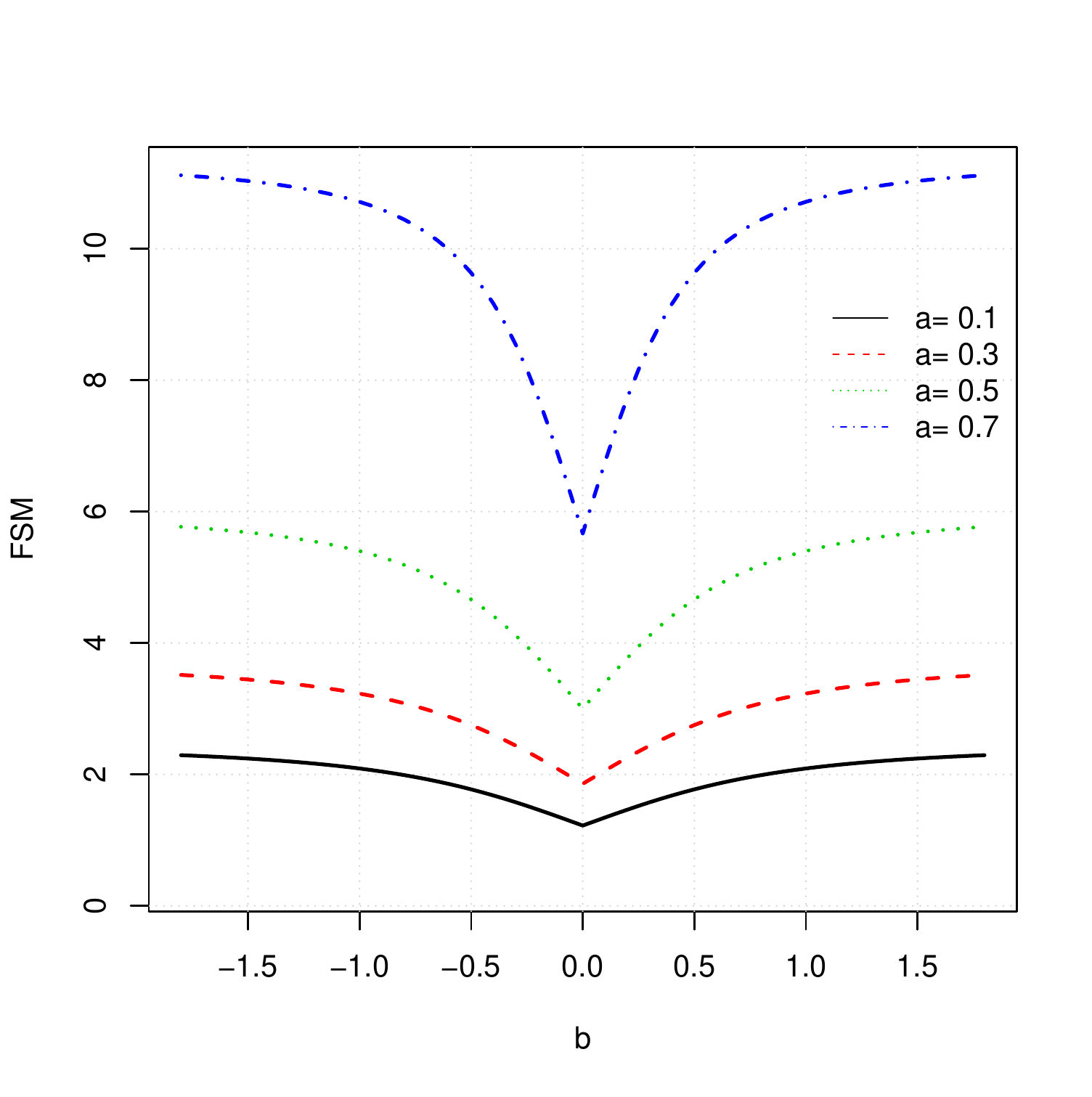}
 	\vspace{-0.8cm}
 	\caption{\label{ill.example}{\it \small The illustrative VFAR(1) model. Left: $||\bA||_{\cL}$ as a function of $a$ and $b$, plotted against $b$ for different $a.$ Right: $\cM(\mF_{\bX})$ as a function of $a$ and $b$, plotted against $b$ for different $a.$}}
 \end{figure*}

We consider an illustrative example with
\begin{equation}
\label{vfar1.ex.A}
\bA(u,v)=\left[\begin{array}{cc}a\psi_1(u)\psi_1 (v)& b\psi_1(u)\psi_2(v)\\0 & a\psi_2(u)\psi_2(v)\end{array}\right],
\bX_t(u)=\left[\begin{array}{c}x_{t1}\psi_1(u)\\x_{t2}\psi_2(u)\end{array}\right] 
\mbox{ and }\bvarepsilon_t(u)=\left[\begin{array}{c}e_{t1}\psi_1(u)\\e_{t2}\psi_2(u)\end{array}\right],
\end{equation}
where
$(e_{t1},e_{t2})^{\T} \overset{\text{i.i.d.}}{\sim}N({\bf 0}, \bI_2)$
and $\|\psi_j\|=1$ for $j=1,2.$ Section~\ref{ap.ill.ex} of the Supplementary Material provides details to calculate $\rho(\bA)$ (spectral radius of $\bA$), $||\bA||_{\cL}$ and $\cM(\mF_{\bX})$ for the example we consider. In particular, $\rho(\bA)=|a|<1$ corresponds to a stationary solution to equation~(\ref{vfar1.ex}). Figure~\ref{ill.example} 
visualizes $||\bA||_{\cL}$ and $\cM(\mF_{\bX})$ for various values of $a \in (0,1).$ We observe a few apparent patterns. First, increasing $a$ results in a value for larger $||\bA||_{\cL}.$ As $|b|$ grows large enough, the condition of $||\bA||_{\cL}<1$ will be violated, but equation~(\ref{vfar1.ex}) still have a stationary solution. Second, processes with stronger temporal dependence, i.e. with larger values of $a$ or $|b|,$ have larger values of $\cM(\mF_{\bX})$ and will be considered less stable. In our high-dimensional VFAR(1) modelling, it makes more sense to use $\cM(\mF_{\bX})$ rather than $||\bA||_{\cL}$ as a measure of stability of the process. 

\begin{definition}
\label{def.sub.fsm}
For all $k$-dimensional subprocesses of $\{\bX_{t}(\cdot)\}_{t\in\eZ},$ i.e. $\big\{\big(X_{tj}(\cdot)\big):j \in J\big\}_{t \in \eZ},$ for $J \subseteq\{1,\dots,p\}$ and $|J| \leq k,$ we define the corresponding functional stability measure by
\begin{equation}
\label{eq.sub.fsm}
\cM_k(\mF_{\bX}) = 2\pi \cdot \underset{\theta\in [-\pi, \pi],\|\bPhi\|_0 \le k,\bPhi \in \cH_0^p}{\text{ess} \sup} \frac{\langle \bPhi, \mF_{\bX,\theta}(\bPhi)\rangle}{\langle \bPhi, \bXi_0(\bPhi)\rangle}, \ k=1,\dots,p.
\end{equation}
\end{definition}

It is obvious from definitions in  Condition~\ref{cond.bd.fsm} and (\ref{eq.sub.fsm}) that
$$
\cM_1(\mF_{\bX}) \le \cM_2(\mF_{\bX}) \le \dots \le \cM_p(\mF_{\bX}) = \cM(\mF_{\bX}) <\infty,
$$
which will be used in the subsequent analysis.

\subsection{Concentration bounds on $\widehat \bXi_0$ and $\widehat \bSigma_0$}
\label{con.bd.Sigma}
Based on $n$ observed $p$-dimensional realizations, $\bX_1(\cdot),\dots, \bX_n(\cdot),$ generated from a stationary process satisfying Condition~\ref{cond.bd.fsm}, we construct an empirical estimator for $\bSigma_{h}$ by
\begin{equation}
    \label{est.autocov.fn}
    \widehat \bSigma_{h}(u,v)=\frac{1}{n-h}\sum_{t=1}^{n-h} \bX_t(u)\bX_{t+h}(v)^{\T}, (u,v)\in \cU^2.
\end{equation}
The estimator $\widehat\bXi_{h}$ for $\bXi_{h}$ can be defined by replacing $\bSigma_{h}(\cdot,\cdot)$ with $\widehat\bSigma_{h}(\cdot,\cdot).$

\begin{condition}
	\label{cond.cov.func}
	(i) The marginal-covariance functions, $\Sigma_{jj}(u,v)$'s, are continuous on $\cU^2$ and bounded uniformly on $j\in\{1,\ldots,p\}$ and $(u,v)\in \cU^2;$
	(ii) $\lambda_0 = {\max}_{1\leq j \leq p}\int_{\cU} \Sigma_{jj}(u,u)du < \infty$.
\end{condition}

The following theorem provides the concentration bounds on $\widehat \bXi_0$ under the quadratic and bilinear forms. These concentration results serve as a starting point to establish some nonasymptotic upper bounds on the error $\widehat \bSigma_0-\bSigma_0$ and the convergence rate of the VFAR estimate in Section~\ref{sec.vfar}.

\begin{theorem}
	\label{thm_con_op}
	Suppose that Conditions~\ref{cond.bd.fsm} and \ref{cond.cov.func} hold. Then for any given vectors $\bPhi_1,$ 
	$\bPhi_2 \in \cH_0^p$ with $\max(\|\bPhi_1\|_0, \|\bPhi_2\|_0) \le k,$ $k=1,\dots,p,$ there exists some universal constant $c > 0$ such that for any $\eta > 0,$
	\begin{equation}
	\label{thm_Phi_1}
	P\left\{\left|\frac{\big\langle \bPhi_1, (\widehat{\bXi}_0 - \bXi_0)(\bPhi_1)\big\rangle} {\big\langle \bPhi_1, \bXi_0(\bPhi_1)\big\rangle}\right| > \cM_k(\mF_{\bX})\eta\right\} \le 2\exp\Big\{- c n \min(\eta^2, \eta)\Big\},
	\end{equation}
	and
	\begin{equation}
	\label{thm_Phi_12}
	P\left\{\left|\frac{\big\langle \bPhi_1, (\widehat{\bXi}_0 - \bXi_0)(\bPhi_2)\big\rangle} {\big\langle \bPhi_1, \bXi_0(\bPhi_1)\big\rangle +
	\big\langle \bPhi_2, \bXi_0(\bPhi_2)\big\rangle}\right| >   \cM_k(\mF_{\bX}) \eta\right\} \le 4\exp\Big\{- c n \min(\eta^2, \eta)\Big\}.
	\end{equation}	
\end{theorem}

\vspace{1ex}

The concentration inequalities in (\ref{thm_Phi_1}) and (\ref{thm_Phi_12}) suggest that the temporal dependence may affect the concentration properties via $\cM_k(\mF_{\bX})$ and in turn the nonasymptotic results in two different ways, depending on which term in the tail bounds is dominant. We next apply the results in Theorem~\ref{thm_con_op} to establish some concentration inequalities and nonasymptotic upper bounds on $\widehat \bSigma_0-\bSigma_0$ 
under different functional matrix norms.

\begin{theorem}
	\label{thm_max_bound}
	Suppose that Conditions~\ref{cond.bd.fsm} and \ref{cond.cov.func} hold. Then there exists some universal constant $\tilde c > 0$ such that for any $\eta > 0$ and each $j,k =1,\ldots,p,$
	\begin{equation}
	\label{thm_Sigma_compt}
	P\left\{\big\|\widehat{\Sigma}_{0,jk} - \Sigma_{0,jk}\big\|_{\cS} >   2 \cM_1(f_{\bX})\lambda_0 \eta\right\} \le 4 \exp\Big\{- \tilde c n \min(\eta^2, \eta)\Big\}, 
	\end{equation}
	and 
	\begin{equation}
	\label{thm_Sigma_max}
	P\left\{\big\|\widehat{\bSigma}_0 - \bSigma_0\big\|_{\max} >   2 \cM_1(\mF_{\bX})\lambda_0 \eta\right\} \le 4p^2 \exp\Big\{- \tilde c n \min(\eta^2, \eta)\Big\}.
	\end{equation}
	In particular, if the sample size $n$ satisfies the bound
	$n \geq \rho^2\log p,$ where $\rho$ is some constant with $\rho > \sqrt{2} \tilde c^{-1/2},$ then with probability greater than $1- 4p^{2- \tilde c\rho^2}$, the estimate $\widehat{\bSigma}_0$ satisfies the bound
	\begin{equation}
	\label{bd_Sigma_max}
	\big\|\widehat{\bSigma}_0 - \bSigma_0\big\|_{\max} \le 2 \cM_1(\mF_{\bX})\lambda_0 \rho \sqrt{\frac{\log p}{n}}.
	\end{equation}
\end{theorem}

\begin{theorem}
	\label{thm_F_bound}
	Suppose that Conditions~\ref{cond.bd.fsm} and \ref{cond.cov.func} hold. Then 
	 there exists some universal constant $\tilde c > 0$ such that
	for any $\eta > 0$
	\begin{equation}
	\label{thm_Sigma_F}
	P\left\{\big\|\widehat{\bSigma}_0 - \bSigma_0\big\|_{F} >   2 \cM_1(\mF_{\bX}) \lambda_0 \eta \right\} \le \frac{p^2}{\eta^2 n} \big (16 \tilde c^{-1} + 128 \tilde c^{-2}n^{-1} \big).
	\end{equation}
	In particular, if the sample size $n$ satisfies the bound $n>128(\tilde\rho^2\tilde c^{2}-16\tilde c)^{-1},$  where $\tilde\rho$ is some positive constant with $\tilde\rho>4\tilde c^{-1/2},$ then with probability greater than $1 - \tilde\rho^{-2} \big (16 \tilde c^{-1} + 128 \tilde c^{-2}n^{-1} \big),$ the estimate $\widehat\bSigma_0$ satisfies the bound 
	\begin{equation}
	\label{bd_Sigma_F}
	\big\|\widehat{\bSigma}_0 - \bSigma_0\big\|_{F} \le  2 \cM_1(\mF_{\bX}) \lambda_0 \tilde\rho\sqrt{\frac{p^2}{n}}.
	\end{equation}
\end{theorem}

The concentration properties of $\widehat \Sigma_{0,jk}$ under the independence and ``$j=k$" setting were studied in \cite{qiao2018a}. By contrast, when $\bX_1(\cdot), \dots, \bX_n(\cdot)$ are serially dependent, (\ref{thm_Sigma_compt}) provides concentration bounds on $\widehat\Sigma_{0,jk}$ with free choices of $(j,k)$ and more explicit control of the constants. In the bounds established in Theorems~\ref{thm_max_bound} and \ref{thm_F_bound}, the effect of dependence is captured by $\cM_1(\mF_{\bX})$ with larger values yielding a slower convergence rate. When the processes of $\{\bX_t(\cdot)\}_{t \in \eZ}$ are stable and $\cM_1(\mF_{\bX})$ remains constant with respect to $(n,p),$ we obtain the error rates  under functional elementwise maximum and Frobenius norms as $\big\|\widehat{\bSigma}_0 - \bSigma_0\big\|_{\max}=O_P\big\{(\log p/n)^{1/2}\big\}$ and $\big\|\widehat{\bSigma}_0 - \bSigma_0\big\|_F=O_P\big\{(p^2/n)^{1/2}\big\},$ respectively. These convergence rates are of the same order as those of the sample covariance matrix for i.i.d. scalar observations \cite[]{bickel2008}.
See also \cite{basu2015a} for a similar elemenetwise convergence result in a scalar time series context.

In the following proposition, we present similar concentration bounds on the autocovariance matrix function $\widehat{\bXi}_h$ at lag $h \neq 0.$
\begin{proposition}
	\label{thm_sigma_k}
	Suppose that Conditions~\ref{cond.bd.fsm} and \ref{cond.cov.func} hold. Then for any given vectors $\bPhi_1,$ 
	$\bPhi_2 \in \cH_0^p$ with $\max(\|\bPhi_1\|_0, \|\bPhi_2\|_0) \le k,$ $k=1,\dots,p,$ there exists some universal constant $c > 0$ such that for any $\eta > 0$ and $h \neq 0,$
	\begin{equation}
	\label{thm_sig_k1}
	P\left\{\left|\frac{\big\langle \bPhi_1, (\widehat{\bXi}_h - \bXi_h)(\bPhi_1)\big\rangle} {\big\langle \bPhi_1, \bXi_0(\bPhi_1)\big\rangle}\right| > 2 \cM_k(\mF_{\bX})\eta\right\} \le 4\exp\Big\{- c n \min(\eta^2, \eta)\Big\},
	\end{equation}
	and
	\begin{equation}
	\label{thm_sig_k2}
	P\left\{\left|\frac{\big\langle \bPhi_1, (\widehat{\bXi}_h - \bXi_h)(\bPhi_2)\big\rangle}
    	{\big\langle \bPhi_1, \bXi_0(\bPhi_1)\big\rangle +  
    	\big\langle \bPhi_2, \bXi_0(\bPhi_2)\big\rangle}\right| >    2\cM_k(\mF_{\bX}) \eta\right\} \le 8\exp\Big\{- c n \min(\eta^2, \eta)\Big\}.
	\end{equation}	
\end{proposition}

Proposition~\ref{thm_sigma_k} provides the concentration bounds on $\widehat \bXi_h$ under the quadratic and bilinear forms and can be used to derive the exponential convergence rates of $\widehat{\bSigma}_h$ under functional elementwise maximum and Frobenius norms as
$\big\|\widehat{\bSigma}_h - \bSigma_h\big\|_{\max}=O_P\left\{\cM_1(\mF_{\bX})(\log p/n)^{1/2}\right\}$ and $ \big\|\widehat{\bSigma}_h - \bSigma_h\big\|_F=O_P\left\{\cM_1(\mF_{\bX})(p^2/n)^{1/2}\right\},$ respectively. In particular, the diagonalwise concentration results for $\widehat\bSigma_h$ are useful to address error-contaminated curve time series problems under an autocovariance framework \cite[]{bathia2010,qiao2018c}.

\subsection{Concentration bounds under a FPCA framework}
\label{sec.bd.fpca}
For each $j=1,\dots, p,$ we assume that $X_{tj}(\cdot)$ admits the Karhunen-Lo\'eve expansion, i.e. 
$
    X_{tj}(\cdot)=\sum_{l=1}^{\infty}\xi_{tjl}\phi_{jl}(\cdot),
$
which forms the foundation of FPCA.
The coefficients $\xi_{tjl}=\langle X_{tj}, \phi_{jl}\rangle,$ $l\geq 1,$ namely {\it functional principal component} (FPC) scores, correspond to a sequence of random variables with $E(\xi_{tjl})=0,$ $\text{Var}(\xi_{tjl})=\lambda_{jl}$ and $\cov(\xi_{tjl},\xi_{tjl'})=0$ if $l \neq l'.$ The eigen-pairs $\{(\lambda_{jl},\phi_{jl})\}_{l \geq 1}$ satisfy the eigen-decomposition $ \langle\Sigma_{0,jj}(u,\cdot), \phi_{jl}(\cdot) \rangle=\lambda_{jl}\phi_{jl}(u)$ with $\lambda_{j1} \geq \lambda_{j2}\geq \cdots.$ We say that $X_{tj}(\cdot)$ is $d_j$-dimensional if $\lambda_{jd_{j}} \neq 0$ and $\lambda_{j(d_j+1)}=0$ for some integer $1\leq d_j<\infty.$ If $d_j=\infty,$ all the eigenvalues are nonzero and $X_{tj}(\cdot)$ is a truly infinite dimensional functional object.

To implement FPCA based on realizations, $X_{1j}(\cdot), \dots, X_{nj}(\cdot),$ we first compute the sample estimator for $\Sigma_{0,jj}$ by $\widehat\Sigma_{0,jj}(u,v)=n^{-1}\sum_{t=1}^nX_{tj}(u)X_{tj}(v).$ Performing the eigen-decomposition on $\widehat\Sigma_{0,jj},$ i.e. $\langle \widehat\Sigma_{0,jj}(u,\cdot), \widehat\phi_{jl}(\cdot)\rangle=\widehat \lambda_{jl}\widehat\phi_{jl}(u), l\geq 1,$ leads to the estimated eigen-pairs $(\widehat\lambda_{jl},\widehat\phi_{jl})$ and estimated FPC scores $\widehat\xi_{tjl}=\langle X_{tj}, \widehat\phi_{jl}\rangle$ for $l \geq 1.$

\subsubsection{Eigenvalues and eigenfunctions}
\label{sec.bd.eigen}
To present the relevant concentration results under a FPCA framework, we impose the following regularity condition.

\begin{condition}
\label{cond_eigen}
For each $j=1,\dots,p,$ all the nonzero eigenvalues of $\Sigma_{0,jj}$ are different, i.e. $\lambda_{j1} >  \lambda_{j2} > \cdots > \lambda_{jd_j},$ and there exist some positive constants $c_0$ and $\alpha>1$ such that $ \lambda_{jl}-\lambda_{j(l+1)} \geq c_0 l^{-\alpha-1}$ for $l =1,\dots,d_j.$
\end{condition}

Condition~\ref{cond_eigen} is standard in functional data analysis literature, see e.g. \cite{hall2007} and \cite{qiao2018c}. The parameter $\alpha$ controls the lower bounds for spacings between adjacent eigenvalues with larger values of $\alpha$ allowing tighter eigen-gaps. This condition also implies that $\lambda_{jl} \geq c_0 \alpha^{-1} l^{-\alpha}$ as $\lambda_{jl}=\sum_{k=l}^{\infty}\{\lambda_{jk}-\lambda_{j(k+1)}\}\geq c_0\sum_{k=l}^\infty k^{-\alpha-1}.$ 

\begin{theorem} 
\label{thm_lambda}
Suppose that Conditions~\ref{cond.bd.fsm}--\ref{cond_eigen} hold. Then there exists some universal constant $c_1>0$ such that for each $j=1,\dots,p, l=1, \dots,d_j$ and any $\eta >0$,
\begin{equation}
\label{thm_con_lambda}
P\left\{\left |\frac{\widehat \lambda_{jl} - \lambda_{jl}}{\lambda_{jl}}\right|  >  \cM_1(\mF_{\bX}) \eta  + \rho_1 l^{2\alpha + 1} \cM_{1}^2(\mF_{\bX}) \eta^2\right\} \le 4 \exp\Big\{- c_1 n \min(\eta^2,\eta)\Big\},
\end{equation}
where $\rho_1 = 16\sqrt{2} c_0^{-2} \alpha \lambda_0^2.$ In particular, let $M$ be a positive integer possibly depending on $(n,p)$ with $\min_{1\leq j \leq p}\lambda_{jM}>0.$ If the sample size  $n$ satisfies the bound $n\ge \tilde\rho_2^2\log(pM)\max\big\{(\tilde\rho_1-1)^{-2}\rho_1^2M^{4\alpha+2}\cM_{1}^2(\mF_{\bX}),1\big\},$ where
$\tilde\rho_1$ and $\tilde\rho_2$ are some constants with $\tilde\rho_1>1$ and $\tilde\rho_2>c_1^{-1/2},$ then with probability greater than $1- 4(pM)^{1- c_1\tilde\rho_2^2}$, the estimates $\big\{\widehat\lambda_{jl}:j=1,\dots,p, l=1,\dots,M\big\}$ satisfy the bound
\begin{equation}
\label{bd_lambda_max}
\max_{1 \le j\le p, 1 \le l \le M}\left |\frac{\widehat \lambda_{jl} - \lambda_{jl}}{\lambda_{jl}}\right| \le \tilde\rho_1 \tilde\rho_2\cM_{1}(\mF_{\bX})\sqrt{\frac{\log (pM)}{n}}.
\end{equation}
\end{theorem}

The concentration inequalities for the relative errors of $\{\widehat\lambda_{jl}\}$ are presented in (\ref{thm_con_lambda}) with the effect of temporal dependence being encoded by $\cM_{1}(\mF_{\bX}).$ We provide three remarks for the above theorem. First, when $\bX_1(\cdot), \dots, \bX_n(\cdot)$ are independent, the concentration bounds on the absolute errors of $\{\widehat\lambda_{jl}\}$ can be found in \cite{qiao2018a}. By contrast, Theorem~\ref{thm_lambda} does not require the condition for upper bounds on eigenvalues. Moreover, when $\cM_{1}(\mF_{\bX})$ remains constant with respect to $(n,p),$ (\ref{thm_con_lambda}) ensures a sharper bound on $\widehat\lambda_{jl}-\lambda_{jl}$ as long as $\lambda_{jl}$ converges to zero as $l$ grows with $n.$ Second, when $X_{tj}(\cdot)$ is $d_j$-dimensional for $j=1,\dots,p$ with $d=\max_{j}d_j<\infty,$ it is easy to show that a convergence rate similar to (\ref{bd_lambda_max}) can be achieved, i.e. $\max_{1\leq j\leq p, 1\leq l\leq d_j}\big|(\widehat\lambda_{jl}-\lambda_{jl})/\lambda_{jl}\big|=O_P\left[\cM_1(\mF_{\bX})\{\log (pd)/n\}^{1/2}\right].$ Third, in the special case of $\lambda_{jl}=0$ for $l>d_j,$  comparing with the results in (\ref{thm_con_lambda}), we can obtain faster exponential decay rate for $\widehat\lambda_{jl}$ due to the property of first order degeneracy. See also similar arguments on the fast convergence for zero eigenvalues in \cite{bathia2010} and \cite{lam2012}.

\begin{theorem}
\label{thm_principal}
Suppose that Conditions~\ref{cond.bd.fsm}--\ref{cond_eigen} hold. Then there exists some universal constant $c_2>0$ such that for each $j=1,\dots,p, l=1, \dots,d_j,$ any given function $g \in \cH$ and any $\eta > 0,$
	\begin{equation}
	\label{thm_con_phi}
	P\left\{\Big\|\widehat \phi_{jl} - \phi_{jl}\Big\| > 4 \sqrt{2}  \cM_1(\mF_{\bX}) \lambda_0 c_0^{-1} l^{\alpha + 1}  \eta\right\} 
	\le  4 \exp\Big\{ - c_2 n \min(\eta^2, \eta)\Big\},
	\end{equation}
	and 
	\begin{equation}
	\label{thm_con_phi_g}
	\begin{split}
	&P\left\{\left|\big \langle \widehat \phi_{jl} - \phi_{jl}, g\big \rangle \right| \ge
	\rho_2\|g^{-jl}\|_{\lambda}  \cM_1(\mF_{\bX}) \lambda_{jl}^{1/2}l^{\alpha + 1}  \eta
	+ \rho_3 \|g\|\cM_1^2(\mF_{\bX})  l^{2(\alpha+1)}\eta^2 \right\} \\
	&\le 8 \exp\Big\{ - c_2 n\min(\eta^2, \eta)\Big\} + 4\exp\Big\{ - c_2  \cM_1^{-2}(\mF_{\bX}) n l^{-2(\alpha + 1)}\Big\},
	\end{split}
	\end{equation}
where $g(\cdot)=\sum_{m=1}^{\infty}g_{jm}\phi_{jm}(\cdot),$ $\|g^{-jl}\|_{\lambda} = \big({\sum}_{m: m\neq l}\lambda_{jm} g_{jm}^2 \big)^{1/2},$
$
\rho_{2} = 2 c_0^{-1}  \lambda_0 
$
and 
$\rho_{3} = 4(6+2\sqrt{2}) c_0^{-2} \lambda_0^2$ with $c_0 \leq 4\cM_1(\mF_{\bX})\lambda_0l^{\alpha+1}.$
In particular, let $M$ be a positive integer possibly depending on $(n,p)$ with $\min_{1\leq j \leq p}\lambda_{jM}>0.$ If the sample size $n$ satisfies the lower bound
$$n \geq \tilde \rho_4^2 \log(pM)\max\Big[\big\{\tilde\rho_3-\rho_2{\max}_{1\leq j\leq p,1\leq l\leq M}(||g^{-jl}||_{\lambda}\lambda_{j1})\big\}^{-2}\rho_3^2||g||^2M^{2\alpha+2}\cM_{1}^2(\mF_{\bX}),1\Big],$$
where $\tilde\rho_3$ and $\tilde\rho_4$ are some constants with $\tilde\rho_3>\rho_2{\max}_{j,l}(||g^{-jl}||_{\lambda}\lambda_{j1})$ and $\tilde\rho_4>c_2^{-1/2},$ then with probability greater than $1-8(pM)^{1-c_2\tilde\rho_4^2},$ 
the estimates $\big\{\widehat\phi_{jl}(\cdot):j=1,\dots,p, l=1,\dots, M\big\}$ satisfy the bounds
\begin{equation}
\label{bd_phi_max}
\max_{1 \le j\le p, 1 \le l \le M}\left\| \widehat \phi_{jl} - \phi_{jl}\right\| \le 4\sqrt{2}\lambda_0c_0^{-1} \tilde\rho_4\cM_{1}(\mF_{\bX}) M^{\alpha+1}\sqrt{\frac{\log (pM)}{n}},
\end{equation}
and
\begin{equation}
\label{bd_phig_max}
\max_{1 \le j\le p, 1 \le l \le M}\left|\big \langle \widehat \phi_{jl} - \phi_{jl}, g\big \rangle \right| \le \tilde\rho_3 \tilde\rho_4\cM_{1}(\mF_{\bX}) M^{\alpha+1}\sqrt{\frac{\log (pM)}{n}}.
\end{equation}
\end{theorem}

We give two comments for Theorem~\ref{thm_principal}. First, the nonasymptotic results for $\{\widehat\phi_{jl}(\cdot)\}$ in (\ref{thm_con_phi})--(\ref{bd_phig_max})
are commonly affected by the temporal
dependence through $\cM_{1}(\mF_{\bX}).$ For independent curve observations, the concentration bounds on $\{\widehat\phi_{jl}(\cdot)\}$ under the $L_2$ norm were developed in \cite{qiao2018a}. By contrast, our concentration results hold for a larger class of stationary processes, while the independence setting can be viewed as a special situation with $\cM_{1}(\mF_{\bX})=1.$ Second, in the special case of $g \in \{\phi_{jm}\}_{m \geq 1},$ the result in (\ref{thm_con_phi_g}) corresponds to the optimal exponential decay rate in the following sense. (i) When $g=\phi_{jm}$ for $m \neq l,$ (\ref{thm_con_phi_g}) implies that $\big|\lambda_{jm}^{-1/2}\lambda_{jl}^{-1/2}\langle\widehat\phi_{jl}-\phi_{jl}, \phi_{jm}\rangle\big|=O_P\big(l^{\alpha+1}n^{-1/2}\big);$ (ii) When $g=\phi_{jl},$ we have $\|g^{-jl}\|_{\lambda}=0,$ thus obtaining a faster convergence rate of $\big|\langle\widehat\phi_{jl}-\phi_{jl}, \phi_{jl}\rangle\big|=O_P\big(l^{2(\alpha+1)}n^{-1}\big).$ Provided the fact that $\langle\widehat \phi_{jl} - \phi_{jl}, \phi_{jl}\rangle = - 2^{-1} \|\widehat \phi_{jl} - \phi_{jl}\|^2,$ this convergence rate is consistent to that of $\|\widehat \phi_{jl} - \phi_{jl}\|=O_P\big(l^{\alpha+1}n^{-1/2}\big),$ implied by (\ref{thm_con_phi}). 

\subsubsection{Covariance between FPC scores} 
\label{sec.bd.score}

For each $j,k=1,\dots,p,$ $l=1,\dots d_j$ and $m=1,\dots, d_k,$ denote the covariance between FPC scores $\xi_{tjl}$ and $\xi_{tkm}$ by $\sigma_{jklm}$ and its sample estimate by $\widehat \sigma_{jklm}=n^{-1}\sum_{t=1}^n\widehat\xi_{tjl}\widehat\xi_{tkm}$ We next present the concentration bounds on $\{\widehat\sigma_{jklm}\}.$
\begin{theorem}
\label{thm_cov}
	Suppose that Conditions~\ref{cond.bd.fsm}--\ref{cond_eigen} hold. Then there exists some universal constant $c_3>0$ such that  
	\begin{description}
		\item[] (i) for each $j = 1,\ldots,p,$ $l=1,\ldots,d_j$ and any $\eta >0,$
		\begin{equation}
		\label{con_bd_score1}
		P\left\{\left |\frac{\widehat \sigma_{jjll} - \sigma_{jjll}}{\lambda_{jl}}\right|  >  \cM_1(\mF_{\bX}) \eta  + \rho_1 \cM_{1}^2(\mF_{\bX}) l^{2\alpha + 1}\eta^2\right\} \le 4 \exp\Big\{- c_3n \min(\eta^2,\eta)\Big\};
		\end{equation}
		\item[] (ii) for each $j,k = 1,\ldots,p,$ $l=1, \dots, d_j,$ $m=1, \dots, d_k,$ but $j \neq k$ or $l \neq m,$ and any $\eta >0$,
		\begin{equation}
		\label{con_bd_score2}
		\begin{split}
		&P\left\{ \left| \frac{\widehat \sigma_{jklm} - \sigma_{jklm}}{\lambda_{jl}^{1/2} \lambda_{km}^{1/2}}\right| \ge \rho_4 \cM_1(\mF_{\bX}) (l+m)^{\alpha + 1}  \eta
		+ \rho_5 \cM_1^2(\mF_{\bX})  (l+m)^{3\alpha + 2}\eta^2\right\} \\
		\le& 20 \exp\big\{ - c_3 n\min(\eta^2, \eta)\big\} + 8\exp\big\{ - c_3  \cM_1^{-2}(\mF_{\bX}) n l^{-2(\alpha + 1)}\big\},
		\end{split}
		\end{equation}
		where $\rho_1 = 16\sqrt{2} c_0^{-2} \alpha \lambda_0^2,$ $\rho_4$ and $\rho_5$ are defined in (\ref{thm_cov_rho}) in the Appendix.
	\end{description}
In particular, let $M$ be a positive integer possibly depending on $(n,p)$ with $\min_{1\leq j \leq p}\lambda_{jM}>0.$ If the sample size $n$ satisfies the lower bound 
$$n \geq \tilde\rho_6^2\log(pM)\max\big\{(\tilde \rho_5-2^{\alpha+1}\rho_4)^{-2}2^{6\alpha+4}\rho_5^2M^{4\alpha+2}\cM_1^2(\mF_{\bX}),1\big\},$$
where $\tilde\rho_5$ and $\tilde\rho_6$ are some constants with $\tilde\rho_5>2^{\alpha+1}\rho_4$ and $\tilde \rho_6 >c_3^{-1/2},$ then with probability greater than $1-20(pM)^{1 - c_3 \tilde\rho^{2}_6},$ the estimates $\big\{\widehat\sigma_{jklm}:j,k=1, \dots, p, l,m, =1, \dots, M\big\}$ satisfy the bound
\begin{equation}
\label{bd_score_max}
\underset{\underset{1 \le l,m \le M}{1\le j,k\le p}}{\max}\left| \frac{\widehat \sigma_{jklm} - \sigma_{jklm}}{\lambda_{jl}^{1/2} \lambda_{km}^{1/2}}\right| \le \tilde\rho_5 \tilde\rho_6\cM_1(\mF_{\bX}) M^{\alpha + 1} \sqrt{\frac{\log (pM)}{n}}.
\end{equation}
\end{theorem}

Likewise, the relative error bounds in Theorem~\ref{thm_cov} are also governed by $\cM_1(\mF_{\bX}).$ We provide three comments here. First, when $j=k$ and $l=m$ with $\sigma_{jjll}=\lambda_{jl},$ the concentration inequality in (\ref{con_bd_score1}) reduces to (\ref{thm_con_lambda}) corresponding to a faster exponential decay rate. Second, in the special case of no temporal dependence, \cite{qiao2018a} developed concentration results in terms of the absolute errors of ${\widehat\sigma_{jlkm}}.$ By contrast, we derive the concentration bounds on the relative errors in Theorem~\ref{thm_cov}. If we further assume that $\cM_1(\mF_{\bX})$ remains constant with respect to $(n,p)$ and $\lambda_{jl} \asymp l^{-\alpha}$ for each $j,$ as required in \cite{qiao2018a}, then (\ref{con_bd_score2}) would result in a sharper convergence rate of $|\widehat\sigma_{jklm}-\sigma_{jklm}|=O_P\big\{(l+m)n^{-1/2}\big\}$ as opposed to
$O_P\big\{(l+m)^{\alpha+1}n^{-1/2}\big\},$ established in \cite{qiao2018a}. Third, we can relax Condition~\ref{cond_eigen} by allowing different decay rates of lower bounds on the eigen-gaps across $j,$ denoted by $\alpha_j$ for $j=1,\dots,p,$ the resulting relative error rate becomes $\big|\lambda_{jl}^{-1/2}\lambda_{km}^{-1/2}(\widehat\sigma_{jklm}-\sigma_{jklm})\big|=O_P\big\{(l^{\alpha_j+1}+m^{\alpha_k+1})n^{-1/2}\big\}.$

\section{Vector functional autoregressive models} 
\label{sec.vfar}
Inspired from the standard VAR formulation, we propose a VFAR model of lag $L$, namely VFAR($L$), which is able to characterize linear interdependencies among multiple curve time series, as follows
\begin{equation}
\label{vfar1}
\bX_{t} = \sum_{h = 1}^L \boldsymbol{\cal A}_{h}(\bX_{t-h})  + \bvarepsilon_{t}, ~~t =L+1, \ldots,n,
\end{equation}
where the errors, $\bvar_t(\cdot)=\big(\varepsilon_{t1}(\cdot), \dots, \varepsilon_{tp}(\cdot)\big)^T, t=L+1, \dots, n,$ are 
i.i.d. from a $p$-dimensional mean zero Gaussian process, independent of $\bX_{t-1}(\cdot), \bX_{t-2}(\cdot), \dots,$ and $\boldsymbol{\cal A}_{h}: \cH^p \rightarrow \cH^p$ is the transition matrix operator at lag $h \in \{1, \dots,L\},$ induced from the kernel of the transition matrix function at lag $h,$ $\bA_h=\big(A_{hjk}(\cdot,\cdot)\big)_{1\leq j,k\leq p}\in {\eS}^{p \times p},$ through $\boldsymbol{\cal A}_{h}(\bX_{t-h})(u)= \int_{\cU}\bA_{h}(u,v) \bX_{t-h}(v)dv.$ The structure of transition matrix functions provides insights into the complex temporal and cross-sectional interrelationship  amongst $p$ curve time series. To make the problem of fitting (\ref{vfar1}) feasible in ``large $p,$ small $n$" scenarios, we assume functional sparsity in $\bA_1, \dots, \bA_L,$ that is most of the components in $\big\{X_{(t-h)j}(\cdot): h=1,\dots,L, j=1,\dots,p\big\}$ are unrelated to $X_{tj}(\cdot)$ for $j=1,\dots,p.$

Due to the infinite dimensional nature of functional data, for each $j,$ we take a standard dimension reduction approach through FPCA to approximate $X_{tj}(\cdot)$ using the leading $q_j$ principal components, i.e. $X_{tj}(\cdot)\approx \sum_{l=1}^{q_j}\xi_{tjl}\phi_{jl}(\cdot)=\bxi_{tj}^{\T}\bphi_{jl}(\cdot),$ where $\bxi_{tj}=(\xi_{tj1}, \dots, \xi_{tjq_j})^{\T}$, $\bphi_{j}(\cdot)=\big(\phi_{j1}(\cdot), \dots, \phi_{jq_j}(\cdot)\big)^{\T}$ and $q_j$ is chosen large enough to provide a reasonable approximation to the trajectory $X_{tj}(\cdot).$ 

Once the FPCA has been performed for each $X_{tj}(\cdot),$ we let $\bV_{hj} \in {\eR}^{(n-L)\times q_j}$ with its row vectors given by $\bxi_{(L+1-h)j}, \dots, \bxi_{(n-h)j}$ and
$\bPsi_{hjk}=\int_{\cU}\int_{\cU}\bphi_k(v)A_{hjk}(u,v)\bphi_j(u)^{\T}dudv \in {\eR}^{q_k \times q_j}.$ Then further derivations in Section~\ref{ap.vfar.deriv} of the Supplementary Material lead to the matrix representation of (\ref{vfar1}) as
\begin{equation}
\label{vfar3}
\bV_{0j} = \sum_{h = 1}^L \sum_{k=1}^p
\bV_{hk}\bPsi_{hjk} + \bR_j + \bE_j, ~~j=1,\dots,p,
\end{equation}
where $\bR_{j}$ and $\bE_j$ are $(n-L)$ by $q_j$ error matrices whose row vectors are formed by the truncation and random errors, respectively. Hence, we can rely on the block sparsity pattern in $\{\bPsi_{hjk}:h=1,\dots,L,j,k=1,\dots,p\}$ to recover the functional sparsity structure in $\{A_{hjk}:h=1,\dots,L,j,k=1,\dots,p\}.$ It is also worth noting that (\ref{vfar3}) can be viewed as a $\big(\sum_{j=1}^p q_j\big)$-dimensional VAR(L) model with the error vector consisting of both the approximation and random errors.

\subsection{Estimation procedure}
\label{sec.vfar.est}
The estimation procedure proceeds in the following three steps.

{\bf Step~1}. We perform FPCA based on observed curves, $X_{1j}(\cdot), \dots, X_{nj}(\cdot)$ and thus obtain estimated eigenfunctions $\widehat\bphi_{j}(\cdot) = \big(\widehat\phi_{jl}(\cdot), \dots, \widehat\phi_{jq_j}(\cdot)\big)^{\T}$ and FPC scores $\widehat\bxi_{tj} = \big(\widehat\xi_{tj1},\dots, \widehat\xi_{tjq_j}\big)^{\T}$ for $j=1, \dots, p.$ See Section~\ref{ap.select.tune} of the Supplementary Material for the selection of $q_j$'s in practice.

{\bf Step 2}. Motivated from the matrix representation of a VFAR(L) model in (\ref{vfar3}), we propose a penalized {\it least squares} (LS) approach, which minimizes the following optimization criterion over $\{\bPsi_{hjk}: h=1,\dots, L, k=1,\dots,p\}:$
\begin{equation}
    \label{vfar.crit}
    \frac{1}{2}\left\|\widehat\bV_{0j}-\sum_{h=1}^L \sum_{k=1}^p \widehat\bV_{hk}\bPsi_{hjk}\right\|_F^2 + \gamma_{nj} \sum_{h=1}^L \sum_{k=1}^p  \left\|\widehat\bV_{hk}\bPsi_{hjk}\right\|_F,
\end{equation}
where $\widehat\bV_{hj},$ the estimate of $\bV_{hj},$ is a $(n-L)$ by $q_j$ matrix with its $i$-th row vector given by $\widehat\bxi_{(L+i-h)j}$ for $i=1, \dots, (n-L),$ and $\gamma_{nj}$ is a non-negative regularization parameter. The $\ell_1/\ell_2$ type of standardized group lasso penalty \cite[]{simon2012} in (\ref{vfar.crit}) forces the elements of $\bPsi_{hjk}$ to either all be zero or non-zero. Potentially, one could modify (\ref{vfar.crit}) by adding an unstandardized group lasso penalty
\cite[]{yuan2006} in the form of $\gamma_{nj} \sum_{h=1}^L \sum_{k=1}^p  \left\|\bPsi_{hjk}\right\|_F$ to produce block sparsity in $\{\bPsi_{hjk}\}.$ However, orthonormalization within each group would correspond to the uniformly
most powerful invariant test for inclusion of a group \cite[]{simon2012}, hence we use a standardized group lasso penalty here. See  \cite{fan2015} for a similar ``fit penalized" framework by implementing an $\ell_1/\ell_2$-penalized LS approach to fit a functional additive model. In Section~\ref{ap.alg.vfar} of the Supplementary Material, we develop a block version of {\it fast iterative shrinkage-thresholding algorithm} (FISTA), which mirrors recent gradient-based techniques \cite[]{beck2009, dono2015}, to solve the optimization problem in (\ref{vfar.crit}) with the solution given by $\{\widehat \bPsi_{hjk}\}.$ The proposed FISTA algorithm is easy to implement and converges fast, thus is suitable for solving large-scale optimization problems.

{\bf Step~3}. Finally, we estimate elements of transition matrix functions in (\ref{vfar1}) by
\begin{equation}
\label{est.A}
\widehat A_{hjk}(u,v)=\widehat\bphi_k(v)^{\T}\widehat\bPsi_{hjk}\widehat\bphi_j(u), ~~h=1,\dots,L, ~ j,k=1,\dots,p.
\end{equation}

\subsection{Functional network Granger causality}
\label{sec.ngc}
In this section, we extend the definition of {\it network Granger causality} (NGC) under a VAR framework \cite[]{Bluk2005} to the functional domain and then use the extended definition under our proposed VFAR framework to understand the causal relationship among multiple curve time series.

In analogy to the NGC formulation, a {\it functional NGC} (FNGC) model consists of $p$ nodes, one for each functional variable, and a number of edges with directions connecting a subset of nodes. Specifically, curve times series of $\{X_{tk}(\cdot)\}_{t\in\eZ}$ is defined to be Granger causal for that of $\{X_{tj}(\cdot)\}_{t\in\eZ}$ or equivalently there is an edge from node $k$ to node $j,$ if $A_{hjk}(u,v) \neq 0$ for some $(u,v)\in \cU^2$ or $h\in \{1,\dots, L\}.$ Then our proposed FNGC model can be represented by a directed graph $G=(V,E)$ with vertex set $V=\{1, \dots, p\}$ and edge set
$$E=\left\{(k,j): A_{hjk}(u,v) \neq 0 \text{ for some } (u,v) \in \cU^2 \text{ or } h \in \{1,\dots, L\}, (j,k) \in V^2\right\}.$$ Hence, to explore the FNGC structure and the direction of influence from one node to the other, we need to develop an approach to estimate $E,$ i.e. identifying the locations of non-zero entries in $\widehat \bA_1,\dots, \widehat\bA_L,$ the details of which are presented in Section~\ref{sec.vfar.est}.

\subsection{Theoretical properties}
\label{sec.vfar.theory}
According to Section~\ref{ap.vfar1.rep} of the Supplementary Material, all VFAR($L$) models in (\ref{vfar1}) can be reformulated as a VFAR(1) model. Without loss of generality, we consider a VFAR(1) model in the form of $$\bX_t(u)=\int_{\cU}\bA(u,v)\bX_{t-1}(v)dv + \bvarepsilon_t(u), ~~ t=1, \dots,n.$$ 

To simplify our notation in this section, we focus on the setting where the $q_j$'s are the same across $j=1,\dots,p.$ However, our theoretical results extend naturally to the more general setting. In our empirical studies, we select different $q_j$'s, see Section~\ref{ap.select.tune} of the Supplementary Material for details. Let $\widehat\bZ=(\widehat\bV_{11}, \dots, \widehat\bV_{1p})^{
\T} \in {\eR}^{(n-1)\times pq},$ $
\bPsi_{j} = (\bpsi_{1j1}^\T, \dots, \bpsi_{1jp}^\T)^\T \in \eR^{pq \times q}.
$ and 
$\widehat\bD=\text{diag}\big(\widehat\bD_1,\dots, \widehat\bD_p\big) \in {\eR}^{pq\times pq},$ where $\widehat \bD_k=\big\{(n-1)^{-1}\widehat\bV_{1k}^{\T}\widehat\bV_{1k}\big\}^{1/2} \in {\eR}^{q \times q}$ for $k=1,\dots,p.$ Then minimizing (\ref{vfar.crit}) over $\bPsi_j\in {\eR}^{pq\times q}$ is equivalent to minimizing the following criterion over $\bB_{j} \in {\eR}^{pq\times q}$
\begin{eqnarray}
\label{vfar.crit.theory}
\langle\langle \widehat \bY_j, \bB_{j}\rangle\rangle + 
\frac{1}{2} \langle\langle \bB_{j}, \widehat \bGamma\bB_{j}\rangle\rangle  + \gamma_{nj}  \|\bB_{j}\|_{1}^{(q)},
\end{eqnarray}
where $\widehat \bY_j = (n-1)^{-1}\widehat\bD^{-1} \widehat{\bZ}^\T\widehat{\bV}_{0j},$  $\widehat \bGamma= (n-1)^{-1}\widehat\bD^{-1}\widehat{\bZ}^\T\widehat{\bZ} \widehat\bD^{-1}.$ Let $\widehat\bB_{j}$ be the minimizer of (\ref{vfar.crit.theory}), then $\widehat\bPsi_{j}=\widehat\bD^{-1}\widehat\bB_{j}$ with its $k$-th row block given by $\widehat\bPsi_{jk}$ and $\widehat\bA =\big\{\widehat A_{jk}(\cdot,\cdot)\big\}$ with its $(j,k)$-th entry being $\widehat A_{jk}(u,v)=\widehat \bphi_{k}(v)^{\T}\widehat \bPsi_{jk}\widehat\bphi_{j}(u)$ for $j,k=1,\dots,p$ and $(u,v)\in \cU^2.$ 

Before stating the condition on the entries of $\bA=\{A_{jk}(\cdot,\cdot)\},$ we begin with
some notation. For the $j$-th row of $\bA,$ we denote the set of non-zero functions by $S_j=\big\{k \in \{1,\dots,p\}: ||A_{jk}||_{\cS} \neq 0 \big\}$ and its cardinality by $s_j=|S_j|$ for $j=1,\dots,p.$ We also denote the maximum degree or row-wise cardinality by $s=\max_j s_j$ (possibly depends on $n$ and $p$), corresponding to the maximum number of non-zero functions in any row of $\bA.$ 

\begin{condition} 
\label{cond.fvar.bias}	
For each $j=1,\dots,p$ and $k\in S_j,$ $\tA_{jk}(u,v) = \sum_{l,m=1}^\infty a_{jklm} \phi_{jl}(u) \phi_{km}(v)
$
and there exist some positive constants $\beta > \alpha/2+1 $ and $\mu_{jk}$ such that $|a_{jklm}| \le \mu_{jk}(l+m)^{- \beta - 1/2}$ for $l,m \geq 1.$ 
\end{condition}

For each $(j,k)$, the basis with respect to which coefficients $\{a_{jklm}\}_{l,m \geq 1}$ are defined is determined by $\{\phi_{jl}(\cdot)\}_{l\geq 1}$ and $\{\phi_{km}(\cdot)\}_{m\geq 1}$ The basis $\phi_{j1}(\cdot), \phi_{j2}(\cdot),\dots $ is canonical in functional data problems, since it provides the unique basis with respect to which $X_{tj}(\cdot)$ can be expressed under a Karhunen-Lo\'eve expansion with uncorrelated coefficients. The parameter $\beta$ in Condition~\ref{cond.fvar.bias} determines the decay rate of the upper bounds for coefficients $\{a_{jklm}\}_{l,m \geq 1}.$ The assumption $\beta>\alpha/2+1$ can be interpreted as requiring each $A_{jk}(\cdot,\cdot)$ be sufficiently smooth relative to marginal-covariance functions, $\Sigma_{0,jj}(\cdot,\cdot)$ and $\Sigma_{0,kk}(\cdot,\cdot),$ the smoothness of which are indicated by the corresponding spectral decompositions with $\min(\lambda_{jl},\lambda_{kl}) \geq c_0\alpha^{-1}l^{-\alpha}$ from Condition~\ref{cond_eigen}. See \cite{hall2007} and \cite{qiao2018c} for similar smoothness conditions in functional linear models.

Next, we establish the consistency of the VFAR estimate based on the following sufficient conditions including a restricted eigenvalue (RE) condition and two deviation conditions in Conditions~\ref{cond.fvar.RE}--\ref{cond.fvar.max.error}, respectively. Finally, using the concentration results presented in Section~\ref{sec.bd.fpca}, we show that all stable VFAR models satisfy these three conditions with high probability. 

\begin{condition}
\label{cond.fvar.RE}
The symmetric matrix $\widehat \bGamma \in {\eR}^{pq \times pq}$ satisfies the RE condition with tolerance $\tau_1>0$ and curvature $\tau_2>0$ if
\begin{equation}
\label{cond_RE}
\btheta^\T \widehat \bGamma \btheta \geq  \tau_2 \|\btheta\|^2 -\tau_1 \|\btheta\|_1^2 \quad \forall \btheta \in {\eR}^{pq}.
\end{equation}
\end{condition}

\begin{condition}
\label{cond.fvar.eigen}	
	There exist some positive constants $C_{\lambda}$ and $C_{\phi}$ such that 
    \begin{equation}
    \label{fvar.eigen}	
    \begin{split}
    \max_{1 \le j \le p, 1 \le l \le q} \left|\frac{\widehat{\lambda}_{jl}^{-1/2} - \lambda_{jl}^{-1/2}}{\lambda_{jl}^{-1/2}}\right| \le C_{\lambda}  \cM_1(\mF_{\bX})\sqrt{\frac{\log (pq)}{n}},\\
    \max_{1 \le j \le p, 1 \le l \le q}\|\widehat \phi_{jl} - \phi_{jl}\| \le C_{\phi} \cM_1(\mF_{\bX})q^{\alpha + 1} \sqrt{\frac{\log (pq)}{n}}.
    \end{split}
    \end{equation}
\end{condition}

\begin{condition}
	\label{cond.fvar.max.error}
	There exits some positive constant $C_{E}$ such that
	\begin{equation}
	\label{fvar.max.error}
	\Big\|\widehat \bY_j- \widehat \bGamma\bB_j\Big\|_{\max}^{(q)} \leq C_E \cM_1(\mF_{\bX}) s_j \Big\{q^{\alpha + 2}\sqrt{\frac{\log (pq)}{n}} + q^{-\beta +1} \Big\},~j=1,\ldots,p.
	\end{equation}
	\end{condition}

We are now ready to present the theorem on the convergence rate of the VFAR estimate in a deterministic design, i.e. we assumes a fixed realization of $\bX_1(\cdot), \dots, \bX_n(\cdot).$
\begin{theorem}
\label{thm.vfar}
    Suppose that Conditions \ref{cond.bd.fsm}--\ref{cond.fvar.max.error} hold with $\tau_2\ge 32\tau_1 q^2 s.$ Then, for any regularization parameter, $\gamma_{nj} \ge 2C_E \cM_1(\mF_{\bX}) s_j\big\{q^{\alpha + 2}(\log(pq)/n)^{1/2} + q^{-\beta +1} \big\},$ $\gamma_n=\max_{1\leq j\leq p}\gamma_{nj}$ and $q^{\alpha/2}s\gamma_n\to 0$ as $n, p \to \infty,$ any minimizer $\widehat\bB_j$ of (\ref{vfar.crit.theory}) satisfies
$$
\|\widehat \bB_{j} - \bB_{j}\|_F \le \frac{24 s_j^{1/2}\gamma_{nj}}{\tau_2}, 
~~\|\widehat \bB_{j} - \bB_{j}\|_{1}^{(q)} \le \frac{96s_j\gamma_{nj}}{\tau_2} ~ \text{ for }j =1,\ldots,p,
$$
and the estimated transition matrix function, $\widehat \bA,$ satisfies 
\begin{equation}
\label{err.A}
\|\widehat \bA - \bA\|_\infty \leq \frac{96\alpha^{1/2}q^{\alpha/2}s\gamma_{n}}{c_0^{1/2}\tau_2}\Big\{1+o(1)\Big\}.
\end{equation}
\end{theorem}

The convergence rate of $\widehat\bA$ under the functional matrix $\ell_{\infty}$ norm is governed by two sets of parameters: (1) dimensionality parameters, sample size ($n$), number of functional variables ($p$), and maximum row-wise cardinality ($s$) in $\bA$; (2) internal parameters, functional stability measure $\big(\cM_1(\mF_{\bX})\big),$ the truncated dimension of curve time series ($q$), tolerance ($\tau_1$), curvature ($\tau_2$), decay rate of the lower bounds for eigenvalues ($\alpha$) and decay rate of the upper bounds for transition matrix basis coefficients ($\beta$). In the following, we provide three remarks for Theorem~\ref{thm.vfar}. First, it is easy to see that larger values of $\alpha$ (tighter eigen-gaps) or $\cM_1(\mF_{\bX})$ (less stable process of $\big\{\bX_t(\cdot)\big\}$) or $s$ (denser structure in $\bA$) yield a slower convergence rate, while enlarging $\beta$ or $\tau_2$ will increase the entrywise smoothness in $\bA$ or the curvature of the RE condition, respectively, thus resulting in a faster rate. Second, the convergence rate consists of two terms corresponding to the variance-bias tradeoff as commonly considered in nonparametric statistics. Specifically, the variance is of the order $O_P\big[\cM_1(\mF_{\bX})s^2q^{(3\alpha+4)/2}\{\log(pq)/n\}^{1/2}\big]$ and the bias term is bounded by $O\big\{\cM_1(\mF_{\bX})s^2q^{(\alpha-2\beta+2)/2}\big\}.$ To balance both terms, we can choose the optimal truncation dimension $q$ satisfying $\log(pq)q^{2\alpha+2\beta+2}\asymp n.$ Third, in the special case when each $X_{tj}(\cdot)$ is finite dimensional, although the truncation step is no longer required, the FPC scores still need to be estimated. The convergence rate then reduces to $O_P\big\{\cM_1(\mF_{\bX})s^2(\log(p)/n)^{1/2}\big\},$ which is slightly different from the rate of the high-dimensional VAR estimate in \cite{basu2015a}.

Next, we turn to the case of random designs, in which $\bX_t(\cdot),\dots,\bX_n(\cdot)$ are drawn from random ensembles. We need to verify that Conditions~\ref{cond.fvar.RE}--\ref{cond.fvar.max.error} are satisfied with high probability in the following Propositions~\ref{res.re}--\ref{prop_Error_max}, respectively. We now introduce a sufficient condition for Proposition~\ref{res.re}.

\begin{condition}
\label{cond_min_bound}
For $\bSigma_0=\big(\bSigma_{0,jk}(\cdot,\cdot)\big)_{1\leq j,k\leq p},$ we denote a diagonal matrix function by $\bD_0 = \text{diag}(\Sigma_{0,11}, \dots, \Sigma_{0,pp}) \in {\eS}^{p \times p}$ with its induced operator given by $\bcD_0.$ The infimum $\underline{\mu}$ of the functional Rayleigh quotient of $\bXi_0$ relative to $\bcD_0,$ defined as follows, is bounded below by zero, i.e.
\begin{eqnarray}
\underline{\mu} = \inf_{\bPhi \in \widebar{\cH}_0^p} \frac{\big \langle\bPhi, \bXi_0(\bPhi) \big \rangle}
{\big \langle\bPhi, \bcD_0(\bPhi) \big \rangle}  > 0,
\end{eqnarray}
where $\widebar{\cH}_0^p = \{\bPhi \in \cH^p: \big \langle\bPhi, \bcD_0(\bPhi) \big \rangle \in (0, \infty)\}.$
\end{condition}

Here the value $\underline{\mu},$ chosen as the curvature $\tau_2$ in the proof of Proposition~\ref{res.re}, can be understood as requiring the minimum eigenvalue of the correlation matrix function for $\bX_t(\cdot)$ to be bounded below by zero. In particular, if $X_{tj}(\cdot)$ is $d_j$-dimensional for $j=1,\dots,p$ with $d=\max_jd_j<\infty,$ it is easy to show that $\underline{\mu}$ reduces to the minimum eigenvalue of the correlation matrix for the $\sum_jd_j$-dimensional vector,  $\bxi_t=(\xi_{t11}, \dots, \xi_{t1d_1}, \dots, \xi_{tp1}, \dots, \xi_{tpd_p})^{\T}.$ Correspondingly, we can define the supremum $\overline{\mu}$ of the functional Rayleigh quotient of $\bXi_0$ relative to $\bcD_0$ by
$
\overline{\mu} = \sup_{\bPhi \in \widebar{\cH}_0^p} \frac{ \langle\bPhi, \bXi_0(\bPhi) \rangle}
{ \langle\bPhi, \bcD_0(\bPhi) \rangle},
$
which corresponds to the maximum eigenvalue of the correlation matrix function for $\bX_t(\cdot)$ or the correlation matrix for $\bxi_t$ depending on whether the dimension of each $X_{tj}(\cdot)$ is infinite or not. 

\begin{proposition} (Verify Condition~\ref{cond.fvar.RE})
	\label{res.re}
	Suppose that Conditions \ref{cond.bd.fsm}-\ref{cond_eigen} and \ref{cond_min_bound} hold. 
	Then there exist three positive constants $C_{\Gamma}$, $c_4$ and $c_5$ such that, for $n \succsim \log p+\log q,$ 
	\begin{eqnarray*}
		\btheta^T \widehat \bGamma \btheta \ge  \underline{\mu} \big\|\btheta\big\|^2_2 - C_\Gamma \cM_1(\mF_{\bX}) q^{\alpha +1}\sqrt{\frac{\log (pq)}{n}}\big\|\btheta\big\|_1^2
	\end{eqnarray*}
with probability greater than $1 - c_4(pq)^{-c_5}$. 
\end{proposition}
\begin{proposition} (Verify Condition~\ref{cond.fvar.eigen})
\label{prop_Error_eigen} 
Suppose that Conditions~\ref{cond.bd.fsm}--\ref{cond_eigen} hold. Then there exists four positive constants $C_{\phi}$, $C_{\lambda}$, $c_4$ and $c_5,$ such that, for $n \succsim \log p+\log q,$ the two deviation bounds in (\ref{fvar.eigen}) hold with probability greater than $1 -c_4(pq)^{-c_5}.$ 
\end{proposition}

\begin{proposition} (Verify Condition~\ref{cond.fvar.max.error})
\label{prop_Error_max} 
    Suppose that Conditions \ref{cond.bd.fsm}-\ref{cond.fvar.bias} hold. Then there exist three positive constants $C_E$, $c_4$ and $c_5,$ such that, for $n \succsim \log p+\log q,$ the deviation bound in (\ref{fvar.max.error}) holds
	with probability greater than $1 -c_4(pq)^{-c_5}.$ 
\end{proposition}

Propositions \ref{res.re}--\ref{prop_Error_max} can be proved by applying the concentration results in Theorems~\ref{thm_lambda}--\ref{thm_cov}. Here we choose suitable common constants $c_4, c_5$ and sufficiently large $n,$ as stated in Propositions \ref{res.re}--\ref{prop_Error_max}, such that the joint probability for the three events corresponding to the nonasymptotic upper bounds in 
(\ref{bd_lambda_max}), (\ref{bd_phi_max}) and (\ref{bd_score_max}), respectively, is greater than $1-c_4(pq)^{-c_5}.$ As a consequence, if the sample size $n \succsim \log p+\log q,$ then with probability greater than $1-c_4(pq)^{-c_5},$ the estimate $\widehat\bA$ satisfies the error bound in (\ref{err.A}).

\subsection{Simulation studies}
\label{sec.sim}
In this section, we conduct a number of simulations to compare the finite sample performance of our proposed method to potential competitors. 

In each simulated scenario, we generate functional variables by $X_{tj}(u)=\bs(u)^{\T}\btheta_{tj},$ $j=1,\dots,p, u \in \cU=[0,1],$ where $\bs(\cdot)$ is a $5$-dimensional orthonormal Fourier basis function and each $\btheta_{t}=(\btheta_{t1}^{\T},\dots,\btheta_{tp}^{\T} )^{\T} \in \eR^{5p}$ is generated from a stationary VAR(1) process, $\btheta_t=\bB\btheta_{t-1}+\bfeta_{t},$ with transition matrix $\bB \in \eR^{5p\times 5p},$ whose $(j,k)$-th block is given by $\bB_{jk}, j,k=1,\dots,p$ and innovations $\bfeta_{t}$'s being randomly sampled from $N({\bf 0}, \bI_{5p}).$ The observed values, $W_{tjs},$ are then generated, with measurement errors, from $W_{tjs}=X_{tj}(u_s)+e_{tjs},$ at $T=50$ equally spaced time points, $0=u_1,\dots,u_T=1$ with errors $e_{tjs}$'s being randomly sampled from $N(0,0.5^2).$ In our simulations, we generate $n=100$ or $n=200$ observations of $p=40$ or $p=80$ functional variables, and we aim to show that, although our method is developed for fully observed curve time series, it still works well even with measurement error.

According to Section~\ref{ap.vfar1.sim} of the Supplementary Material, $\bX_t(\cdot)$ follows from a VFAR(1) model, $\bX_t(u)=\int_{\cU}\bA(u,v)\bX_{t-1}(v)dv+\bvarepsilon_t(u),$ where $\varepsilon_{tj}(u)=\bs(u)^{\T}\bfeta_{tj}$ and autocoefficient functions satisfy $A_{jk}(u,v)=\bs(u)^{\T}\bB_{jk}\bs(v)$ for $j,k=1,\dots,p, (u,v) \in \cU^2.$ Hence,
the functional sparsity structure in $\bA$ can be correspondingly characterized by the block sparsity pattern in $\bB.$ We consider two different scenarios for $\bB$ as follows.
\begin{itemize}
    \item[(i)] {\bf Block sparse}. We generate a block sparse $\bB$ without any special structure. Specifically, we generate  $\bB_{jk}=w_{jk}\bC_{jk}$ for $j,k=1\dots,p,$ where entries in $\bC_{jk}$ are randomly sampled from $N(0,1)$ and $w_{jk}$'s are generated from $\{0,1\}$ under the constraint of $\sum_{k=1}^p w_{jk}=5$ for each $j,$ such that the same row-wise block sparsity level for $\bB$ can be produced. To guarantee the stationarity of $\{\bX_t(\cdot)\},$ we rescale $\bB$ by $\kappa\bB/\rho(\bB),$ where $\kappa$ is generated from Unif[0.5,1].
    
    \item[(ii)] {\bf Block banded}. We generate a block banded $\bB,$ with entries in $\bB_{jk}$ being randomly sampled from $N(0,1)$ if $|j-k|\leq 2,$ and being zero at other locations. $\bB$ is then rescaled as described in (i).
\end{itemize}

We perform regularized FPCA on each function and use $5$-fold cross-validation to choose $q_j,$ the details of which are discussed in Sections~\ref{ap.select.tune} and \ref{ap.rfpca} of the Supplementary Material. Typically $q_j=$ 4, 5 or 6 are selected in our simulations. We compare our proposed $\ell_1/\ell_2$-penalized LS estimate using all selected principal components, namely $\ell_1/\ell_2$-$\text{LS}_{\text{a}},$ to its two competitors. One method, $\ell_1/\ell_2$-$\text{LS}_{2},$ relies on minimizing $\ell_1/\ell_2$-penalized LS based on the first two estimated principal components, which capture partial curve information. The other approach, $\ell_1$-$\text{LS}_{1},$ project the functional data into a standard format by computing the first estimated FPC score and then implement an $\ell_1$ regularization approach \cite[]{basu2015a} for VAR estimation on this data. We examine the sample performance of three approaches, $\ell_1/\ell_2$-$\text{LS}_{\text{a}},$ $\ell_1/\ell_2$-$\text{LS}_{2},$ and $\ell_1$-$\text{LS}_{1}$ in terms of model selection consistency and estimation accuracy.   
\begin{itemize}
    \item {\bf Model selection}. We plot the true positive rates against false positive rates, respectively defined as 
    $\frac{\#\{(j,k):  ||\widehat A_{jk}^{(\bgamma_n)}||_{\cS}\neq0 \text{ and } ||A_{jk}||_{\cS}\neq 0\}}{\#\{(j,k): ||A_{jk}||_{\cS}\neq 0\}}$ and 
    $\frac{\#\{(j,k):  ||\widehat A_{jk}^{(\bgamma_n)}||_{\cS}\neq0 \text{ and } ||A_{jk}||_{\cS}= 0\}}{\#\{(j,k): ||A_{jk}||_{\cS}= 0\}},$ over a sequence of $\bgamma_n=(\gamma_{n1}, \dots, \gamma_{np})$ values to produce a ROC curve. We compute the {\it area under the ROC curve} (AUROC) with values closer to one indicating better performance in recovering the functional sparsity structure in $\bA.$
    
    \item {\bf Estimation error}. We calculate the relative estimation accuracy for $\widehat\bA$ by $\|\widehat{\bA} - \bA\|_F/\|\bA\|_F,$ where $\widehat\bA$ is the regularized estimate based on the optimal regularization parameter selected by minimizing AIC or BIC. See Section~\ref{ap.select.tune} of the Supplementary Material for details.
\end{itemize}

To investigate the support recovery consistency, Table~\ref{auroc.table} reports the average AUROCs of three comparison methods for both settings. 
In all simulations, we observe that $\ell_1/\ell_2$-$\text{LS}_{\text{a}}$ with most of curve information being captured, provides highly significant improvements over its two competitors and $\ell_1$-$\text{LS}_{1}$ gives the worst results. See Section~\ref{ap.sim} of the Supplementary Material for the graphical illustration of results, which is consistent to our findings in Table~\ref{auroc.table}. 

\begin{table}[ht]
	\caption{\label{auroc.table} The mean and standard error (in parentheses) of AUROCs over 100 simulation runs. The best values are in bold font.}
	\begin{center}
	\vspace{-0.1cm}
		\resizebox{6.1in}{!}{
			\begin{tabular}{ccccccc}		
				\hline 
				 & &  Model~(i)& & & Model~(ii) &  \\\hline
				$(n,p)$& $\ell_1/\ell_2$-$\LS_a$  & $\ell_1/\ell_2$-$\LS_2$ & $\ell_1$-$\LS_1$ & $\ell_1/\ell_2$-$\LS_a$ & $\ell_1/\ell_2$-$\LS_2$ & $\ell_1$-$\LS_1$ \\

				$(100,40)$ & {\bf 0.840(0.018)} & 0.690(0.019) & 0.591(0.023) & {\bf 0.872}(0.016) & 0.719(0.022) & 0.609(0.024)\\
				$(100,80)$ & {\bf 0.829(0.015)} & 0.682(0.017) & 0.585(0.015) & {\bf 0.869}(0.014) & 0.714(0.017) & 0.600(0.017) \\ 
				$(200,40)$ & {\bf 0.951(0.011)} & 0.764(0.020) & 0.616(0.021) & {\bf 0.971}(0.006) & 0.795(0.018) & 0.639(0.023) \\ 
				$(200,80)$ & {\bf 0.948(0.010)} & 0.770(0.017) & 0.626(0.015) & {\bf 0.969}(0.005) & 0.799(0.014) & 0.644(0.015) \\
				\hline 
			\end{tabular}
		}	
	\end{center}
		\vspace{-0.5cm}
\end{table}

To evaluate the estimation accuracy, Table~\ref{err.table} presents numerical results of relative errors of regularized estimates. We also report the performance of the LS estimate in the oracle case, where we know locations of non-zero entries of $\bA$ in advance. Several conclusions can be drawn from Table~\ref{err.table}. First, in all scenarios, the proposed BIC-based $\ell_1/\ell_2$-$\text{LS}_{\text{a}}$ method provides the highest estimation accuracy among all the comparison methods. Second, the performance of AIC-based methods severally deteriorate in comparison with their BIC-based counterparts. Given the high dimensional and functional natural of the model structure, computing effective degrees of freedom for the VFAR estimate leads to a very challenging task and requires further investigation. Third, $\LS_{\text{oracle}}$ estimates give much worse results than BIC-based regularized estimates. This is not surprising, since even in the ``large $n,$ small $p,$" scenario, e.g. $n=100, p=40$ under Model~(i), implementing LS requires estimating $5\times5^2=125$ parameters based on only 100 observations, which intrinsically results in a high dimensional estimation problem.


\begin{table}[ht]
  	\caption{\label{err.table} The mean and standard error (in parentheses) of relative estimation errors of $\widehat\bA$ over 100 simulation runs. The best values are in bold font.}
	\begin{center}
	\vspace{-0.1cm}
		\resizebox{6.1in}{!}{
  		\begin{tabular}{ccccccc}
  			\hline 
  			Model &$(n,p)$ & $\gamma_{jn}$ & $\ell_1/\ell_2$-$\LS_{\text{a}}$  & $\ell_1/\ell_2$-$\LS_2$ & $\ell_1$-$\LS_1$ & $\LS_{\text{oracle}}$\\
  			\hline
\multirow{8}{*}{(i)}&\multirow{2}{*}{(100,40)} & $\aic$& 1.783(0.044) & 1.455(0.032) & 1.047(0.008) & \multirow{2}{*}{1.483(0.044)} \\ 
                          && $\bic$& {\bf 0.971(0.007)} & 0.997(0.002) & 1.002(0.002) & \\ 
&\multirow{2}{*}{(100,80)} & $\aic$& 1.546(0.037) & 1.576(0.030) & 1.124(0.013) & \multirow{2}{*}{1.530(0.050)} \\ 
 &                         & $\bic$& {\bf 0.991(0.003)} & 0.999(0.001) & 1.002(0.001) &  \\ 
&\multirow{2}{*}{(200,40)} & $\aic$& 1.419(0.035) & 1.159(0.016) & 1.020(0.004) & \multirow{2}{*}{0.990(0.029)}  \\ 
&                          & $\bic$& {\bf 0.850(0.012)} & 0.989(0.003) & 0.999(0.002) & \\ 
&\multirow{2}{*}{(200,80)} & $\aic$& 1.544(0.045) & 1.350(0.022) & 1.032(0.004) & \multirow{2}{*}{1.016(0.033)}  \\ 
 &                         & $\bic$& {\bf 0.915(0.013)} & 0.994(0.002) & 1.000(0.001) & \\ \hline
\multirow{8}{*}{(ii)}&\multirow{2}{*}{(100,40)} & $\aic$ & 1.679(0.040) & 1.400(0.026) & 1.039(0.008) & \multirow{2}{*}{1.363(0.039)}  \\ 
 &                         & $\bic$ & {\bf 0.957(0.009)} & 0.995(0.002) & 1.000(0.002) & \\ 
&\multirow{2}{*}{(100,80)} & $\aic$ & 1.435(0.028) & 1.489(0.020) & 1.105(0.012) & \multirow{2}{*}{1.383(0.036)}    \\ 
  &                        & $\bic$ & {\bf 0.983(0.004)} & 0.998(0.001) & 1.001(0.001) &  \\  
&\multirow{2}{*}{(200,40)} & $\aic$ & 1.329(0.028) & 1.128(0.013) & 1.014(0.004) & \multirow{2}{*}{0.909(0.024)}   \\ 
 &                         & $\bic$ & {\bf 0.824(0.009)} & 0.985(0.003) & 0.998(0.002) & \\ 
&\multirow{2}{*}{(200,80)} & $\aic$ & 1.428(0.027) & 1.296(0.013) & 1.024(0.004) & \multirow{2}{*}{0.926(0.019)} \\ 
 &                         & $\bic$ & {\bf 0.881(0.010)} & 0.990(0.002) & 0.998(0.001) &  \\
  			\hline 
			\end{tabular}
		}	
	\end{center}
		\vspace{-0.5cm}
\end{table}

\section{Discussion}
\label{sec.discussion}

We identify several potential directions for future research. First, our theoretical results are developed for multivariate stationary Gaussian processes. For non-Gaussian processes, one needs to control higher order dependence possibly by using the higher order spectra. It is interesting to develop suitable concentration results for multivariate non-Gaussian and/or non-stationary processes, which would pose non-trivial theoretical challenges. The second extension considers applying our derived concentration results on other important statistical models for large-scale curve time series, e.g. functional additive models \cite[]{fan2015} and functional extension of factor models \cite[]{lam2012}. Third, it is well known that $\widehat\bSigma_{0}$ is not a consistent estimator for $\bSigma_0$ in ``large $p,$ small $n$" regimes. It is of great interest to develop consistent estimators through various types of regularization under the high dimensional, functional and dependent setting we consider in the paper. These topics are beyond the scope of the current paper and will be pursued elsewhere.


\section*{Appendix}
\setcounter{equation}{0}
\renewcommand{\theequation}{A.\arabic{equation}}

\appendix
\section{Technical proofs}
\label{asec.pf}
In the following, write $ \langle \phi, K\rangle $, $\langle K, \phi \rangle $ and $ \langle \phi_1,\langle K, \phi_2\rangle \rangle $ for $$
\int_{\cU} K(u,v)\phi(u)du,~~ \int_{\cU} K(u,v) \phi(v)dv~~\mbox{and}~~ \int_{\cU} \int_{\cU} K(u,v) \phi_1(u) \phi_2(v)dudv,
$$
respectively.

\subsection{Proof of Theorem \ref{thm_con_op}}
\label{asec.pf.thm1}
(i) Define $\bY = \big(\langle \bPhi_1,\bX_1\rangle, \ldots, \langle \bPhi_1,\bX_n\rangle\big)^\T$, then $\bY \sim N(\0, \bQ)$, where $Q_{rs} = \big\langle \bPhi_1, \bXi_{r-s}(\bPhi_1)\big\rangle$ for $r,s =1,\ldots,n$.
Note that $\big\langle\bPhi_1, \widehat{\bXi}_0(\bPhi_1)\big\rangle = {n^{-1}} \bZ^{\T} \bQ \bZ$, where $\bZ \sim N(\0,\bI_{n}),$ and $\big\langle\bPhi_1, {\bXi}_0(\bPhi_1)\big\rangle = \tE\big(n^{-1}\bZ^{\T} \bQ \bZ\big).$
By the Hanson-Wright inequality of \cite{rudelson2013}, we obtain that
$$
P\left\{\Big|\big\langle \bPhi_1, \big(\widehat{\bXi}_0- \bXi_0\big)(\bPhi_1)\big\rangle\Big|>\epsilon\right\} \le 2 \exp\left\{- c \min\left( \frac{n^2 \epsilon^2}{\|\bQ\|_F^2}, \frac{n\epsilon}{\|\bQ\|}\right)\right\}
$$
for some constant $c> 0.$ By $\|\bQ\|_F^2/n \le \|\bQ\|^2$ and letting $\epsilon = \eta\|\bQ\|$, we obtain that
\begin{equation}
\label{coneq_Phi_1}
P\left\{\Big|\big\langle \bPhi_1, \big(\widehat{\bXi}_0- \bXi_0\big)(\bPhi_1)\big\rangle\Big|> \eta\|\bQ\|\right\} \le 2 \exp\left\{- c n\min\left( \eta^2,\eta\right)\right\}
\end{equation}
for some universal constant $c> 0$.

Next we derive an upper bound on the operator norm $\|\bQ\|.$ Specifically, for any $\bw = (w_1,\ldots,w_n)^T \in \eR^n$ with $ ||\bw||=1,$ define $G_{\bw}(\theta) = \sum_{r =1 }^n w_r \exp(-i r\theta)$ and its conjugate by $G^*_{\bw}(\theta).$ Then we obtain that
\begin{eqnarray*}
	\bw^{\T}\bQ\bw &=& \sum_{r = 1}^n \sum_{s =1 }^n w_r w_s \big\langle \bPhi_1, \bXi_{r-s}(\bPhi_1)\big\rangle\\
	& =& \sum_{r = 1}^n \sum_{s =1 }^n w_r w_s \int_{-\pi}^{\pi}\big\langle\bPhi_1, \mF_{\bX,\theta}(\bPhi_1)\big\rangle \exp \{i(r-s)\theta\} d\theta \\
	&=& \int_{-\pi}^{\pi}\big\langle\bPhi_1, \mF_{\bX,\theta}(\bPhi_1)
	\big\rangle G_{\bw}(\theta) G^{*}_{\bw}(\theta)d\theta,
\end{eqnarray*}
where the second line follows from the inversion formula (\ref{inv.formula}).
For a fixed $p$-dimensional vector $\bPhi\in \cH^p$, denote 
$
\cM (\mF_{\bX},{\bPhi}) = 2\pi \cdot {\text{ess} \sup}_{\theta\in [-\pi, \pi] }\big|\langle \bPhi,\mF_{\bX,\theta}(\bPhi) \rangle\big|.
$
Since $\big\langle\bPhi_1, \mF_{\bX,\theta}(\bPhi_1)\big\rangle$ is Hermitian and $\int_{-\pi}^{\pi}G_{\bw}(\theta) G^{*}_{\bw}(\theta)d\theta = 2\pi$, we have
$
\|\bQ\| \le \cM(\mF_{\bX},\bPhi_1).
$
Then it follows from the definition in (\ref{eq.sub.fsm}) that
$$
\|\bQ\| \le \cM(\mF_{\bX},\bPhi_1) \le \cM_k(\mF_{\bX}) \big\langle \bPhi_1, \bXi_0(\bPhi_1)\big\rangle.
$$
This result, together with (\ref{coneq_Phi_1}) implies (\ref{thm_Phi_1}).

(ii) Note that 
\begin{eqnarray*}
	4\big\langle \bPhi_1, (\widehat{\bXi}_0 - \bXi_0)(\bPhi_2)\big\rangle
	&\leq &　\big\langle \widetilde{\bPhi}_1, (\widehat{\bXi}_0 - \bXi_0)(\widetilde{\bPhi}_1)\big\rangle 
	 - \big\langle \widetilde{\bPhi}_2, (\widehat{\bXi}_0 - \bXi_0)(\widetilde{\bPhi}_2)\big\rangle.
\end{eqnarray*}
where $\widetilde{\bPhi}_1 = \bPhi_1 + \bPhi_2,$ $\widetilde{\bPhi}_2 = \bPhi_1 - \bPhi_2$
and
$\cM(\mF_{\bX},\widetilde{\bPhi}_i) \le 2\{\cM(\mF_{\bX},\bPhi_1)+ \cM(\mF_{\bX}, \bPhi_2) \}$ for $i=1, 2.$ 
Combing these with results in (i) leads to
\begin{eqnarray*}
&&P\Big[\Big|\big\langle \bPhi_1, (\widehat{\bXi}_0 - \bXi_0)(\bPhi_2)\big\rangle\Big|> \{\cM(\mF_{\bX},\bPhi_1)+ \cM(\mF_{\bX}, \bPhi_2)\}\eta\Big]\\
&\leq& \sum_{i=1}^2P\Big[\Big|\big\langle \widetilde\bPhi_i, (\widehat{\bXi}_0 - \bXi_0)(\widetilde\bPhi_i)\big\rangle\Big|> \cM(\mF_{\bX},\widetilde\bPhi_i)\eta\Big]\le 4 \exp\Big\{- c n\min\left( \eta^2,\eta\right)\Big\}
\end{eqnarray*}
for some universal constant $c> 0.$ This, together with, $\cM(\mF_{\bX},\bPhi_i) \le \cM_k(\mF_{\bX}) \big\langle \bPhi_i, \bXi_0(\bPhi_i)\big\rangle$ for $i=1, 2,$ implies (\ref{thm_Phi_12}), which completes the proof.  $\square$

\subsection{Proofs of Theorems \ref{thm_max_bound} and \ref{thm_F_bound}}
\label{asec.pf.thm2}

First, we derive the concentration bound on $\|\widehat \Sigma_{0,jk} - \Sigma_{0,jk}\|_{\cS}$ for each $j$ and $k$.
Let $d_j = \big |\{l: \lambda_{jl}> 0\}\big|$ for $j =1,\ldots,p$ and 
$ \Delta_{jklm} = (\lambda_{jl} \lambda_{km})^{-1/2}\big\langle  \phi_{jl}, \langle \widehat \Sigma_{0,jk} - \Sigma_{0,jk}, \phi_{km} \rangle \big \rangle$ 
 for $l = 1,\ldots,d_j, m = 1,\ldots, d_k,$ $j, k= 1,\ldots,p.$ If $d_j < \infty$ (or $d_k < \infty$), define 
 $\Delta_{jklm} = 0 $ for $l \ge d_{j+1}$ (or $m \ge d_{k+1}$).
Clearly we have that
$
\big\|\widehat{\Sigma}_{0,jk} - \Sigma_{0,jk} \big\|_{\mathcal{S}}^2
= \sum^\infty_{l,m =1} \lambda_{jl}\lambda_{km} \Delta_{jklm}^2.
$
Then by Jensen's inequality, we have that
\begin{equation}
\tE\Big\{ \big\|\widehat{\Sigma}_{0,jk} - \Sigma_{0,jk} \big\|_{\mathcal{S}}^{2q}\Big\} 
\le  \Big( \sum_{l,m= 1}^{\infty} \lambda_{jl} \lambda_{km}\Big)^ {q-1} \sum_{l,m = 1}^{\infty} 
\lambda_{jl} \lambda_{km} \tE\big|\Delta_{jklm}  \big|^{2q}\le  \lambda_0^{2q} \sup_{l,m} \tE\big|\Delta_{jklm}  \big|^{2q}.
\label{Sigma.bd}
\end{equation}

For any given $(j, k,l,m),$ let $$
\bPhi_{1} = (0, \ldots,0,\lambda_{jl}^{-1/2}\phi_{jl}, 0,\ldots,0)^{\T} \text{ and } \bPhi_{2} = (0, \ldots,0,\lambda_{km}^{-1/2}\phi_{km}, 0,\ldots,0)^{\T}.
$$
By the definition of $\Delta_{jklm}$ and orthonormality of $\{\phi_{jl}(\cdot)\}$ and $\{\phi_{km}(\cdot)\}$ for each $j,k=1,\dots,p,$ we have 
$
\Delta_{jklm}  = \langle\bPhi_{1}, (\widehat{\bXi}_0 - \bXi_0)(\bPhi_{2})\big\rangle,
$
$\langle \bPhi_{1},\bXi_0(\bPhi_{1})\rangle=\langle \bPhi_{2},\bXi_0(\bPhi_{2})\rangle=1.$ Applying (\ref{thm_Phi_12}) in Theorem~\ref{thm_con_op}, we can obtain that
	\begin{equation}
	\label{con_Delta}
	P\Big\{ \big|\Delta_{jklm}  \big| > 2 \cM_1(\mF_{\bX}) \eta 
	\Big\} \le 4 \exp\Big\{- c n \min(\eta^2, \eta)\Big\},
	\end{equation}
for $j, k= 1,\ldots,p,$ $l = 1,\ldots,\Lambda_j$ and $m = 1,\ldots, \Lambda_k.$ It then follows from Lemma~\ref{lemma.moment}  of the Supplementary Material that for each integer $q \ge 1$,
$$
\left\{2 \cM_1(\mF_{\bX})  \right\}^{-2q} \tE\big|\Delta_{jklm}  \big|^{2q} \le q! 4(4c^{-1}n^{-1})^q + 4(2q)!(4c^{-1}n^{-1})^{2q}.
$$
This together with (\ref{Sigma.bd}) implies that 
\begin{equation}
\label{Sigma_moment}
	\left(2 \cM_1(\mF_{\bX}) \lambda_0 \right)^{-2q} \tE\Big\{\big\|\widehat{\Sigma}_{0,jk} - \Sigma_{0,jk} \big\|_{\mathcal{S}}^{2q}\Big\}  \le  q! 4(4c^{-1}n^{-1})^q + (2q)!4(4c^{-1}n^{-1})^{2q}.
\end{equation}
Finally it follows from that Lemma~\ref{lemma.moment} of the Supplementary Material that there exists some universal constant $\tilde c > 0 $ such that 
\begin{eqnarray*}
	P\Big \{\big\|\widehat{\Sigma}_{0,jk} - \Sigma_{0,jk} \big\|_{\mathcal{S}} \ge 2\cM_1(\mF_{\bX}) \lambda_0 \eta   \Big \}\le 4 \exp\Big\{- \tilde c n \min(\eta^2, \eta)\Big\}.
\end{eqnarray*}

Using the definition of 
$\|\widehat \bSigma_0 - \bSigma_0\|_{\max} = {\max}_{1 \le j,k \le p} \|\widehat \Sigma_{0,jk} - \Sigma_{0,jk}\|_{\cS}$ and applying the union bound of probability, we obtain that
\begin{eqnarray*}
	P\Big \{ \|\widehat \bSigma_0 - \bSigma_0\|_{\max} \ge 2 \cM_1(\mF_{\bX}) \lambda_0 \eta   \Big \}\le 4 p^2\cdot\exp\Big\{ - \tilde c n \min(\eta^2, \eta)\Big\}.
\end{eqnarray*}

Let $\eta = \rho \sqrt{\log p /n} \leq 1$ and $ \rho^2 \tilde c > 2$, which can be achieved for sufficiently large $n$.  We obtain that
\begin{eqnarray*}
	P\left \{ \|\widehat \bSigma_0 - \bSigma_0\|_{\max} \ge 2 \cM_1(\mF_{\bX}) \lambda_0 \rho \sqrt{\frac{\log p}{n}}   \right \}\le 4 p^{2 - \tilde c\rho^2}.
\end{eqnarray*} 

By similar arguments, we can prove Theorem \ref{thm_F_bound}. See Section~\ref{asec.pf.thm3} of the Supplementary Material for details. The proof is complete. $\square$



\subsection{Proofs of Theorems \ref{thm_lambda} and \ref{thm_principal}}
To simplify our notation, for each $j,k=1,\dots,p,$ we will denote $\Sigma_{0,jk}$ and $\widehat\Sigma_{0,jk}$ by $\Sigma_{jk}$ and $\widehat\Sigma_{jk}, $ respectively, in our subsequent proofs. Let $\delta_{jl} = {\min}_{1 \le k \le l}\{\lambda_{jk} - \lambda_{j(k+1)}\}$ and $\widehat \Delta_{jk} = \widehat \Sigma_{jk} - \Sigma_{jk}$ for $j,k=1, \dots,p$ and $l=1,2 \dots.$ It follows from (4.43) and Lemma~4.3 of \cite{Bbosq1} that
\begin{equation}
\label{eigen.bd}
\underset{l \ge 1}{\sup}~ |\widehat \lambda_{jl} - \lambda_{jl}| \le \|\widehat \Delta_{jj}\|_{\cS}~~\mbox{and}~~\underset{l \ge 1}{\sup} ~\delta_{jl}\|\widehat \phi_{jl} - \phi_{jl}\| \le 2\sqrt{2} \|\widehat \Delta_{jj}\|_{\cS}.
\end{equation}
Moreover, we can express $\widehat \lambda_{jl} - \lambda_{jl}$ and $\widehat \phi_{jl} - \phi_{jl},$ as stated in Lemma~\ref{lemma.lambda.phi} in Section~\ref{sec.pf.lemma.lambda.phi} of the Supplementary Material.

{\bf Proof of Theorem~\ref{thm_lambda}}. By Lemma~\ref{lemma.lambda.phi}, we obtain that 
\begin{eqnarray*}
\frac{\widehat \lambda_{jl} - \lambda_{jl}}{\lambda_{jl}} = \frac{\big\langle \phi_{jl}, \langle \widehat \Delta_{jj}, \phi_{jl} \rangle \big \rangle}{\lambda_{jl}}  + \frac{R_{jl}}{\lambda_{jl}}, ~j =1,\cdots,p, l = 1,\ldots,L.
\end{eqnarray*}
Note that $\lambda_{jl} = \big\langle \phi_{jl}, \langle \Sigma_{jj}, \phi_{jl} \rangle \big \rangle$. It follows from (\ref{thm_Phi_1}) in Theorem~\ref{thm_con_op} that for any $\eta > 0$,
\begin{equation}
\label{con.bd.lambda.ratio}
    P\left\{\left|\frac{\big\langle \phi_{jl}, \langle \widehat \Delta_{jj}, \phi_{jl} \rangle \big \rangle}{\lambda_{jl}}\right| > \cM_1(\mF_{\bX}) \eta\right\} \le 2 \exp\Big\{- c n \min(\eta^2,\eta)\Big\}.
\end{equation}
We next turn to the term $\big|R_{jl}/\lambda_{jl}\big|$. By (\ref{eigen.bd}), Lemma~\ref{lemma.lambda.phi} and Condition~\ref{cond_eigen} with $\delta_{jl} \ge c_0 l^{-\alpha -1}$ and $\lambda_{jl} \ge c_0 \alpha^{-1} l^{-\alpha},$ we have
$$
\left|\frac{R_{jl}}{\lambda_{jl}}\right| \le 4\sqrt{2}c_0^{-2} \alpha l^{2\alpha+1} \|\widehat \Delta_{jj}\|_{\cS}^2. 
$$
It then follows from (\ref{thm_Sigma_compt}) in Theorem~\ref{thm_max_bound} that there exists some constant $\tilde c> 0 $ such that for any $\eta > 0$
\begin{equation}
\label{con.bd.R}
P\left\{ \left|\frac{R_{jl}}{\lambda_{jl}}\right| > 4\sqrt{2}c_0^{-2} \alpha l^{2\alpha+1} \{2 \cM_1(\mF_{\bX}) \lambda_0\eta\}^2\right\} \le 2 \exp\Big\{- \tilde c n \min(\eta^2,\eta)\Big\}.
\end{equation}
Let $c_1 = \min(c, \tilde c)$. It follows from $\rho_1 = 16\sqrt{2} c_0^{-2} \alpha \lambda_0^2$, (\ref{con.bd.lambda.ratio}) and (\ref{con.bd.R}) that 
$$
P\left\{\left |\frac{\widehat \lambda_{jl} - \lambda_{jl}}{\lambda_{jl}}\right|  >  \cM_1(\mF_{\bX}) \eta  + \rho_1 l^{2\alpha + 1} \cM_{1}^2(\mF_{\bX}) \eta^2\right\} \le 4 \exp\Big\{- c_1 n \min(\eta^2,\eta)\Big\}.
$$

Applying the union bound of probability, we obtain that 
$$
P\left\{\underset{1 \leq j\leq p, 1\leq l \leq M}{\max}\left |\frac{\widehat \lambda_{jl} - \lambda_{jl}}{\lambda_{jl}}\right|  >  \cM_1(\mF_{\bX}) \eta  + \rho_1 l^{2\alpha + 1} \cM_{1}^2(\mF_{\bX}) \eta^2\right\} \le 4pM \exp\Big\{- c_1 n \min(\eta^2,\eta)\Big\}.
$$
Let $\eta = \tilde\rho_2 \sqrt{\log (pM)/n} \leq 1$ and $1+\rho_1 M^{2\alpha+1}\cM_{1}(\mF_{\bX})\eta\leq \tilde\rho_1$, both of which can be achieved for sufficiently large $n.$ We obtain that
$$
P\left\{\underset{1 \leq j\leq p, 1\leq l \leq M}{\max}\left |\frac{\widehat \lambda_{jl} - \lambda_{jl}}{\lambda_{jl}}\right|  >  \tilde\rho_1\tilde\rho_2\cM_1(\mF_{\bX}) \sqrt{\frac{\log (pM)}{n}} \right\} \le 4(pM)^{1-c_1\rho_2^2}.
$$
The proof is complete. $\square$

\vspace{1ex}
{\bf Proof of Theorem \ref{thm_principal}.} 
It follows from (\ref{eigen.bd}), Condition~\ref{cond_eigen} with $\delta_{jl} \ge c_0 l^{-\alpha -1}$ and (\ref{thm_Sigma_compt}) in Theorem~\ref{thm_max_bound} that there exists some universal constant $c_1=\min(c,\tilde c)$ such that for any $\eta>0,$ (\ref{thm_con_phi}) holds.

Next we will prove the concentration bound on $\langle \widehat \phi_{jl} - \phi_{jl}, g\rangle.$ It follows from the expansion $g(\cdot) = \sum_{m = 1}^\infty g_{jm} \phi_{jm}(\cdot)$ and (\ref{express.phi}) in Lemma~\ref{lemma.lambda.phi} that
\begin{eqnarray*}
\big \langle \widehat \phi_{jl} - \phi_{jl}, g\big \rangle &=& 
 \sum_{m: m \neq l}(\widehat \lambda_{jl} - \lambda_{jm})^{-1} g_{jm} 
\big \langle \widehat \phi_{jl}, \langle \widehat \Delta_{jj}, \phi_{jm} \rangle \big \rangle + 
g_{jl} \big \langle \widehat\phi_{jl} - \phi_{jl}, \phi_{jl} \big\rangle\\
& = &\sum_{m: m \neq l}\Big \{(\widehat \lambda_{jl} - \lambda_{jm})^{-1} -(\lambda_{jl} - \lambda_{jm})^{-1} \Big\} g_{jm} 
\big \langle \widehat \phi_{jl}, \langle \widehat \Delta_{jj}, \phi_{jm} \rangle \big \rangle \\
 && + \sum_{m: m \neq l}(\lambda_{jl} - \lambda_{jm})^{-1} g_{jm} 
\big \langle \widehat \phi_{jl} - \phi_{jl}, \langle \widehat \Delta_{jj}, \phi_{jm} \rangle \big \rangle \\
&& +  \sum_{m: m \neq l}(\lambda_{jl} - \lambda_{jm})^{-1} g_{jm} 
\big \langle \phi_{jl}, \langle \widehat \Delta_{jj}, \phi_{jm} \rangle \big \rangle 
+ g_{jl} \big \langle \widehat\phi_{jl} - \phi_{jl}, \phi_{jl} \big\rangle\\
& = & I_{1} +I_{2} + I_3 + I_4.
\end{eqnarray*}

Let $\Omega_{d_j} = \{ 2 \|\widehat \Delta_{jj}\|_{\cS} \le \delta_{jd_j}\}.$ It follows from Condition~\ref{cond_eigen} and (\ref{thm_Sigma_compt}) with the choice of $\eta=\{4\cM_1(\mF_{\bX})\lambda_0l^{\alpha+1}\}^{-1}c_0\leq 1$ that 
\begin{equation}
\label{bd.Delta}
    P\big(\Omega_{d_j}^C\big) \leq P\big( \|\widehat \Delta_{jj}\|_{\cS}\geq 2^{-1}c_0 l^{-\alpha-1}\big) \leq 4 \exp\left\{-16^{-1}\tilde c\cM_1^{-2}(\mF_{\bX})\lambda_0^{-2}c_0^2 l^{-2(\alpha + 1)} n\right\}.
\end{equation}
On the event $\Omega_{d_j},$ we can see that ${\sup}_{l \leq d_j}|\widehat\lambda_{jl}-\lambda_{jl}|\leq \lambda_{jd_j}/2,$ which implies $2^{-1}\lambda_{jl} \leq \widehat\lambda_{jl} \leq 2\lambda_{jl}.$ Moreover, $|\widehat \lambda_{jl}-\lambda_{jl}| \leq 2^{-1}|\lambda_{jl}-\lambda_{jm}|$ for $1\leq l \neq m \leq d_j$ and hence $
|\widehat \lambda_{jl} - \lambda_{jm}|\ge 2^{-1} | \lambda_{jl} - \lambda_{jm}|$ for $j=1,\ldots,p.$ By Condition~\ref{cond_eigen}, $|\lambda_{jl} - \lambda_{jm}| \ge c_0 l^{- \alpha -1}$ for $1\leq m \neq l \leq d_j.$ Using the above results, we have 
\begin{eqnarray*}
|I_1|^2 &\le& \big(\widehat \lambda_{jl} - \lambda_{jl}\big)^{2} \sum_{m: m \neq l}(\widehat \lambda_{jl} - \lambda_{jm})^{-2}(\lambda_{jl} - \lambda_{jm})^{-2} g_{jm}^2 \| \widehat \Delta_{jj}\|_{\cS}^2\\
&\le& 4\big(\widehat \lambda_{jl} - \lambda_{jl}\big)^{2}\| \widehat \Delta_{jj}\|_{\cS}^2\sum_{m: m \neq l}(\lambda_{jl} - \lambda_{jm})^{-4} g_{jm}^2\\
&\le& 4 c_0^{-4}\|g^{-jl}\|^2 l^{4(\alpha+1)} \big(\widehat \lambda_{jl} - \lambda_{jl}\big)^{2}  \|\widehat \Delta_{jj}\|_{\cS}^2,
\end{eqnarray*}
where $||g^{-jl}||=(\sum_{m:m \neq l} g_{jm}^2)^{1/2}.$
This together with (\ref{eigen.bd}) implies that, on the event $\Omega_{d_j},$
$$ |I_1| \le 2 c_0^{-2} \|g^{-jl}\| l^{2(\alpha+1)}\big \| \widehat \Delta_{jj}\big\|_{\cS}^2.$$
Similarly, we can show that 
$$
|I_2| \le c_0^{-1} \|g^{-jl}\|  l^{\alpha+1} \|\widehat \phi_{jl} - \phi_{jl}\| \|\widehat \Delta_{jj}\|_{\cS} \le 2\sqrt{2}c_0^{-2} \|g^{-jl}\| l^{2(\alpha+1)}  \|\widehat \Delta_{jj}\|_{\cS}^2.
$$
Moreover, by the result $||\widehat\phi_{jl}-\phi_{jl}||^2=\langle\widehat\phi_{jl}-\phi_{jl},-2\phi_{jl}\rangle$ and (\ref{eigen.bd}) we have
$$|I_4| = 2^{-1}|g_{jl}| \|\widehat \phi_{jl} - \phi_{jl}\|^2 \le 4 c_0^{-2} |g_{jl}| l^{2(\alpha + 1)} \|\widehat \Delta_{jj}\|_{\cS}^2.
$$
Combing the above upper bound results, we have
$$
|I_1| + |I_2| +|I_4| \le (6+2\sqrt{2}) c_0^{-2} \|g\| l^{2(\alpha + 1)} \|\widehat \Delta_{jj}\|_{\cS}^2.
$$

Let $\widetilde{\lambda}_{g} = {\sum}_{m: m \neq l} \lambda_{jm} (\lambda_{jl} - \lambda_{jm})^{-2} g_{jm}^2 \leq c_0^{-2}  l^{2(\alpha + 1)}   \|g^{-jl}\|_{\lambda}^2.$ Then it follows from (\ref{thm_Phi_12}) in Theorem~\ref{thm_con_op} and the fact $\lambda_{jl}+\lambda_{jm} \geq 2\lambda_{jl}^{1/2}\lambda_{jm}^{1/2}$ that
\begin{equation}
\label{I3.con.bd}
P\left\{ \Big|\lambda_{jl}^{-1/2}\widetilde{\lambda}_{g}^{-1/2} I_3\Big | \ge 2 \cM_1(\mF_{\bX}) \lambda_0 \eta\right \} \le 4 \exp\Big\{ - c n \min(\eta^2, \eta)\Big\}.
\end{equation}
Define $\Omega_{1,\eta} = \Big \{ \|\widehat \Delta_{jj}\|_{\cS} \le 2 \cM_1(\mF_{\bX}) \lambda_0 \eta\Big\}$ and
$$
\Omega_{2,\eta} = \Big\{ \big|I_3\big| \le 2 c_0^{-1}  \lambda_0 \lambda_{jl}^{1/2} \|g^{-jl}\| _{\lambda}\cM_1(\mF_{\bX}) l^{\alpha + 1} \eta\Big \}.
$$

Let
$
\rho_{2} = 2 c_0^{-1}  \lambda_0
$
and 
$\rho_{3} = 4(6+2\sqrt{2}) c_0^{-2} \lambda_0^2.$
Under the event $\Omega_{d_j} \cap \Omega_{1,\eta}\cap \Omega_{2,\eta},$ we obtain that 
$$
\left|\big \langle \widehat \phi_{jl} - \phi_{jl}, g\big \rangle \right| \le
\rho_2\|g^{-jl}\|_{\lambda}  \cM_1(\mF_{\bX}) \lambda_{jl}^{1/2}l^{\alpha + 1}  \eta
+ \rho_3 \|g\|\cM_1^2(\mF_{\bX})  l^{2(\alpha+1)}\eta^2.
$$

Let $c_2 = \min(16^{-1}\lambda_0^{-2}c_0^2\tilde c,c_1)$. 
It follows from (\ref{thm_Sigma_compt}) in Theorem~\ref{thm_max_bound} and (\ref{I3.con.bd}) that $$P\big(\Omega_{1,\eta}^C\cup \Omega_{2,\eta}^C\big) 
\le 8 \exp\Big\{ - c_2n \min(\eta^2, \eta)\Big\}.$$ This together with (\ref{bd.Delta}) completes the proof of (\ref{thm_con_phi_g}).
Finally, letting $\eta = \tilde\rho_4 \sqrt{\frac{\log (pM)}{n}}<1,$ and following the similar technique used in the proof of (\ref{bd_lambda_max}), we can obtain (\ref{bd_phi_max}) and (\ref{bd_phig_max}), which completes the proof. $\square$

\subsection{Proof of Theorem \ref{thm_cov}}
For the special case of $(j,k,l,m)$ with $j=k,$ provided that $\widehat \sigma_{jjlm} = \widehat \lambda_{jl}I(l = m)$ and $\sigma_{jjlm} = \lambda_{jl}I(l = m)$ for $j=1,\ldots,p$ and $l,m =1,\ldots,d_j,$ (\ref{con_bd_score1}) follows directly from Theorem~\ref{thm_lambda}.

For general cases of $(j,k,l,m)$ with $j \neq k,$ 
$
\widehat{\sigma}_{jklm} = n^{-1}\sum_{t = 1}^n \widehat{\xi}_{tjl} \widehat{\xi}_{tkm}
$
and  $\sigma_{jklm}= E(\xi_{tjk} \xi_{tlm}).$  
Let $\widehat r_{jl} = \widehat \phi_{jl} - \phi_{jl},$ then error of $\widehat{\sigma}_{jklm} - \sigma_{jklm}$ can be decomposed as
\begin{eqnarray*}
\label{thm_cov_1}
\widehat{\sigma}_{jklm} - \sigma_{jklm} 
& = & \big\langle \widehat r_{jl}, \langle \widehat \Sigma_{jk}, \widehat r_{km}\rangle \big \rangle + 
 \Big(\big\langle \widehat r_{jl}, \langle \widehat \Delta_{jk}, \phi_{km}\rangle \big \rangle + 
 \big\langle \phi_{jl}, \langle \widehat \Delta_{jk}, \widehat r_{km}\rangle \big \rangle \Big)\\
 &&  + \Big(\big\langle \widehat r_{jl}, \langle \Sigma_{jk}, \phi_{km}\rangle \big \rangle  
 + \big\langle \phi_{jl}, \langle \Sigma_{jk}, \widehat r_{km}\rangle \big \rangle \Big) + 
\big\langle \phi_{jl}, \langle \widehat \Delta_{jk}, \phi_{km}\rangle \big \rangle \\
	& = & I_1 + I_2 + I_3 + I_4.
\end{eqnarray*}

Let $\Omega_1 = \Big\{ \|\widehat \Delta_{jk}\|_{\cS} \le \lambda_0\Big\}$ and
$
\Omega_{jk,\eta} = \Big \{ \|\widehat \Delta_{jk}\|_{\cS} \le 2 \cM_1(\mF_{\bX}) \lambda_0 \eta\Big\}.
$ 
On the event $\Omega_1 \cap \Omega_{jj,\eta} \cap \Omega_{kk,\eta} \cap \Omega_{jk,\eta},$ it follows from Condition~\ref{cond_eigen} with $\lambda_{jl} \geq c_0 \alpha^{-1} l^{-\alpha}$, (\ref{eigen.bd}) and Lemma~\ref{lemma.sub.ieq} in Section~\ref{sec.lemma.sub.ieq} of the Supplementary Material  that
\begin{equation}
\label{thm_cov_2}
\begin{split}
\left |\frac{I_1}{\lambda_{jl}^{1/2}\lambda_{km}^{1/2}} \right| 
&\le c_0^{-1}\alpha (lm)^{\alpha/2}\|\widehat r_{jl}\| \big(\|\widehat \Delta_{jk}\| + \|\Sigma_{jk}\|_{\cS}\big)\|\widehat r_{km}\|\\
&\le 64 c_0^{-3} \lambda_0^3 \alpha \cM_{1}^2(\mF_{\bX})(l+m)^{3\alpha + 2} \eta^2, 
\end{split}
\end{equation}
\begin{equation}
\label{thm_cov_3}
\begin{split}
\left|\frac{I_{2}}{\lambda_{jl}^{1/2} \lambda_{km}^{1/2}}\right| 
&\le 2\sqrt{2} c_0^{-2} \alpha (lm)^{\alpha/2}\|\widehat \Delta_{jk}\|_{\cS}\Big( l^{\alpha + 1}\|\widehat \Delta_{jj}\|_{\cS} + m^{\alpha + 1}\|\widehat \Delta_{kk}\|_{\cS}\Big) \\
&\le  16\sqrt{2} c_0^{-2} \alpha \lambda_0^2 \cM_1^2(\mF_{\bX})(l+m)^{2\alpha + 1} \eta^2
\end{split}
\end{equation}
For the term $I_4$, it follows from Theorem~\ref{thm_con_op} and $\lambda_{jl}+\lambda_{jm} \geq 2\lambda_{jl}^{1/2}\lambda_{jm}^{1/2}$ that 
\begin{equation}
\label{thm_cov_4}
P\left\{ \left|\frac{I_4}{\lambda_{jl}^{1/2}\lambda_{km}^{1/2}}\right| \ge 2 \cM_1(\mF_{\bX}) \lambda_0 \eta\right \} \le 4 \exp\Big\{ - c n \min(\eta^2, \eta)\Big\}.
\end{equation}
Finally, we consider the term $I_{3}$. 
Let
$
\rho_{1}^* = 2 c_0^{-1}  \lambda_0^2
$
and 
$\rho_{2}^* = 4(6+2\sqrt{2}) c_0^{-5/2}\alpha^{1/2} \lambda_0^{5/2}.$
By Lemma~\ref{lemma.sub.ieq}, we have that $\|\langle \Sigma_{jk}, \phi_{km}\rangle\| \le \lambda_{km}^{1/2} \lambda_0^{1/2}$
and $\|\langle \phi_{jl},\Sigma_{jk}\rangle\| \le\lambda_{jl}^{1/2} \lambda_0^{1/2}.$
These results together with (\ref{thm_con_phi_g}) in Theorem~\ref{thm_principal} and Condition~\ref{cond_eigen} with $\lambda_{jl} \geq c_0 \alpha^{-1} l^{-\alpha}$ imply that 
\begin{equation}
\label{thm_cov_5}
\left|\frac{I_3}{\lambda_{jl}^{1/2} \lambda_{km}^{1/2}} \right| \le
\rho_1^* \cM_1(\mF_{\bX}) (l+m)^{\alpha + 1}  \eta
+ \rho_2^* \cM_1^2(\mF_{\bX})  (l+m)^{(5\alpha + 4)/2}\eta^2.
\end{equation}
holds with probability greater than $$1-8 \exp\big\{ - c_2n \min(\eta^2, \eta)\big\} - 
4\exp\big\{ - c_2  \cM_1^{-2}(\mF_{\bX}) n l^{-2(\alpha + 1)}\big\}.
$$

Combining the results in (\ref{thm_cov_2})--(\ref{thm_cov_5}) and by (\ref{thm_Sigma_compt}) in Theorem~\ref{thm_max_bound}, we obtain that there exists three positive constants $\rho_4$, $\rho_5$ and $c_3$ such that 
\begin{eqnarray}
&&P\left\{ \left| \frac{\widehat \sigma_{jklm} - \sigma_{jklm}}{\lambda_{jl}^{1/2} \lambda_{km}^{1/2}}\right| \ge \rho_4 \cM_1(\mF_{\bX}) (l+m)^{\alpha + 1}  \eta
+ \rho_5 \cM_1^2(\mF_{\bX})  (l+m)^{3\alpha + 2}\eta^2\right\} \nonumber \\
&& \le 20 \exp\big\{ - c_3 n \min(\eta^2, \eta)\big\} + 8\exp\big\{ - c_3  \cM_1^{-2}(\mF_{\bX}) n l^{-2(\alpha + 1)}\big\},
\end{eqnarray}
where $c_3=\min(4^{-1}\tilde c,c_2),$
\begin{eqnarray}
\label{thm_cov_rho}
\rho_4 = 2 \lambda_0 +2c_0^{-1} \lambda_0^2~\mbox{and}~\rho_5 = 64c_0^{-3} \alpha\lambda_0^3 + 16\sqrt{2} c_0^{-2}\alpha \lambda_0^2 + 4(6+2\sqrt{2}) c_0^{-5/2}\alpha^{1/2}\lambda_0^{5/2}.
\end{eqnarray}

Finally, letting $\eta = \tilde\rho_6 \sqrt{\frac{\log (pM)}{n}}<1$ and following the similar technique used in the proof of (\ref{bd_lambda_max}), we can obtain (\ref{bd_score_max}), which completes the proof. $\square$

\subsection{Proof of Theorem~\ref{thm.vfar}}
Since $\widehat \bB_{j} \in \eR^{pq \times q}$ is the minimizer of (\ref{vfar.crit.theory}), we have 
$$
-\langle\langle \widehat \bY_j, \widehat \bB_{j}\rangle\rangle + 
\frac{1}{2} \langle\langle \widehat \bB_{j}, \widehat \bGamma \widehat\bB_{j}\rangle\rangle 
 + {\gamma_{nj}}  \|\widehat\bB_{j}\|_{1}^{(q)}
 \le -\langle\langle \widehat \bY_j, \bB_{j}\rangle\rangle + 
 \frac{1}{2} \langle\langle \bB_{j}, \widehat \bGamma \widehat\bB_{j}\rangle\rangle  
 + {\gamma_{nj}}  \|\bB_{j}\|_{1}^{(q)}.
$$
Letting $\bDelta_j = \widehat \bB_{j} - \bB_{j}$ and  $S_j^{c}$ be the complement of $S_j$ in the set $\{1, \dots,p\},$ we have
\begin{eqnarray*}
	\frac{1}{2}\langle\langle \bDelta_j, \widehat \bGamma \bDelta_j \rrangle 
	&\leq& \llangle \bDelta_j, \widehat \bY_j - \widehat \bGamma \bB_{j}\rrangle + 
	{\gamma_{nj}}  \Big( \|\bB_{j}\|_{1}^{(q)} - \|\bB_{1j} + \bDelta_j\|_{1}^{(q)}\Big)\\
	&\leq & \llangle \bDelta_j, \widehat \bY_j - \widehat \bGamma \bB_{j}\rrangle + 
	{\gamma_{nj}} \Big( \|\bB_{jS_j}\|_{1}^{(q)}- \|\bB_{jS_j} + \bDelta_{jS_j}\|_{1}^{(q)} - \|\bDelta_{jS_j^c}\|_{1}^{(q)}\Big)\\
	&\leq& \llangle \bDelta_j, \widehat \bY_j - \widehat \bGamma \bB_j\rrangle + 
	{\gamma_{nj}} \Big( \|\bDelta_{jS_j}\|_{1}^{(q)} - \|\bDelta_{jS_j^c}\|_{1}^{(q)}\Big)
\end{eqnarray*}
By Lemma~\ref{lemma.norm.ieq} in Section~\ref{sec.lemma.norm.ieq} of the Supplementary Material, Condition~\ref{cond.fvar.max.error} and the choice of $\gamma_{nj},$ we have 
$$
\big|\llangle \bDelta_j, \widehat \bY_j - \widehat \bGamma \bB_j\rrangle \big| 
\le \|\widehat \bY_j - \widehat \bGamma \bB_j\|_{\max}^{(q)} \|\bDelta_j\|_{1}^{(q)} \leq \frac{\gamma_{nj}}{2} \big( \|\bDelta_{jS_j}\|_{1}^{(q)} + \|\bDelta_{jS_j^c}\|_{1}^{(q)}\big).
$$
Combing the above two results, we have
\begin{eqnarray*}
	0 \leq \frac{1}{2}\langle\langle \bDelta_j, \widehat \bGamma \bDelta_j \rrangle \le 
	\frac{3\gamma_{nj}}{2} \|\bDelta_{jS_j}\|_{1}^{(q)}   - \frac{\gamma_{nj}}{2} \|\bDelta_{jS_j^c}\|_{1}^{(q)},
\end{eqnarray*}
which implies $\|\bDelta_{jS_j^c}\|_{1}^{(q)} \le 3\|\bDelta_{jS_j}\|_{1}^{(q)}$ and therefore $\|\bDelta_j\|_{1}^{(q)} \le 4\|\bDelta_{jS_j}\|_{1}^{(q)} \le 4 \sqrt{s_j} \|\bDelta_{j}\|_{F}.$ 
This result together with Condition~\ref{cond.fvar.RE} and $\tau_2\ge  32\tau_1 q^2 s_j$ implies that \begin{equation}
\langle\langle \bDelta_j, \widehat \bGamma \bDelta_j \rrangle 
\ge \tau_2 \|\bDelta_j\|_F^2 - \tau_1 q^2\big\{\|\bDelta_j\|_{1}^{(q)}\big\}^2 
\ge \big(\tau_2 - 16\tau_1 q^2 s_j\big)\|\bDelta_j\|_F^2 \ge \frac{\tau_2}{2}\|\bDelta_j\|_F^2.
\end{equation}
Therefore,
$$
\frac{\tau_2}{4}\|\bDelta_j\|_F^2 \le \frac{3}{2} \gamma_{nj} \|\bDelta_j\|_{1}^{(q)} \le 6 \gamma_{nj}s_j^{1/2} \|\bDelta_{j}\|_{F},
$$
which implies that
\begin{equation}
\label{err.B}
\|\bDelta_j\|_F \le \frac{24s_j^{1/2}\gamma_{nj}}{\tau_2} \text{ and  } \|\bDelta_j\|_{1}^{(q)} \le \frac{96 s_j\gamma_{nj}}{\tau_2}.    
\end{equation}
as is claimed in Theorem 7.

Next we prove the upper bound of $\widehat\bA - \bA.$ 

For $k\in S_j,$ it follows from  $\bPsi_{jk}=\int\int\bphi_k(v)A_{jk}(u,v)\bpsi_j(u)^{\T}dudv,$ Condition~\ref{cond.fvar.bias} with $A_{jk}(u,v)=\bphi_k(v)^{\T}\ba_{jk}\bphi_j(u) + (\sum_{l,m=1}^{\infty}-\sum_{l,m=1}^{q})a_{jklm}\phi_{jl}(u)\phi_{km}(v)$ and orthonormality of $\{\phi_{jl}(\cdot)\}_{l\geq 1}$ and $\{\phi_{km}(\cdot)\}_{m\geq 1}$ that  $\|\bPsi_{jk}\|_F = ||\ba_{jk}||_F=\big\{\sum_{l,m=1}^q \mu_{jk}^2 (l+m)^{-2\beta-1}\big\}^{1/2}\leq \big\{\mu_{jk}^2\int_{1}^q\int_{1}^q (x+y)^{-2\beta-1}dxdy\big\}^{1/2}=O(\mu_{jk}).$ For $k \in S_j^c,$ we have $\bPsi_{jk}={\bf 0}.$ Hence
\begin{equation}
    \label{bd.Psi}
    \|\bPsi_{j}\|_1^{(q)} = \sum_{k=1}^p ||\bPsi_{jk}||_F =O\big(\sum_{k \in S_j}\mu_{jk}\big) =O(s_j).
\end{equation}

Observe that 
$
\widehat \bPsi_j-\bPsi_j=\widehat\bD^{-1}\widehat\bB_j-\bD^{-1}\bB_j=(\widehat \bD^{-1}-\bD^{-1})\bB_j + \bD^{-1}(\widehat\bB_j-\bB_j) + (\widehat \bD^{-1}-\bD^{-1})(\widehat\bB_j-\bB_j).
$
It follows from the diagonal structure of $\widehat\bD^{-1}$ and $\bD^{-1}$ that
\begin{equation}
\label{err.Psi0}
\begin{split}
\|\widehat \bPsi_{j} - \bPsi_j\|_1^{(q)}\le & \|(\widehat \bD^{-1} - \bD^{-1})\|_{\max} \|\bB_j\|_1^{(q)} + \|\bD^{-1}\|_{\max}\|\widehat \bB_j - \bB_j\|_{1}^{(q)}\\
&+\|(\widehat \bD^{-1} - \bD^{-1})\|_{\max} \|\widehat \bB_j - \bB_j\|_1^{(q)}.
\end{split}
\end{equation}
By Conditions~\ref{cond_eigen}, \ref{cond.fvar.eigen} and the fact $\widehat\bD_k=\text{diag}\big(\widehat\lambda_{k1}^{1/2}, \dots, \widehat\lambda_{kq}^{1/2}\big),$ 	$\bD_k=\text{diag}\big(\lambda_{k1}^{1/2}, \dots, \lambda_{kq}^{1/2}\big),$
we have $\|(\widehat \bD^{-1} - \bD^{-1})\|_{\max} \leq \alpha^{1/2}c_0^{-1/2}q^{\alpha/2}C_{\lambda}\cM(\mF_{\bX}) \sqrt{\frac{\log (pq)}{n}}$ and $\|\bD^{-1}\|_{\max} \le \alpha^{1/2}c_0^{-1/2} q^{\alpha/2}.$ By Condition~\ref{cond.cov.func} and (\ref{bd.Psi}), we have $||\bB||_1^{(q)}\leq ||\bD||_{\max}^{(q)}\|\bPsi_{j}\|_1^{(q)}= O(\lambda_0^{1/2}s_j).$ These results together with (\ref{err.B}) implies that 
\begin{equation}
    \label{err.Psi}
    \|\widehat \bPsi_{j} - \bPsi_j\|_1^{(q)} \leq \frac{96\alpha^{1/2}q^{\alpha/2}s_j\gamma_{nj}}{c_0^{1/2}\tau_2}\Big\{1+o(1)\Big\},
\end{equation}
where the constant comes from the second term in (\ref{err.Psi0}), since the first and third terms are of smaller orders relative to the second term. 

For each $j,k=1,\dots,p,$ note that 
\begin{eqnarray*}	
\widehat \tA_{jk}(u,v) - \tA_{jk}(u,v)	&=& \widehat\bphi_k(v)^\T \widehat\bPsi_{jk} \widehat\bphi_j(u) -\bphi_k(v)^\T \bPsi_{jk} \bphi_j(u) + \tR_{jk}(u,v)\\
&=&\widehat \bphi_k(v)^{\T} \widehat \bPsi_{jk} \left\{\widehat \bphi_j(u) -\bphi_j(u) \right\} +
 \left\{\widehat \bphi_k(v)-\bphi_k(v)\right\}^{\T} \widehat \bPsi_{jk}\bphi_j(u)  \\
	& & +~ \bphi_k(v)^{\T} (\widehat \bPsi_{jk} -\bPsi_{jk})\bphi_j(u) + \tR_{jk}(u,v),
\end{eqnarray*}

We bound the first three terms. By Lemma~\ref{lemma.err.A} in Section~\ref{sec.pf.lemma.err.A} of the Supplementary Material, we have
\begin{eqnarray}
\label{error.Psi}
&&\left\|\widehat \bphi_k(v)^{\T} \widehat \bPsi_{jk} \left\{\widehat \bphi_j(u) -\bphi_j(u) \right\}\right\|_{\cS}\leq
 q^{1/2} \max_{1 \le l \le q}\|\widehat \phi_{jl} - \phi_{jl}\| \|\widehat \bPsi_{jk} \|_F, \nonumber \\
&&\left\|\left\{\widehat \bphi_k(v)-\bphi_k(v)\right\}^{\T} \widehat \bPsi_{jk}\bphi_j(u) \right\|_{\cS} \leq
q^{1/2} \max_{1 \le m \le q}\|\widehat \phi_{km} - \phi_{km}\| \|\widehat \bPsi_{jk} \|_F,\\
&&\left \| \bphi_k(v)^{\T} (\widehat \bPsi_{jk} - \bPsi_{jk})\bphi_j(u)\right\|_{\cS}= \|\widehat \bPsi_{jk} - \bPsi_{jk}\|_F. \nonumber 
\end{eqnarray}

We then bound the fourth term. By $R_{jk}(u,v)=(\sum_{l,m=1}^{q}-\sum_{l,m=1}^{\infty})a_{jklm}\phi_{jl}(u)\phi_{km}(v),$ we have
\begin{eqnarray*}
||R_{jk}||_{\cS}^2 &=& O(1) \Big\|\sum_{l=1}^{q+1}\sum_{m=1}^{\infty}a_{jklm}\phi_{jl}(u)\phi_{km}(v)\Big\|_{\cS}^2\\
&=&O(1) \sum_{l=1}^{q+1}\sum_{m=1}^{\infty}a_{jklm}^2 \leq O(1)\mu_{jk}^2\sum_{l=1}^{q+1}\sum_{m=1}^{\infty}(l+m)^{-2\beta-1} = O(\mu_{jk}^2q^{-2\beta+1}).
\end{eqnarray*}
This together with Condition~\ref{cond.fvar.bias} implies that
\begin{equation}
    \label{err.R}
    \underset{1\leq j \leq p}{\max}\sum_{k=1}^p||R_{jk}||_{\cS}\leq O\big(q^{-\beta+1/2}\underset{1\leq j \leq p}{\max}\sum_{k \in S_j}\mu_{jk}\big)=O\big(sq^{-\beta+1/2}\big).
\end{equation}

It follows from (\ref{bd.Psi}), (\ref{err.Psi}), (\ref{error.Psi}), (\ref{err.R}) and the fact $\|\widehat \bPsi_{j} \|_1^{(q)} \le \|\widehat \bPsi_{j} - \bPsi_j\|_1^{(q)} + \| \bPsi_{j} \|_1^{(q)} =O(s_j)$ that
\begin{eqnarray*}
\|\widehat \bA - \bA\|_\infty &\le& 2q^{1/2} \max_{\underset{1 \le l \le q}{1 \le j \le p}}\|\widehat \phi_{jl} - \phi_{jl}\|\max_{1 \le j \le p}\|\widehat \bPsi_{j} \|_1^{(q)} 
+ \max_{1 \le j \le p}\|\widehat \bPsi_{j}  - \bPsi_j\|_1^{(q)} + ||\bR||_{\infty}.\\
&\leq &  \frac{96\alpha^{1/2}q^{\alpha/2}s\gamma_{n}}{c_0^{1/2}\tau_2}\Big\{1+o(1)\Big\},
\end{eqnarray*}
where the constant comes from $\max_{j}\|\widehat \bPsi_{j}  - \bPsi_j\|_1^{(q)},$ since other terms are of smaller orders of this term. The proof is complete. $\square$
\linespread{1.07}\selectfont
\bibliography{paperbib}
\bibliographystyle{dcu}

\newpage
\linespread{1.39}\selectfont
\begin{center}
	\title{\large \bf Supplementary Material to ``A General Theory for Large-Scale Curve Time Series via Functional Stability Measure "}	
	
	\bigskip
	\author{{Shaojun Guo and Xinghao Qiao}}
\end{center}
\bigskip

\setcounter{page}{1}
\setcounter{section}{1}
\renewcommand{\theequation}{S.\arabic{equation}}
\setcounter{equation}{0}

This supplementary material contains additional technical proofs in Appendix~\ref{ap.tech.proof}, derivations of functional stability measure for the illustrative VFAR(1) example in Appendix~\ref{ap.ill.ex}, some derivations for VFAR models in Appendix~\ref{ap.deriv}, details of the algorithm to fit VFAR models in Appendix~\ref{ap.alg} and additional simulation results in Appendix~\ref{ap.sim}.

\section{Additional proofs of technical details}
\label{ap.tech.proof}

\subsection{Proof of Theorem \ref{thm_F_bound}}
\label{asec.pf.thm3}
It follows from the definition of $
\|\widehat \bSigma_0 - \bSigma_0\|_{F}^2 = \sum_{j, k =1}^p \|\widehat \Sigma_{0,jk} - \Sigma_{0,jk}\|_{\cS}^2,$   
Chebyshev's inequality and (\ref{Sigma_moment}) with $q=1$ that for any 
$\eta > 0,$
\begin{eqnarray*}
	P\left\{\|\widehat \bSigma_0 - \bSigma_0\|_{F}>  2\cM_1(\mF_{\bX}) \lambda_0\eta \right\} & \le&\frac  {1} {\left(2 \cM_1(\mF_{\bX}) \lambda_0 \right)^{2} \eta^2} \sum_{j, k =1}^p \tE\|\widehat \Sigma_{0,jk} - \Sigma_{0,jk}\|_{\cS}^2 \\
	&\le &  \frac{p^2}{\eta^2} (16 \tilde c^{-1}n^{-1} + 128 \tilde c^{-2}n^{-2}) \\
	&=&  \frac{p^2}{\eta^2 n} (16 \tilde c^{-1} + 128 \tilde c^{-2}n^{-1}).
\end{eqnarray*}
By letting $\eta = \tilde\rho \sqrt{p^2/n}$ with $\rho > 0$, we have that 
\begin{eqnarray*}
	P\left\{\|\widehat \bSigma_0 - \bSigma_0\|_{F}>  2 \cM_1(\mF_{\bX}) \lambda_0\tilde\rho\sqrt{\frac{p^2}{n}} \right\} 
	\le \tilde\rho^{-2} (16 \tilde c^{-1} + 128 \tilde c^{-2}n^{-1}).
\end{eqnarray*}
The proof is complete. $\square$

\subsection{Proof of Proposition 1}
\label{pf.prop1}
Let $\bY_{1,t} = \bX_t + \bX_{t+h}$ , $\bSigma_{\bY_{1}, \ell} (u,v) = \cov\{\bY_{1,t}(u), \bY_{1,(t+ \ell)}(v)\}, \ell \in \eZ, (u,v) \in \cU^2$ and $\bXi_{\bY_1,\ell}$ be the operator induced from the kernel $\bSigma_{\bY_1,\ell}.$ Define the spectral density matrix operator of $\bY_{1,t}$ by
$$
\mF_{\bY_1,\theta}(\bPhi) = \frac{1}{2 \pi} \sum_{\ell = -\infty}^{\infty} \bXi_{\bY_1,\ell}(\bPhi) \exp(- i \ell \theta), ~\theta \in [-\pi, \pi], \bPhi \in \cH^p,
$$
Then we can obtain that
$
\mF_{\bY_1,\theta}(\bPhi) = \{2 + \exp(-ih\theta) + \exp( ih\theta)\}\mF_{\bX,\theta}(\bPhi).
$
Similarly, by letting ${\bY}_{2,t}(u) = \bX_t(u) - \bX_{t+h}(u)$, 
$\bSigma_{\bY_2, \ell} (u,v) = \cov\{\bY_{2,t}(u), \bY_{2,(t+\ell)}(v)\},\ell \in \eZ (u,v), \in \cU^2,$ $\bXi_{\bY_2,\ell}$ be the operator induced from the kernel $\bSigma_{\bY_2,\ell}$ and $\mF_{\bY_2,\theta}$ be the spectral density matrix operator of $\bY_{2},$ $\theta \in [-\pi, \pi],$ we have
$
\mF_{\bY_2,\theta}(\bPhi) = \{2 - \exp(-ih\theta) - \exp( ih\theta)\}\mF_{\bX,\theta}(\bPhi), \bPhi \in \cH^p.
$ Note that
$$
4\big\langle \bPhi_1, (\widehat{\bXi}_h - \bXi_h)(\bPhi_1)\big\rangle = \big\langle \bPhi_1, (\widehat{\bXi}_{\bY_1,0} - \bXi_{\bY_1,0})(\bPhi_1)\big\rangle - \big\langle \bPhi_1, (\widehat{\bXi}_{{\bY}_2,0} - \bXi_{{\bY}_2,0})(\bPhi_1)\big\rangle    
$$
and
$
\cM(\mF_{{\bY}_i},\bPhi_1)\le 4 \cM(\mF_{\bX},\bPhi_1) 
$ for $i=1,2.$
Combing these with results in the proof of (\ref{thm_Phi_1}) leads to
\begin{eqnarray*}
	&&P\Big[\Big|\big\langle \bPhi_1, (\widehat{\bXi}_{h} - \bXi_{h})(\bPhi_1)\big\rangle\Big|> 2\cM(\mF_{\bX},\bPhi_1)\eta\Big]\\
	&\leq& \sum_{i=1}^2 P\Big[\Big|\big\langle \bPhi_1, (\widehat{\bXi}_{\bY_i,0} - \bXi_{\bY_i,0})(\bPhi_1)\big\rangle\Big|> \cM(\mF_{\bY_i},\bPhi_1)\eta\Big] \leq 4 \exp\Big\{- c n\min\left( \eta^2,\eta\right)\Big\},
\end{eqnarray*}
for some constant $c> 0.$ This result, together with, $\cM(\mF_{\bX},\bPhi_1) \le \cM_k(\mF_{\bX}) \big\langle \bPhi_1, \bXi_0(\bPhi_1)\big\rangle$ implies (\ref{thm_sig_k1}).

Note that
\begin{eqnarray*}
	4\big\langle \bPhi_1, (\widehat{\bXi}_h - \bXi_h)(\bPhi_2)\big\rangle
	&\leq &　\big\langle \widetilde{\bPhi}_1, (\widehat{\bXi}_h - \bXi_h)(\widetilde{\bPhi}_1)\big\rangle 
	-\big\langle \widetilde{\bPhi}_2, (\widehat{\bXi}_h - \bXi_h)({\widetilde\bPhi}_2)\big\rangle,
\end{eqnarray*}
where $\widetilde{\bPhi}_1 = \bPhi_1 + \bPhi_2,$ $\widetilde{\bPhi}_2 = \bPhi_1 - \bPhi_2$ and
$\cM(\mF_{\bX},\widetilde{\bPhi}_i) \le 2\{\cM(\mF_{\bX},\bPhi_1)+ \cM(\mF_{\bX}, \bPhi_2) \}$ for $i=1, 2.$ 
Combing these with results and the proof of (\ref{thm_sig_k1}) leads to
\begin{eqnarray*}
	&&P\Big[\Big|\big\langle \bPhi_1, (\widehat{\bXi}_h - \bXi_h)(\bPhi_2)\big\rangle\Big|> 2\{\cM(\mF_{\bX},\bPhi_1)+ \cM(\mF_{\bX}, \bPhi_2)\}\eta\Big]\\
	&\leq& \sum_{i=1}^2P\Big[\Big|\big\langle \widetilde\bPhi_i, (\widehat{\bXi}_h - \bXi_h)(\widetilde\bPhi_i)\big\rangle\Big|> 2\cM(\mF_{\bX},\widetilde\bPhi_i)\eta\Big]\le 8 \exp\Big\{- c n\min\left( \eta^2,\eta\right)\Big\}
\end{eqnarray*}
for some constant $c> 0.$
This, together with, $\cM(\mF_{\bX},\bPhi_i) \le \cM_k(\mF_{\bX}) \big\langle \bPhi_i, \bXi_0(\bPhi_i)\big\rangle$ for $i=1, 2,$ implies (\ref{thm_sig_k2}), which completes the proof.  $\square$

\subsection{Proof of Proposition \ref{res.re}}
It is  easy to see that  $\btheta^T \widehat \bGamma \btheta  = \btheta^T \bGamma \btheta + \btheta^T (\widehat \bGamma - \bGamma) \btheta.$ Hence we have
\begin{eqnarray*}
	\btheta^T \widehat \bGamma \btheta	 \ge \btheta^T \bGamma \btheta -  \|\widehat \bGamma - \bGamma\|_{\max} \|\btheta\|_1^2.
\end{eqnarray*}
By Condition \ref{cond_min_bound}, $\lambda_{\min}(\bGamma) \ge \underline{\mu},$ where $\lambda_{\min}(\bGamma)$ denotes the minimum eigenvalue of $\bGamma.$ Together with Lemma~\ref{lemma_Gamma_max} in Section~\ref{sec.lemma_Gamma_max}, this proposition follows. $\square$ 
\subsection{Proof of Proposition~\ref{prop_Error_eigen}}
Note that on the event $\big\{|\widehat \lambda_{jl}-\lambda_{jl}| \le 2^{-1} \lambda_{jl}\big\},$ we have $\widehat\lambda_{jl}\geq \lambda_{jl}/2,$ $\widehat \lambda_{jl}^{-1/2} \leq \sqrt{2}\lambda_{jl}^{-1/2}$ and
$
|\widehat \lambda_{jl}^{-1/2} - \lambda_{jl}^{-1/2}| \le \frac{\widehat\lambda_{jl}^{-1}|\widehat\lambda_{jl}-\lambda_{jl}|\lambda_{jl}^{-1}}{\widehat\lambda_{jl}^{-1/2}+\lambda_{jl}^{-1/2}}\le 2\lambda_{jl}^{-3/2}|\widehat \lambda_{jl} - \lambda_{jl}|,
$ which implies that $\left|\frac{\widehat\lambda_{jl}^{-1/2}-\lambda_{jl}^{-1/2}}{\lambda_{jl}^{-1/2}}\right| \leq 2\left|\frac{\widehat\lambda_{jl}-\lambda_{jl}}{\lambda_{jl}}\right|.$ Then it follows from Theorems~\ref{thm_lambda} and \ref{thm_principal} that there exists positive constants $C_{\lambda}, C_{\phi},$ $c_4$ and $c_5$ such that the first and second deviation bounds in (\ref{fvar.eigen}) respectively hold with probability greater than $1 -c_4(pq)^{-c_5}.$ The proof is complete. $\square$

\subsection{Proof of Proposition~\ref{prop_Error_max}}
Notice that 
\begin{eqnarray}
\label{term1}
\widehat \bY_j - \widehat \bGamma \bB_j
& = & \Big\{(n-1)^{-1}\widehat \bD^{-1} \widehat \bZ^\T \widehat \bV_j - (n-1)^{-1}\bD^{-1} \tE(\bZ^\T \bV_j) \Big\}\\
& & 	+ (n-1)^{-1}\bD^{-1} \tE\Big\{\bZ^\T (\bV_j - \bZ\bD^{-1} \bB_j)\Big\} \nonumber - \big(\widehat \bGamma - \bGamma\big) \bB_j. 
\end{eqnarray}


First, we show the deviation bounds of $\widehat \bD^{-1} (n-1)^{-1}\widehat \bZ^\T \widehat \bV_j - \bD^{-1} \tE((n-1)^{-1}\bZ^\T \bV_j).$  We  decompose this term as 
$\widehat \bD^{-1} \Big\{(n-1)^{-1}\widehat \bZ^\T \widehat \bV_j- \tE((n-1)^{-1}\bZ^\T \bV_j) \Big\}+ (\widehat \bD^{-1} - \bD^{-1}) \tE((n-1)^{-1}\bZ^\T \bV_j).$
It follows from Theorem~\ref{thm_cov} that there exists positive constants $C_1^*,$ $c_4$ and $c_5$ that
\begin{eqnarray}
\label{term2-1}
\sup_{j,k}\left\|\bD_k^{-1} \Big\{(n-1)^{-1}\widehat \bZ_k^\T \widehat \bV_j- \tE\big((n-1)^{-1}\bZ_k^\T \bV_j\big) \Big\}\bD_j^{-1}\right\|_{\max} \le C^*_1 \cM_1(\mF_{\bX}) q^{\alpha + 1} \sqrt{\frac{\log (pq)}{n}},
\end{eqnarray}
with probability great than $1 - c_4 (pq)^{-c_5}.$ Note that $\widehat \bD_k = \mbox{diag}(\widehat \lambda_{k1}^{1/2}, \ldots, \widehat \lambda_{kq}^{1/2})$ and $\bD_k = \mbox{diag}(\lambda_{k1}^{1/2}, \ldots, \lambda_{kq}^{1/2}),$ it follows from Proposition~\ref{prop_Error_eigen} that there exists positive constant $C_2^*,$ such that
\begin{eqnarray}
\label{term2-2}
\left\|\left(\widehat \bD^{-1} - \bD^{-1}\right)\bD\right\|_{\max} \le C_2^* \cM_1(\mF_{\bX}) \sqrt{\frac{\log (pq)}{n}},
\end{eqnarray}
with probability great than $1 - c_4 (pq)^{-c_5}.$ 
By Condition~\ref{cond.cov.func},
we have $\max_j||\bD_j||_F \leq \lambda_0^{1/2}$ and $\|\bD^{-1}\tE\big((n-1)^{-1}\bZ^\T \bV_j\big)\|_{\max}^{(q)} \le q^{1/2}\|\bD^{-1}\tE\big((n-1)^{-1}\bZ^\T \bV_j\big)\bD_j^{-1}\|_{\max}||\bD_j||_F=O(q^{1/2}),$ where the fact that, for $q \times q$ matrix $\bA$ and a diagonal matrix $\bB,$ $\|\bA \bB\|_{F} \le q^{1/2} \|\bA\|_{\max} \|\bB\|_F,$ is used.
These results together with (\ref{term2-1}) and (\ref{term2-2}) imply that there exists $C_3^*$ 
\begin{eqnarray}
\label{term2}
\left\|\widehat \bD^{-1} (n-1)^{-1}\widehat \bZ^\T \widehat \bV_j- \bD^{-1} \tE\big((n-1)^{-1}\bZ^\T \bV_j\big)\right\|_{\max}^{(q)} \le C_3^* \cM_1(\mF_{\bX}) q^{\alpha + 3/2} \sqrt{\frac{\log (pq)}{n}} 
\end{eqnarray}

Second, consider the bias term $
(n-1)^{-1}\bD^{-1} \tE\{\bZ^\T (\bV_j - \bZ\bD^{-1} \bB_j)\}.$ By Section~\ref{ap.vfar.deriv} of the Supplementary Material, $\bR_j$ is a $(n-1) \times q$ matrix  whose row vectors are formed by $\{\br_{tj}, t=2, \dots, n\}$ with $\br_{tj}=(r_{tj1}, \dots, r_{tjq})^{\T}$ and $r_{tjl}=\sum_{k=1}^p \sum_{m=q+1}^{\infty} \langle\phi_{jl}, \langle A_{jk},\phi_{km}\rangle \rangle  \xi_{(t-1)km}$ for $l=1,\dots,q.$  It follows Conditions~\ref{cond.cov.func}, \ref{cond.fvar.bias} and similar arguments in deriving~(\ref{err.R}) and (\ref{term2}) that
there exists some positive constant $C_4^*$ such that
\begin{eqnarray}
\label{term3}
&&\big\|(n-1)^{-1}\bD^{-1} \tE\{\bZ^\T (\bV_j - \bZ\bD^{-1} \bB_j)\}\big\|_{\max}^{(q)}\nonumber\\
&\leq& q^{1/2}\big\|(n-1)^{-1}\bD^{-1} \tE(\bZ^\T\bR_j)\bD_j^{-1}\big\|_{\max} ||\bD_j||_F 
\le C_4^* s_j q^{-\beta + 1 }.
\end{eqnarray}

Third, it follows from Lemma~\ref{lemma_Gamma_max}, Lemma~15 in the Supplementary Material of \cite{qiao2018a} and $||\bB||_1^{(q)}= O\big(\lambda_0^{1/2}s_j\big)$ that there exist some positive constants $C_5^*$ such that 
\begin{eqnarray}
\label{term4}
\left \|\big(\widehat \bGamma - \bGamma\big)\bB_j \right\|_{\max}^{(q)} \leq \left \|\widehat \bGamma - \bGamma \right\|_{\max}^{(q)}||\bB_j||_1^{(q)}\le \cM_1(\mF_{\bX}) s_j q^{\alpha + 2}\sqrt{\frac{\log (pq)}{n}},
\end{eqnarray}
with probability great than $1 - c_4 (pq)^{-c_5}.$ 

Combing results in (\ref{term1}), (\ref{term2}), (\ref{term3}) and  (\ref{term4}) implies that there exist positive constants $C_E, c_4$ and $c_5$ such that
$$
||\widehat \bY_j- \widehat \bGamma\bB_j||_{\max}^{(q)} \leq C_E \cM_1(\mF_{\bX}) s_j \Big \{q^{\alpha +2} \sqrt{\frac{\log (pq)}{n}}  + q^{-\beta + 1}\Big \}, ~ j =1,\ldots, p,
$$
with probability greater than $1 -c_4(pq)^{-c_5}.$ The proof is complete. $\square$

\subsection{Lemma~\ref{lemma_Gamma_max} and its proof}
\label{sec.lemma_Gamma_max}
\begin{lemma}
	\label{lemma_Gamma_max}
	Suppose that Conditions \ref{cond.bd.fsm}-\ref{cond_eigen} hold.  Then there exist some positive constants $C_{\Gamma}$, $c_4$ and $c_5$ such that 
	$$
	\big\|\widehat \bGamma - \bGamma\big\|_{\max} \le C_\Gamma \cM_1(\mF_{\bX}) q^{\alpha +1}\sqrt{\frac{\log (pq)}{n}}
	$$
	with probability greater than $1 - c_4(pq)^{-c_5}$.
\end{lemma}
{\bf Proof.} 
Note that 
$$
\|\widehat \bGamma - \bGamma\|_{\max} = \max_{1 \le j,k \le p, 1 \le l,m \le q} \Big| \widehat{\lambda}_{jl}^{-1/2} \widehat{\lambda}_{km}^{-1/2} \widehat \sigma_{jklm} - \lambda_{jl}^{-1/2} \lambda_{km}^{-1/2}\sigma_{jklm}\Big|.
$$
Let $\widehat{s}_{jklm} = \frac{\widehat \lambda_{jl} \widehat \lambda_{km}}{\lambda_{jl} \lambda_{km}}$ for each $(j,k,l,m)$. Then we have 
\begin{eqnarray*}
	\widehat{\lambda}_{jl}^{-1/2} \widehat{\lambda}_{km}^{-1/2} \widehat \sigma_{jklm} - \lambda_{jl}^{-1/2} \lambda_{km}^{-1/2}\sigma_{jklm} 
	= \widehat{s}_{jklm}^{-1/2}\left(\frac{\widehat \sigma_{jklm} - \sigma_{jklm}}{\lambda_{jl}^{1/2} \lambda_{km}^{1/2}} \right)
	+ \Big(\widehat{s}_{jklm}^{-1/2} - 1\Big)\frac{\sigma_{jklm}}{\lambda_{jl}^{1/2} \lambda_{km}^{1/2}}.
\end{eqnarray*}
Let  
$
\Omega_{\lambda} = \left\{\sup_{1\le j \le p, 1\le l \le q}\left|\frac{\widehat \lambda_{jl} - \lambda_{jl}}{\lambda_{jl}}\right| \le 1/5\right\}.
$
Observe that
$$
\widehat s_{jklm} - 1 = \left(\frac{\widehat \lambda_{jl} - \lambda_{jl}}{\lambda_{jl}} + 1\right)\left( \frac{\widehat \lambda_{km} - \lambda_{km} }{\lambda_{km}}\right) + \frac{\widehat \lambda_{jl} - \lambda_{jl}}{ \lambda_{jl}}.
$$
Then under the event $\Omega_{\lambda},$ we have $|\widehat s_{jklm} - 1| \le 1/2,$ and thus $\widehat{s}_{jklm}^{-1/2} \le \sqrt{2}.$ Moreover, provided that fact that $|(1+x)^{-1/2} - 1| \le x$ if $|x| \le 1/2$, we have 
$$
\Big|\widehat{s}_{jklm}^{-1/2} - 1\Big| \le \frac{6}{5} \left(\left|  \frac{\widehat \lambda_{km} - \lambda_{km} }{\lambda_{km}}\right| + \left| \frac{\widehat \lambda_{jl} - \lambda_{jl} }{\lambda_{jl}}\right|\right).
$$

Under the event $\Omega_{\lambda},$ the above results together with the fact of $\sigma_{jklm}\leq \lambda_{jl}^{1/2}\lambda_{km}^{1/2}$ imply that
$$
\|\widehat \bGamma - \bGamma\|_{\max} \le \sqrt{2}\max_{1 \le j,k \le p, 1 \le l,m \le q}\left|\frac{\widehat \sigma_{jklm} - \sigma_{jklm}}{\lambda_{jl}^{1/2} \lambda_{km}^{1/2}} \right| + \frac{12}{5}\max_{1\le j \le p, 1\le l \le q}\left|\frac{\widehat \lambda_{jl} - \lambda_{jl}}{\lambda_{jl}}\right|.
$$
Then it follows from Theorems~\ref{thm_lambda} and \ref{thm_cov} that there exist some positive constants $C_\Gamma$, $c_4$ and $c_5$ such that 
$$
\|\widehat \bGamma - \bGamma\|_{\max} \le C_\Gamma \cM_1(\mF_{\bX}) q^{\alpha +1}\sqrt{\frac{\log (pq)}{n}}
$$
with probability greater than $1 - c_4(pq)^{-c_5}$. The proof is complete. $\square$

\subsection{Lemma~\ref{lemma.moment} and its proof}
The following lemma  shows how to derive the tail probability through moment conditions. 

\begin{lemma}
	\label{lemma.moment}
	Let $\tX$ be a random variable. If for some constants $ c_1,c_2 > 0$
	$$
	P\left(|\tX| > t \right ) \le c_1 \exp\{- c_2^{-1}\min(t^2, t)\} ~~ \mbox{for any } t >0,
	$$
	then for any integer $q \ge 1$,
	$$
	\tE(\tX^{2q}) \le q!c_1(4c_2)^q + (2q)!c_1(4c_2)^{2q}.
	$$
	Conversely, if for some positive constants $a_1,a_2$,
	$
	\tE(\tX^{2q}) \le q!a_1a_2^q + (2q)!a_1a_2^{2q}, ~ q \ge 1,
	$ then by letting $c^*_2 = 8\max\{4(a_2 + a^2_2), a_2\}$ and $c_1^* = a_1$, we have that
	$$
	P\left(|\tX| > t \right ) \le c_1^* \exp\{- c^{*-1}_2\min(t^2, t)\} ~~ \mbox{for any } t >0.
	$$
\end{lemma}
{\bf Proof.} This lemma can be proved in a similar way to Theorem 2.3 of \cite{BLM2014} and hence the proof is omitted here. In the proof, the following two 
inequalities are used, i.e. for any $c >0$ and $t > 0$,
$$
\frac{1}{2}\min(t^2,t) \le \frac{t^2}{1 + t} \le  \min(t^2,t),
$$
and
$$
\sqrt{\frac{ct}{2}} + \frac{ct}{2} \le \frac{c(t + \sqrt{t^2 + 4t/c})}{2} \le 
\sqrt{ct} + ct.
$$

\subsection{Lemma~\ref{lemma.lambda.phi} and its proof}
\label{sec.pf.lemma.lambda.phi}
\begin{lemma}
	\label{lemma.lambda.phi}
	For each $j = 1,\ldots,p$ and $l=1,\dots,$ the term of $\widehat \lambda_{jl} - \lambda_{jl}$ can be expressed as
	\begin{eqnarray}
	\label{express.lambda}
	\widehat \lambda_{jl} - \lambda_{jl} = \big\langle \phi_{jl}, \langle \widehat \Delta_{jj}, \phi_{jl} \rangle \big \rangle  + R_{jl},
	\end{eqnarray}
	where $|R_{jl}| \le 2 \|\widehat \phi_{jl} - \phi_{jl}\|\|\widehat \Delta_{jj} \|_{\cS}$.
	Furthermore, if $\inf_{m: m \neq l} |\widehat \lambda_{jl} - \lambda_{jm}| > 0$, then
	\begin{eqnarray}
	\label{express.phi}
	\widehat \phi_{jl} - \phi_{jl} = \sum_{m: m \neq l}(\widehat \lambda_{jl} - \lambda_{jm})^{-1} \phi_{jm} 
	\big \langle \widehat \phi_{jl}, \langle \widehat \Delta_{jj}, \phi_{jm} \rangle \big \rangle 
	+\phi_{jl}\big \langle \widehat\phi_{jl} - \phi_{jl}, \phi_{jl} \big\rangle.
	\end{eqnarray}    
\end{lemma}
{\bf Proof}. This lemma follows directly from Lemma 5.1 of \cite{hall2007} and hence the proof is omitted here. $\square$

\subsection{Lemma~\ref{lemma.sub.ieq} and its proof}
\label{sec.lemma.sub.ieq}
\begin{lemma}
	\label{lemma.sub.ieq}
	For a $p$ by $p$ covariance matrix function, $\bSigma=\big(\Sigma_{jk}(\cdot,\cdot)\big)_{1 \leq j,k \leq p} \in \eS^{p \times p},$ we have
	$$||\Sigma_{jk}||_{\cS} \leq \max_{1\leq j \leq p} \int_{\cU} \Sigma_{jj}(u,u)du=\lambda_0$$
	and $$||\langle \Sigma_{jk}, \phi_{km}\rangle|| \leq \lambda_{km}^{1/2}\lambda_0^{1/2} ~\text{ for } m \geq 1$$
\end{lemma}
{\bf Proof}.
Since every principal 2 by 2 submatrix is positive semi-definite, we have
$
\{\Sigma_{jk}(u,v)\}^2 \leq \Sigma_{jj}(u,u) \Sigma_{kk}(v,v),
$ which implies that $$||\Sigma_{jk}||_{\cS} \leq \sqrt{\int_{\cU}\Sigma_{jj}(u,u)du\int_{\cU}\Sigma_{kk}(v,v)dv} \leq \max_{1\leq j\leq p}\int_{\cU}\Sigma_{jj}(u,u)du=\lambda_0.$$

Provided that $\Sigma_{jk}(u,v)=\sum_{l,m=1}^{\infty}\cov(\xi_{tjl},\xi_{tkm})\phi_{jl}(u)\phi_{km}(v),$ we have
$$||\langle \Sigma_{jk}, \phi_{km}\rangle||^2=\Big(\sum_{l=1}^{\infty}E(\xi_{tjl}\xi_{tkm})\phi_{jl}(u)\Big)^2du \leq \sum_{l=1}^{\infty}E(\xi_{tjl}^2)E(\xi_{tkm}^2)\leq \int_{\cU}\Sigma_{jj}(u,u)du\lambda_{km}\leq \lambda_0\lambda_{km},$$ which completes the proof. $\square.$

\subsection{Lemma~\ref{lemma.err.A} and its proof}
\label{sec.pf.lemma.err.A}
\begin{lemma}
	\label{lemma.err.A}
	For each $j,k=1,\dots,p,$ let $\{\phi_{jl}(\cdot)\}_{1\leq l \leq q}$ and $\{\widehat\phi_{jl}(\cdot)\}_{1 \leq l \leq q}$ correspond to true and estimated eigenfunctions, respectively, and $\widehat\psi_{jklm}$ be the estimate of $\psi_{jklm}$ for $l,m=1,\dots,q.$ Then we have
	\begin{eqnarray*}
		\left\|\sum_{l=1}^q\sum_{m=1}^q\widehat \phi_{km}(\cdot)\widehat \psi_{jklm} \left\{\widehat \phi_{jl}(\cdot) -\phi_{jl}(\cdot) \right\}\right\|_{\cS}^2 & \leq&
		\sum_{l=1}^q\|\widehat \phi_{jl} - \phi_{jl}\|^2 \sum_{l=1}^q\sum_{m=1}^q\widehat \psi_{jklm}^2,\\
		\left \| \sum_{l=1}^q\sum_{m=1}^q\phi_{km}(\cdot)(\widehat \psi_{jklm} - \psi_{jklm})\phi_{jl}(\cdot)\right\|_{\cS}^2&=&  \sum_{l=1}^q\sum_{m=1}^q(\widehat \psi_{jklm} - \psi_{jklm})^2.
	\end{eqnarray*}
\end{lemma}
{\bf Proof}. We prove the first result
\begin{eqnarray*}
	&&  \left\|\sum_{l =1}^q \sum_{m=1}^q\widehat \phi_{km}\widehat \psi_{jklm} (\widehat \bphi_{jl} -\bphi_{jl}) \right\|_\cS^2\\
	& = &\sum_{m =1 }^q \left \|\sum_{l=1}^q \widehat \psi_{jklm} (\widehat \bphi_{jl} -\bphi_{jl}) \right\|^2
	\leq \sum_{m =1 }^q \sum_{l=1}^q \widehat \psi_{jklm}^2 \sum_{l = 1}^q \|\widehat \bphi_{jl} -\bphi_{jl}\|^2,
\end{eqnarray*}
where the first equality from the orthonormality of $\{\widehat \phi_{km}(\cdot)\}_{1 \leq m \leq q}$ and the second inequality comes from Cauchy-Schwarz inequality.
By the orthonormality of of $\{\phi_{km}(\cdot)\}_{1 \leq m \leq q}$  and $\{\phi_{jl}(\cdot)\}_{1 \leq l \leq q},$ we can prove the second result 
\begin{eqnarray*}
	&&  \left\|\sum_{l =1}^q \sum_{m=1}^q \phi_{km}(\widehat \psi_{jklm}-\psi_{jklm})\bphi_{jl} \right\|_\cS^2\\
	& = &\sum_{m =1 }^q \left \|\sum_{l=1}^q (\widehat \psi_{jklm} - \psi_{jklm})\bphi_{jl} \right\|^2
	= \sum_{m =1 }^q \sum_{l=1}^q (\widehat \psi_{jklm}-\psi_{jklm})^2, 
\end{eqnarray*}
which completes the proof. $\square$

\subsection{Lemma~\ref{lemma.norm.ieq} and its proof}
\label{sec.lemma.norm.ieq}
\begin{lemma}
	\label{lemma.norm.ieq}
	Let $\bA,\bB \in \eR^{pq\times q}$ with $j$-th blocks given by $\bA_j,\bB_j \in \eR^{q\times q},$ respectively. We have
	\begin{equation}
	\label{norm.eq1}
	\llangle\bA,\bB\rrangle \leq ||\bB||^{(q)}_{\max} ||\bA||_1^{(q)}.
	\end{equation}
\end{lemma}
{\bf Proof}. By the definition and Cauchy-Schwarz inequality
\begin{eqnarray*}
	\llangle\bA,\bB\rrangle &=& \sum_{j=1}^p \llangle\bA_j,\bB_j\rrangle\\
	&\leq& \sum_{j=1}^p\llangle\bA_j,\bA_j\rrangle^{1/2} \llangle\bB_j,\bB_j\rrangle^{1/2}\\
	&\leq& \underset{j}{\max}||\bB_j||_F\sum_{j=1}^p||\bA_j||_F=||\bB||^{(q)}_{\max}||\bA||_1^{(q)},
\end{eqnarray*}
which completes the proof. $\square$

\section{An illustrative example}
\label{ap.ill.ex}
In the following, for any $\bA, \bB \in \eS^{p \times p}$ and $\bx \in \cH^p,$
write $  \bA \bB,$ $\bA\bx$ and $\bx^{\T}\bA$ for 
$$
\int_{\cU} \bA(u,v') \bB(v',v)dv' ~\mbox{,}~ \int_{\cU} \bA(u,v) \bx(v)dv  ~\mbox{and}~ \int_{\cU} \bx(u)^{\T}\bA(u,v) du,
$$ respectively. For a $p$ by $p$ matrix, $\bC,$ we denote its maximum eigenvalue, spectral radius and operator norm by $\lambda_{\max}(\bC),$ $\rho(\bC)=|\lambda_{\max}(\bC)|$ and $||\bC||=\sqrt{\lambda_{\max}(\bC^{\T}\bC)},$ respectively.


Let $\bx_t=\big(x_{t1},x_{t2}\big)^{\T},$
$\bpsi=\text{diag}\big(\psi_1,\psi_2\big),$
$\bC=\left[\begin{array}{cc}a& b\\0 & a\end{array}\right]$
and $\be_t=\big(e_{t1},e_{t2}\big)^{\T},$ then the model VFAR(1) in
(\ref{vfar1.ex}) and (\ref{vfar1.ex.A}) can be rewritten as
$$
\bpsi(u)\bx_t = \int_{\cU} \bpsi(u)\bC\bpsi(v)\bpsi(v)\bx_{t-1}dv+ \bpsi(u)\be_t,
$$ which leads to a VAR(1) model
\begin{equation}
\label{var1.ex}
\bx_t = \bC \bx_{t-1} + \be_t.
\end{equation}

Provided that $\bA(u,v)=\bpsi(u)\bC\bpsi(v)$ and $||\bC||=\sqrt{\lambda_{\max}(\bC^{\T}\bC)}=\lambda_1$ with
$\bC^{\T}\bC\by=\lambda_1^2\by$ for $||\by||=1,$
it is easy to see that
$$\bA^{\T}\bA=\int \bA^{\T}(u,v')\bA(v',v)dv'=
\int\bpsi(u)\bC^{\T}\bpsi(v')\bpsi(v')\bC\bpsi(v)dv'=\bpsi(u)\bC^{\T}\bC\bpsi(v)$$
and $$\int(\bA^{\T}\bA)(u,v)(\bpsi(v)\by)dv=\int\bpsi(u)\bC^{\T}\bC\bpsi(v)\bpsi(v)\by dv=\bpsi(u)\bC^{\T}\bC\by=\lambda_1^2\bpsi(u)\by.$$ Hence $||\bA||_{\cL}=\sqrt{\lambda_{\max}(\bA^{\T}\bA)}=||\bC||=\lambda_1.$ The left side of Figure~\ref{ill.example} plots $||\bA||_{\cL}$ vs $b$ for different values of $a \in (0,1).$

Let $(\omega_j,\bv_j), j=1,2,$ be the eigen-pairs of $\bC$ satisfying $\bC\bv_j=\omega_j\bv_j.$ Then $$\int_{\cU} \bpsi(u)\bC\bpsi(v)\bpsi(v)\bv_j dv=\int_{\cU}\bA(u,v)\bpsi(v)\bv_j dv=\omega_j\bpsi(u)\bv_j.$$ Hence $\bA$ and $\bC$ share the same eigenvalues, which are $\omega_1=\omega_2=a.$ When $\rho(\bA)=\rho(\bC)=|a|<1,$ (\ref{vfar1.ex}) and (\ref{var1.ex}) correspond to stationary VFAR(1) and VAR(1) models, respectively.

For the VFAR(1) model in (\ref{vfar1.ex}), the spectral density matrix function and the covariance matrix function of $\{\bX_t\}_{t \in \eZ}$ are
\begin{equation}
\label{density.vfar1}
f_{\bX,\theta} = \frac{1}{2\pi}\left(\bSigma_0  + \sum_{h = 1}^\infty \Big\{ \bSigma_0 (\bA^{\T})^{h} \exp(-i h \theta)  +  \bA^{h}\bSigma_0  \exp(i h \theta) \Big\}\right)
\end{equation}
and
\begin{equation}
\label{Sigma0.vfar1}
\bSigma_0=\sigma^2\sum_{h=1}^\infty \bA^{h}(\bA^{h})^{\T},
\end{equation}
respectively.
For the VAR(1) model in (\ref{var1.ex}), the spectral density matrix and the covariance matrix of $\{\bx_t\}_{t \in \eZ}$ are
\begin{equation}
\label{density.var1}
f_{\bx,\btheta}=
\frac{1}{2\pi}\left(\bS_0  + \sum_{h = 1}^\infty \Big\{ \bS_0 (\bC^{\T})^{h} \exp(-i h \theta)  +  \bC^{h}\bS_0  \exp(i h \theta) \Big\}\right)
\end{equation}
and
\begin{equation}
\label{Sigma0.var1}
\bS_0=\sigma^2\sum_{h=1}^\infty \bC^{h}(\bC^{h})^{\T},
\end{equation}
respectively. Noting that $\bA^{h}\bSigma_0=\int_{\cU}\bpsi(u)\bC^{h}\bpsi(v')\bpsi(v')\bS_0\bpsi(v)dv'=\bpsi(u)\bC^{h}\bS_0\bpsi(v)$ and applying similar techniques, we can obtain that $f_{\bX,\theta}=\bpsi(u)f_{\bx,\theta}\bpsi(v)$ and $\bSigma_0=\bpsi(u)\bS_0\bpsi(v).$ The functional stability measure of $\{\bX_t\}_{t \in \eZ}$ under (\ref{vfar1.ex}) is
$$	2 \pi \cdot \underset{\theta\in [-\pi, \pi], \bPhi \in \cH_0^2}{\text{ess}\sup} \frac{\bPhi^{\T}f_{\bX,\theta}\bPhi}{\bPhi^{\T}\bSigma_0\bPhi}= 2 \pi \cdot \underset{\theta\in [-\pi, \pi], \bPhi \in \cH_0^2}{\text{ess}\sup} \frac{(\bpsi\bPhi)^{\T}f_{\bx,\theta}\bpsi\bPhi}{(\bpsi\bPhi)^{\T}\bS_0\bpsi\bPhi},$$
where $\bpsi\bPhi \in {\eR}^2$ with $(\bpsi\bPhi)_j=\langle\bphi_j,\bPhi\rangle, j=1,2.$ Hence the functional stability measure of $\{\bX_t\}_{t \in \eZ}$ under (\ref{vfar1.ex}) is the same as that of $\{\bx_t\}_{t \in \eZ}$ under (\ref{var1.ex}), i.e. the essential supremum of the maximal eigenvalue of $2\pi\bS_0^{-1/2}f_{\bx,\theta}\bS_0^{-1/2}$ over $\theta \in [-\pi,\pi].$ Some calculations yield $f_{\bx,\theta}$ and $\bS_0$ as follows.

By (\ref{density.var1}), we have
\begin{eqnarray*}
	f_{\bx,\btheta}&=&
	\frac{1}{2\pi}\left(\bS_0+ \bS_0\left[\begin{array}{cc} \sum_{h=1}^{\infty}a^{h}\exp(-ih\theta)& 0\\
		\sum_{h=1}^{\infty} h a^{h-1}\exp(-ih\theta)b &\sum_{h=1}^{\infty}a^{h}\exp(-ih\theta)
	\end{array}\right]\right.\\
	&& \left.+ \left[\begin{array}{cc} \sum_{h=1}^{\infty}a^{h}\exp(ih\theta)& \sum_{h=1}^{\infty} h a^{h-1}\exp(ih\theta)b\\
		0 &\sum_{h=1}^{\infty}a^{h}\exp(ih\theta)
	\end{array}\right]\bS_0\right)\\
	&=&\frac{1}{2\pi}\left(\bS_0+\bS_0\left[\begin{array}{cc} \frac{\alpha\exp(-i\theta)}{1-a\exp(-i\theta)}& 0\\
		\frac{b\exp(-i\theta)}{(1-a\exp(-i\theta))^2}&\frac{a\exp(-i\theta)}{1-a\exp(-i\theta)} 
	\end{array}\right]
	+
	\left[\begin{array}{cc} \frac{a\exp(i\theta)}{1-a\exp(i\theta)}& \frac{b\exp(i\theta)}{(1-a\exp(i\theta))^2}\\
		0&\frac{a\exp(i\theta)}{1-a\exp(i\theta)} 
	\end{array}\right]\bS_0
	\right).
\end{eqnarray*}

By (\ref{Sigma0.var1}), we have
\begin{eqnarray*}
	\bS_0&=&\left[\begin{array}{cc}
		\sum_{h=0}^{\infty} a^{2h} + \sum_{h=0}^{\infty}h^2a^{2h-2}b^2& \sum_{h=0}^{\infty} ha^{2h-1}b\\ 
		\sum_{h=0}^{\infty} ha^{2h-1}b & \sum_{h=0}^{\infty} a^{2h}\end{array}\right]\\
	&=&\left[\begin{array}{cc}\frac{1}{1-a^2} + \frac{(a^2+1)b^2}{(1-a^2)^3}& \frac{ab}{(1-a^2)^2}\\\frac{ab}{(1-a^2)^2} & \frac{1}{1-a^2}\end{array}\right].
\end{eqnarray*}
The right side of Figure~\ref{ill.example} plots functional stability measures of $\{\bX_t\}_{t \in \eZ}$ vs $b$ for different values of $a \in (0,1).$

\section{Derivations for VFAR models}
\label{ap.deriv}
\subsection{Matrix representation of a VFAR($L$) model in (\ref{vfar3})}
\label{ap.vfar.deriv}
Note that the VFAR(L) model in (\ref{vfar1}) can be equivalently represented as
\begin{equation}
\label{vfar2}
X_{tj}(u) = \sum_{h = 1}^L \sum_{k=1}^p
\langle A_{hjk}(u,\cdot), X_{(t-h)k}(\cdot)\rangle + \varepsilon_{tj}(u), ~~t = L+1,\ldots,n,j=1\dots,p.
\end{equation}
It then follows from the Karhunen-Lo\'eve expansion that (\ref{vfar2}) can be rewritten as
\begin{eqnarray*}
	\sum_{l=1}^{\infty} \xi_{tjl}\phi_{jl}(u) &=& \sum_{h=1}^{L} \sum_{k=1}^p \sum_{m=1}^{\infty} \langle A_{hjk}(u,\cdot), \phi_{km}(\cdot)\rangle \xi_{(t-h)km} + \varepsilon_{tj}(u).
\end{eqnarray*}
This, together with orthonormality of $\{\phi_{jm}(\cdot)\}_{m \geq 1},$ implies that
$$	\xi_{tjl} = \sum_{h=1}^{L} \sum_{k=1}^p \sum_{m=1}^{q_k} \langle\phi_{jl},\langle A_{hjk},\phi_{km}\rangle \rangle\ \xi_{(t-h)km}+ r_{tjl} + \epsilon_{tjl},
$$
where $r_{tjl}=\sum_{h=1}^{L} \sum_{k=1}^p \sum_{m=q_k+1}^{\infty} \langle\phi_{jl}, \langle A_{hjk},\phi_{km}\rangle \rangle  \xi_{(t-h)km}$ and 
$\epsilon_{tjl}=\langle \phi_{jl},\varepsilon_{tj}\rangle$ for $l=1,\dots, q_j,$ represent the approximation and random errors, respectively.
Let $\br_{tj}=(r_{tj1}, \dots, r_{tjq_j})^{\T}$ and  $\beps_{tj}=(\epsilon_{tj1}, \dots, \epsilon_{tjq_j})^{\T}.$ Let $\bR_j, \bE_j$ be $(n-L) \times q_j$ matrices whose row vectors are formed by $\{\br_{tj}, t=L+1, \dots, n\}$ and $\{\beps_{tj}, t=L+1, \dots, n\}$ respectively. Then (\ref{vfar2}) can be represented in the matrix form of (\ref{vfar3})

\subsection{VFAR(1) representation of a VFAR($L$) model}
\label{ap.vfar1.rep}
We can represent a $p$-dimensional VFAR($L$) model in (\ref{vfar1}) as a $pL$-dimensional VFAR(1) model in the form of $\widetilde\bX_t(u)=\int_{\cU}\widetilde \bA_{1}(u,v)\widetilde\bX_{t-1}(v)dv+ \widetilde\bvarepsilon_{t-1}(u), u \in \cU$
with $\widetilde \bX_{t}=\left[\begin{array}{c}\bX_{t}\\\bX_{t-1}\\\vdots\\\bX_{t-L+1}\end{array}\right]\in {\cH}^{pL},$ $\widetilde \bA_1=\left[\begin{array}{ccccc}\bA_{1} & \bA_{2} & \cdots & \bA_{L-1} & \bA_{L}\\ \bI_{p} & {\bf 0} & \cdots & {\bf 0} & {\bf 0}\\{\bf 0} & \bI_{p} & \cdots & {\bf 0} & {\bf 0}\\\vdots & \vdots & \ddots & \vdots & \vdots\\{\bf 0} & {\bf 0} & \cdots & \bI_{p} & {\bf 0}\end{array}\right] \in {\eS}^{pL \times pL},$ $\widetilde \bvarepsilon_{t}=\left[\begin{array}{c}\bvarepsilon_{t}\\
\bvarepsilon_{t-1}\\\vdots\\ \bvarepsilon_{t-L+1}\end{array}\right]\in {\cH}^{pL}$ and $\bI_p=\big(C_{jk}(\cdot,\cdot)\big)_{1 \leq j,k \leq p}\in {\eS}^{p \times p}$ with $C_{jk}(u,v)=I(j=k)I(u=v).$ See a similar VAR(1) representation of a VAR($L$) model in \cite{basu2015a}. 

\subsection{VFAR(1) representation of the simulation example}
\label{ap.vfar1.sim}
Noting that $\btheta_t=\bB\btheta_{t-1}+\bfeta_t,$ we have $\btheta_{tj}=\sum_{k=1}^p\bB_{jk}\btheta_{(t-1)k}+\bfeta_{tj}$ for $j=1,\dots,p.$ Multiplying both sides by $\bs(u)^{\T}$ and applying $\int_{\cU}\bs(v)\bs(v)^{\T}dv=\bI_5,$ we obtain that $\bs(u)^{\T}\btheta_{tj}=\int_{\cU}\sum_{k=1}^p\bs(u)^{\T}\bB_{jk}\bs(v)\bs(v)^{\T}\btheta_{(t-1)k}dv+ \bs(u)^{\T}\bfeta_{tj}.$ Letting $A_{jk}(u,v)=\bs(u)^{\T}\bB_{jk}\bs(v),$ $X_{tj}(u)=\bs(u)^{\T}\btheta_{tj}$ and $\varepsilon_{tj}(u)=\bs(u)^{\T}\bfeta_{tj},$ we have $X_{tj}(u)=\sum_{k=1}^p\langle A_{jk}(u,\cdot)X_{(t-1)k}(\cdot)\rangle + \varepsilon_{tj}(u).$

\section{Algorithms in fitting VFAR models}
\label{ap.alg}

\subsection{Selection of tuning parameters}
\label{ap.select.tune}
To fit the proposed VFAR model, we need choose values for three tuning parameters, $q_j$ (the number of selected principal components for $j=1,\dots,p$), $\eta_j$ (the smoothing parameter when performing regularized FPCA, as described in Section~\ref{ap.rfpca}) and $\gamma_{nj}$ (the regularization parameter in (\ref{vfar.crit}) to control the block sparsity level in $\{\widehat\bPsi_{hjk}: h=1,\dots,L, k=1,\dots,p\}$).

We adopt a $K$-fold cross-validated method to choose $(q_j,\eta_j)$ for each $j.$ Specifically, let $W_{tjs}$ be observed values of $X_{tj}(u_s)$ at $u_1, \dots, u_T.$ We randomly divide the set $\{1, \dots, n\}$ into $K$ groups, $\cD_1, \dots, \cD_K$ of approximately equal size, with the first group treated as a validation set. Implementing regularized FPCA on the remaining $K-1$ groups, we obtain estimated mean function $\widehat\mu_{jl}^{(-1)}(u),$ FPC scores $\widehat\xi_{tjl,\eta_j}^{(-1)}$ and eigenfunctions $\widehat \phi_{jl}^{(-1)}(u;\eta_j)$ for $l=1, \dots, q_j.$ The predicted curve for the $t$-th sample in group one can be computed by $\widehat W_{tjs}^{(1)}=\widehat\mu_{jl}^{(-1)}(u_s)+ \sum_{l=1}^{q_j}\widehat\xi_{tjl,\eta_j}^{(-1)}\widehat \phi_{jl}^{(-1)}(u_s;\eta_j)$. This procedure is repeated $K$ times. Finally, we choose $q_j$ and $\eta_j$ as the values that minimize the mean cross-validated error, $$\text{CV}(q_j,\eta_j)=(KT)^{-1}\sum_{k=1}^K \sum_{s=1}^T\sum_{t \in \cD_k} (W_{tjs}-\widehat W_{tjs}^{(k)})^2.$$

There are several possible methods one could adopt to select the regularization parameter $\gamma_{nj}$ for each $j.$ Popular approaches include AIC, BIC and cross-validation. While the third one is computational intensive, we take an approach motivated by the information criterion for sparse additive models \cite[]{voorman2014}. Our proposed AIC and BIC criterion take the form of 
\begin{equation}
\label{aic_vfar}
\text{AIC}_j(\gamma_{nj}) = n \text{log} \{\text{RSS}_j(\gamma_{nj})\} + 2 \text{df}_j(\gamma_{nj})		\end{equation}
and
\begin{equation}
\label{bic_vfar}
\text{BIC}_j(\gamma_{nj}) = n \text{log} \{\text{RSS}_j(\gamma_{nj})\} + \log (n) \text{df}_j(\gamma_{nj}),
\end{equation}
respectively, where $\text{RSS}_j (\gamma_{nj}) = \Big\|\widehat\bV_{0j} - \sum_{h = 1 }^L \sum_{k = 1 }^p \widehat\bV_{hk}\widehat \bPhi_{hjk}^{(\gamma_{nj})}  \Big\|_F^2$ is the residual sum of squares from minimizing (\ref{vfar.crit}) with the regularization parameter $\gamma_{nj}$ and $\text{df}_j(\gamma_{nj})$ is the corresponding value of degrees of freedom. For non-functional data or $q_j=1,$ we can approximate the degrees of freedom by the number of non-zero parameters in $\{\widehat\bPhi_{hjk}:h=1,\dots,L,k=1,\dots,p\}.$ In the functional setting with $q_j>1$, we use
\begin{equation}
\label{df_j}
\text{df}_j(\gamma_{nj}) = \sum_{h = 1 }^L \sum_{k = 1 }^p \left\{ I\Big((h,k): ||\widehat \bPhi_{hjk}^{(\gamma_{nj})}||_F \neq 0\Big) +  (q_{j} q_{k} -1) \frac{||\widehat\bV_{hk}\widehat \bPhi_{hjk}^{(\gamma_{nj})}||_F^2}{||\widehat\bV_{hk}\widehat \bPhi_{hjk}^{(\gamma_{nj})}||_F^2 + \gamma_{nj}} \right\}.
\end{equation}
We seek the value of $\gamma_{nj}$ that minimizes $\text{AIC}_j$ or $\text{BIC}_j.$

\subsection{Regularized FPCA}
\label{ap.rfpca}
In this section, we drop subscripts $j$ for simplicity of notation. Suppose we observe $\bX(\cdot)=(X_1(\cdot),\dots, X_n(\cdot))^{\T}$ on $\cU,$ our goal is to find the first $q$ regularized principal component functions $\{ \phi_{l}(\cdot), l=1, \dots,q\}.$ We obtain the $l$-th leading principal component $\phi_l(\cdot)$ through a smoothing approach, which maximizes the following penalized sample variance [(9.1) in \cite{Bramsay1}]
\begin{equation}
\label{pen.fpca}
\text{PEN}_{\eta}(\phi_l)=\frac{\var(\langle\phi_l, X_i\rangle)}{||\phi_l||^2 + \eta ||\phi''_l||^2},
\end{equation}
subject to $||\phi_l||=1$ and $\langle\phi_l,\phi_{l'}\rangle+\eta\langle\phi_l'',\phi_{l'}''\rangle=0, l'=1,\dots,l-1,$ where $\eta \geq 0$ is a smoothing parameter to control the roughness of $\phi_l(\cdot).$

Suppose that $\bX(u)=\bdelta^{\T}\bbb(u)$ and $\phi_l(u)=\bzeta_l^{\T}\bbb(u)$ where $\bbb(\cdot)$ is a $G$-dimensional B-spline basis function, $\bdelta \in {\eR}^{n\times G}$
and $\bzeta_l \in {\eR}^G$ are the basis coefficients for $\bX(\cdot)$ and $\phi_l(\cdot),$ respectively. Let $\bJ=\int\bbb(u)\bbb(u)^{\T}du,$ $\bU=\bJ\bdelta^{\T}\bdelta\bJ$ and $\bQ=\int\bbb''(u)\bbb''(u)^{\T}du,$ (\ref{pen.fpca}) is equivalent to maximizing
\begin{equation}
\label{pen.fpca1}
\text{PEN}_{\eta}(\phi_l)=\frac{\bzeta_l^{\T}\bU\bzeta_l}{\bzeta_l^{\T}(\bJ+\eta\bQ)\bzeta_l},
\end{equation}
subject to $\bzeta_l^{\T}\bJ\bzeta_l=1$ and $\bzeta_l^{\T}(\bJ+\eta \bQ)\bzeta_{l'}=0, l'=1,\dots, l-1.$ By {\it singular value decomposition} (SVD), we obtain eigen-pairs, $(\bS_1,\bP_1)$ and $(\bS_2,\bP_2)$ such that $\bJ+\eta\bQ=\bP_1\bS_1^{-2}\bP_1^{\T}$ and $\bS_1\bP_1^{\T}\bU\bP_1\bS_1=\bP_2\bS_2^{-2}\bP_2^{\T}.$ Then (\ref{pen.fpca1}) becomes 
$\text{PEN}_{\eta}(\phi_l)=\frac{\bx_l^{\T}\bP_2^{\T}\bS_2\bP_2\bx_l}{\bx_l^{\T}\bx_l},$ where $\bx_l=\bS_1^{-1}\bP_1^{\T}\bzeta_l.$
This suggests us to perform SVD on $\bP_2^{\T}\bS_2\bP_2,$ where we can obtain $\widehat\bx_l,$
$\widehat\bzeta_l=\bP_1\bS_1\widehat\bx_l$ and $\widehat\phi_l(u)=\widehat\bzeta_l^{\T}\bbb(u)/\big(\widehat\bzeta_l^{\T}\bJ\widehat\bzeta_l\big)^{1/2},l=1,\dots, q.$ In practice, we can set $G$ to a pre-specified large enough value, and implement the cross-validation procedure described in Section~\ref{ap.select.tune} to select $q$ and $\eta.$

\subsection{Block FISTA algorithm to solve (\ref{vfar.crit})}
\label{ap.alg.vfar}

The optimization problem in (\ref{vfar.crit}) can be reformulated as follows. 
\begin{eqnarray}
\label{generalproblem}
\min_{\bX \in {\eR}^{r \times q_j}} g(\bX), ~~g(\bX) = f(\bX) + \gamma_{nj} \sum_{k = 1}^{pL} \|\bX_k\|_F,
\end{eqnarray}
where $f(\bX) = 2^{-1}\text{trace}\left\{(\bY - \bB \bX)^{\T} (\bY - \bB \bX)\right\},$ $r = \sum_{h=1}^L\sum_{k=1}^p q_k,$ $\bY \in \eR^{(n-L)\times q_j},$ $\bB \in \eR^{(n-L) \times r},$ and $\bX = (\bX_1^{\T},\ldots,\bX_{pL}^{\T})^{\T} \in \eR^{r \times q_j}$ with $\bX_k \in \eR^{q_k \times q_j}$ for $k=1,\dots,p.$ 
(\ref{generalproblem}) is a convex problem including the smooth part for $\bX_k,$ i.e. $f(\bX)$ and the non-smooth part for $\bX_k$, i.e.$\gamma_{nj} \sum_{k = 1}^{pL} \|\bX_k\|_F$. To solve the minimization problem in (\ref{generalproblem}), we adopt a block version of {\it fast iterative shrinkage-thresholding algorithm} (FISTA) \cite[]{beck2009} combined with a restarting technique \cite[]{dono2015}, namely block FISTA. 

The basic idea behind our proposed block FISTA is summarized as follows. Let $\nabla f(\bX)$ be the gradient of $f(\bX)$ at $\bX.$ We start with an initial value $\bX^{(0)}.$ At the $(m+1)$-th iteration we first try to solve a regularized sub-problem
\begin{eqnarray}
\label{subprob-1}
\min_{\bX \in \eR^{r \times q_j}} \text{trace}\left\{(\nabla f(\bX^{(m)})^{\T} (\bX -  \bX^{(m)})\right\} 
+ (2C)^{-1}\left \|\bX - \bX^{(m)} \right \|_F^2
+ \gamma_{nj} \sum_{k = 1}^{pL} \|\bX_k\|_F,
\end{eqnarray}
where $\bX^{(m)}$ is the $m$-th iterate and $C>0$ is a small constant controlling the stepsize at $(m+1)$-th step. The second term in (\ref{subprob-1}) can be interpreted as a quadratic regularization, which restricts the updated iterate not to be very far from $\bX^{(m)}$. 
The analytical solution to (\ref{subprob-1}) takes the form of
\begin{eqnarray}
\label{thresholdsolution}
\widetilde{\bX}^{(m+1)} = \big(\widetilde{\bX}_k^{(m+1)}\big)~~ \mbox{with}~~\widetilde{\bX}^{(m+1)}_k =  \left( 1 - \gamma_{nj} C{\|\bZ_k^{(m)}\|_F^{-1}}\right)_{+} \bZ_k^{(m)},k=1, \dots, pL,
\end{eqnarray}
where $\bZ^{(m)} = \bX^{(m)} - C\nabla f(\bX^{(m)})=\big(\big(\bZ_1^{(m)}\big)^{\T}, \dots, \big(\bZ_{pL}^{(m)}\big)^{\T}\big)^{\T}$ and $x_+ = \max(0,x).$ (See also (3.a) and (3.b) of Algorithm~\ref{alg1}).


We then take block FISTA \cite[]{beck2009} by adding  an extrapolation step in the algorithm (see also (3.c) and (3.d) of Algorithm~\ref{alg1}): 
$$
\bX^{(m+1)} = \widetilde{\bX}^{(m+1)} + \omega^{(m+1)}(\widetilde{\bX}^{(m+1)} - \widetilde{\bX}^{(m)}),
$$
where the weight $\omega^{(m+1)}$ is specified in Algorithm~\ref{alg1}.  Finally, at the end of each iteration, we evaluate the generalized gradient at $\bX^{(m+1)}$ by computing the sign of 
$$
\text{trace}\left\{(\bX^{(m)} - \widetilde{\bX}^{(m+1)})^{\T}(\widetilde{\bX}^{(m+1)} - \widetilde{\bX}^{(m)})\right\},
$$
which can be thought of a proxy of $\text{trace}\left\{(\nabla g(\bX^{(m)}))^{\T} (\widetilde{\bX}^{(m+1)} -  \widetilde{\bX}^{(m)})\right\}$.
For a positive sign, i.e. the objective function is increasing at $\widetilde{\bX}^{(m+1)}$, we then restart our accelerated algorithm by setting $\bX^{(m+1)} = \bX^{(m)}$ and $ \omega^{(m+1)} = \omega^{(1)}$ \cite[]{dono2015}. This step can guarantee that the objective function $g$ decreases over each iteration. We iterative the above steps until convergence. We summarize the restarting-based block FISTA in Algorithm~\ref{alg1}. In practice, one issue is how to choose the stepsize parameter $C$. In general, the proposed scheme is guaranteed to converge when $C < \big(\lambda_{\max}(\bB^{\T}\bB)\big)^{-1}.$ 
Here we choose $C = 0.9 \big(\lambda_{\max}(\bB^{\T}\bB)\big)^{-1}$, which turns to work well in empirical studies. Alternatively, $C$ can be selected through a line search and one simple backtracking rule.

\begin{center}
	\begin{algorithm}[!t] \caption{\label{alg1}\textbf{Block FISTA for solving (\ref{generalproblem})}}
		\begin{enumerate}
			\item[1.] Input: $C = 0.9 \big(\lambda_{\max}(\bB^{\T}\bB)\big)^{-1}$, $\theta_0 = 1$, $\bX^{(0)} = (\bX_1^{(0)\T},\ldots,\bX_{pL}^{(0)\T})^{\T}= {\bf 0},$
			$\bZ^{(0)} = (\bZ_1^{(0)\T},\ldots,\bZ_{pL}^{(0)\T})^{\T}= {\bf 0},$ $\widetilde{\bX}^{(0)} = (\widetilde{\bX}_1^{(0)\T},\ldots,\widetilde{\bX}_{pL}^{(0)\T})^{\T}= {\bf 0}.$
			\item[2.] For $m = 0, 1, \ldots$ do 	
			\begin{itemize}		
				\item[] (3.a)~ $\bZ^{(m)} = \bX^{(m)} - C\nabla f(\bX^{(m)}),$
				\item[] (3.b)~ $\widetilde{\bX}^{(m+1)}_k =  \left( 1 - \gamma_{nj} C{\|\bZ_k^{(m)}\|_F^{-1}}\right)_+ \bZ_k^{(m)},k=1, \dots, pL,$
				\item[] (3.c)~ $\theta_{m+1} = \big(1+ \sqrt{1+ 4 \theta_m^2}\big)/2,$
				\item[] (3.d)~ $\bX^{(m+1)} = \widetilde{\bX}^{(m+1)} + \frac{\theta_m - 1}{\theta_{m+1}}\big(\widetilde{\bX}^{(m+1)} - \widetilde{\bX}^{(m)}\big),$		
				\item[] (3.e)~ If $\text{trace}\left\{(\bX^{(m)} - \widetilde{\bX}^{(m+1)})^{\T}(\widetilde{\bX}^{(m+1)} - \widetilde{\bX}^{(m)})\right\} > 0$, set 
				$$
				\bX^{(m+1)} = \bX^{(m)}, \theta_{m+1} = 1.
				$$
			\end{itemize}
			\item[] end do until convergence. 
			\item[3.] Output: the final estimator $\bX^{(m+1)}.$
		\end{enumerate}
	\end{algorithm}		
\end{center}

\section{Additional simulation results}
\label{ap.sim}
\begin{figure*}
	\centering
	\includegraphics[width=0.75\textwidth]{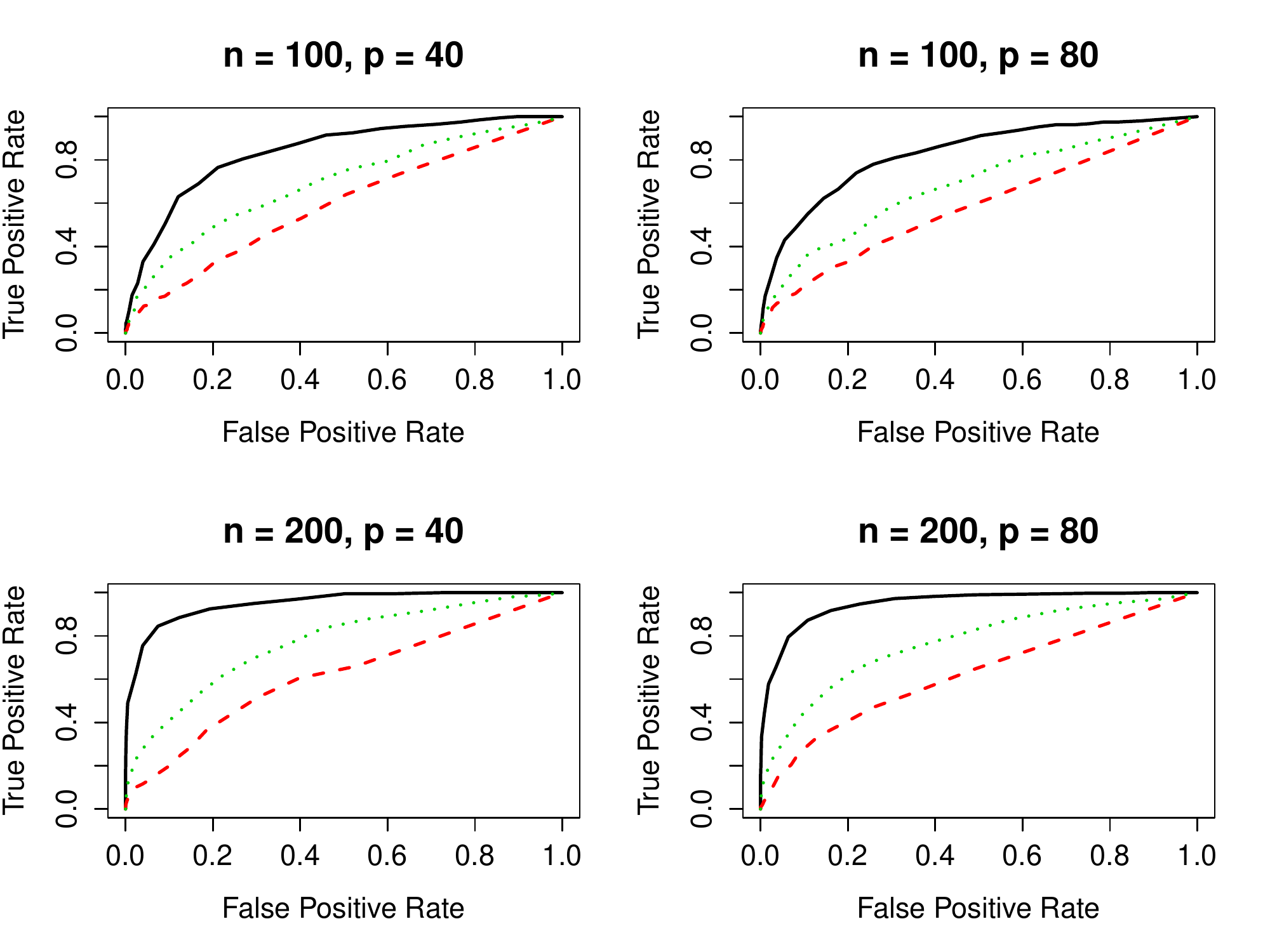}
	\caption{\label{roc.curve.sparse}{Comparisons of median estimated ROC curves over 100 simulation runs. $\ell_1/\ell_2$-$\LS_{\text{a}}$  (black solid), $\ell_1/\ell_2$-$\LS_2$  (green dotted) and $\ell_1$-$\LS_1$  (red dashed) for Model~(i).}}
\end{figure*}

\begin{figure*}
	\centering
	\includegraphics[width=0.75\textwidth]{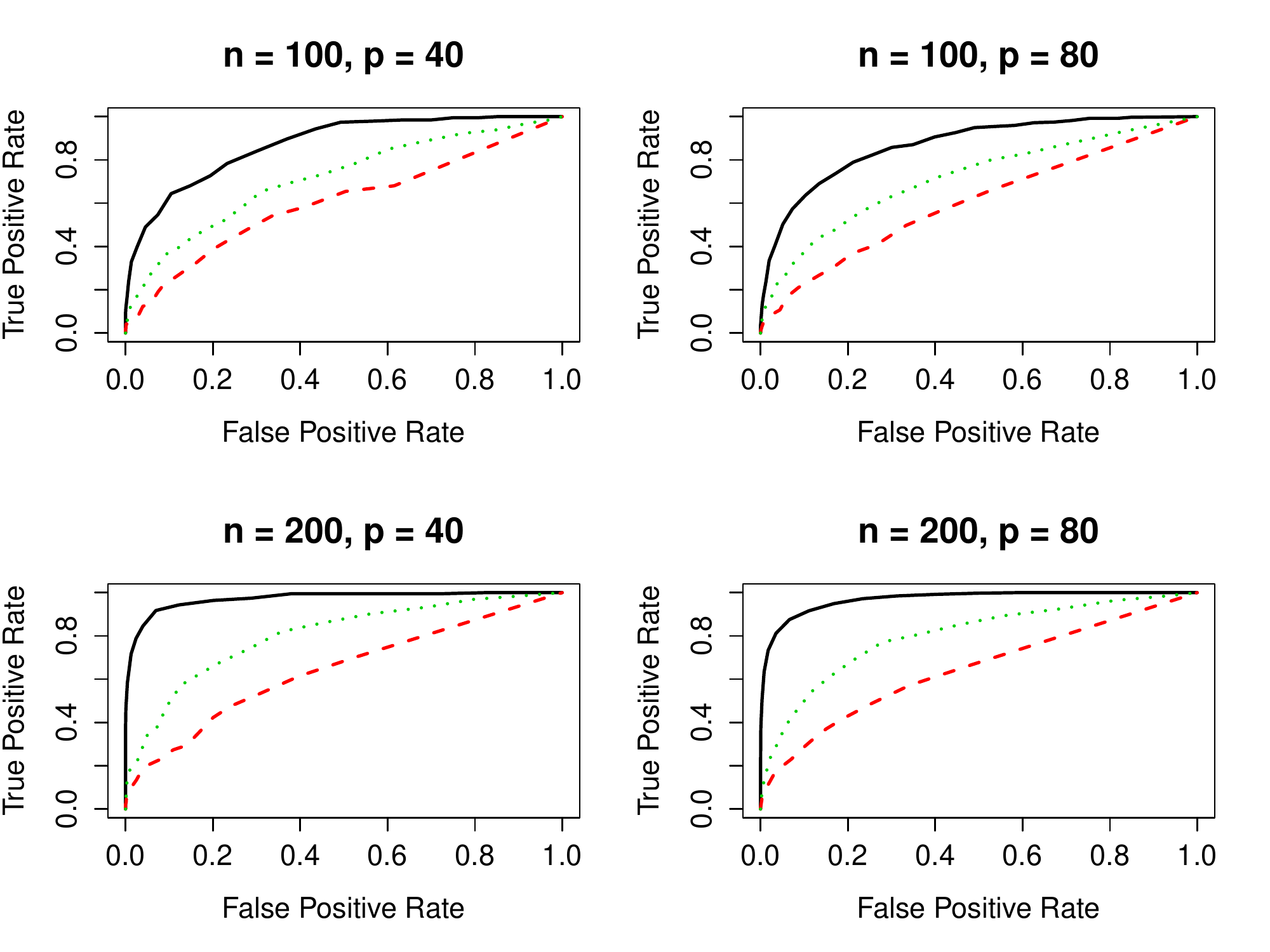}
	\caption{\label{roc.curve.band}{Comparisons of median estimated ROC curves over 100 simulation runs. $\ell_1/\ell_2$-$\LS_{\text{a}}$  (black solid), $\ell_1/\ell_2$-$\LS_2$  (green dotted) and $\ell_1$-$\LS_1$  (red dashed) for Model~(ii).}}
\end{figure*}

Figures~\ref{roc.curve.sparse} and \ref{roc.curve.band} plot the median best ROC curves (we rank ROC curves by the corresponding AUROCs) over the 100 stimulation runs in Models~(i) and (ii), respectively. Again we see that $\ell_1/\ell_2$-$\text{LS}_2,$ which explains partial curve information, although performing better than $\ell_1$-$\text{LS}_1$ is substantially outperformed by $\ell_1/\ell_2$-$\text{LS}_{\text{a}}$ in terms of model selection consistency.

\end{document}